\documentclass[a4paper,10pt,leqno]{amsart}
\title{Euler characteristics of categories and homotopy colimits}
\author{Thomas M. Fiore}
\author{Wolfgang L\"uck}
\author{Roman Sauer}
\address{Thomas M. Fiore \\ Department of Mathematics and Statistics\\
University of Michigan-Dearborn \\ 4901 Evergreen Road \\ Dearborn,
MI 48128 \\ U.S.A.} \email{tmfiore@umd.umich.edu}
\urladdr{http://www-personal.umd.umich.edu/~tmfiore/}
\address{Wolfgang L\"uck \\ Mathematisches Institut der Universit\"at Bonn\\
                Endenicher Allee 60\\
                53115 Bonn \\ Germany}
\email{wolfgang.lueck@him.uni-bonn.de}
\urladdr{http://www.him.uni-bonn.de/lueck}
\address{Roman Sauer \\ Fakult\"at f\"ur Mathematik\\Universit\"at Regensburg\\
Universit\"atsstr. 31\\93053 Regensburg \\Germany}
\email{Roman.Sauer@mathematik.uni-regensburg.de}
\urladdr{http://www.mathematik.uni-regensburg.de/sauer/}

\usepackage{color}

\usepackage{amsmath}
\usepackage{amsthm}
\usepackage{amssymb}
\usepackage{amscd}
\usepackage{amsfonts}  
\usepackage{hyperref}
\usepackage[all]{xy} \SelectTips{cm}{}
\usepackage{graphicx}    
\usepackage{multicol}    

\DeclareMathAlphabet\EuR{U}{eur}{m}{n}
\SetMathAlphabet\EuR{bold}{U}{eur}{b}{n}

\makeindex             

\theoremstyle{plain}
\newtheorem{theorem}{Theorem}[section]
\newtheorem{lemma}[theorem]{Lemma}

\newtheorem{corollary}[theorem]{Corollary}

\theoremstyle{definition}
\newtheorem{definition}[theorem]{Definition}
\newtheorem{construction}[theorem]{Construction}
\newtheorem{example}[theorem]{Example}

\newtheorem{remark}[theorem]{Remark}
\newtheorem{notation}[theorem]{Notation}

{\catcode`@=11\global\let\c@equation=\c@theorem}

\newcommand{\comsquare}[8]                   
{\begin{CD}
#1 @>#2>> #3\\
@V{#4}VV @V{#5}VV\\
#6 @>#7>> #8
\end{CD}
}

\newcommand{\xycomsquareminus}[8]                   
{\xymatrix
{#1 \ar[r]^-{#2} \ar[d]_{#4} &
#3 \ar[d]^{#5}  \\
#6\ar[r]^-{#7} &
#8
}
}

\newcommand{\calb}{{\mathcal B}}
\newcommand{\calc}{{\mathcal C}}
\newcommand{\cald}{{\mathcal D}}

\newcommand{\calg}{{\mathcal G}}
\newcommand{\calh}{{\mathcal H}}
\newcommand{\cali}{{\mathcal I}}
\newcommand{\calj}{{\mathcal J}}

\newcommand{\caln}{{\mathcal N}}

\newcommand{\calq}{{\mathcal Q}}

\newcommand{\calx}{{\mathcal X}}
\newcommand{\caly}{{\mathcal Y}}

\newcommand{\IC}{{\mathbb C}}

\newcommand{\IN}{{\mathbb N}}

\newcommand{\IQ}{{\mathbb Q}}
\newcommand{\IR}{{\mathbb R}}

\newcommand{\IZ}{{\mathbb Z}}

\newcommand{\curs}{\EuR}
\newcommand{\ABELIANGROUPS}{\curs{ABELIAN}\text{-}\curs{GROUPS}}

\newcommand{\CHAINCOMPLEXES}{\curs{CHCOM}}

\newcommand{\CATS}{\curs{CAT}}

\newcommand{\GROUPS}{\curs{GROUPS}}

\newcommand{\MOD}{\curs{MOD}}
\newcommand{\Or}{\curs{Or}}

\newcommand{\SETS}{\curs{SETS}}
\newcommand{\SSETS}{\curs{SSET}}
\newcommand{\SPACES}{\curs{SPACES}}

\newcommand{\aut}{\operatorname{aut}}
\newcommand{\Aut}{\operatorname{Aut}}
\newcommand{\barcon}{\operatorname{bar}}

\newcommand{\colim}{\operatorname{colim}}

\newcommand{\End}{\operatorname{End}}
\newcommand{\ev}{\operatorname{ev}}

\newcommand{\Hom}{\operatorname{Hom}}
\newcommand{\hocolim}{\operatorname{hocolim}}

\newcommand{\id}{\operatorname{id}}

\newcommand{\ind}{\operatorname{ind}}

\newcommand{\iso}{\operatorname{iso}}

\newcommand{\mor}{\operatorname{mor}}

\newcommand{\ob}{\operatorname{ob}}
\newcommand{\op}{\operatorname{op}}

\newcommand{\pr}{\operatorname{pr}}

\newcommand{\res}{\operatorname{res}}

\newcommand{\rk}{\operatorname{rk}}

\newcommand{\Split}{\operatorname{Split}}

\newcommand{\pt}{\{\bullet\}}

\newcommand{\uor}[1]{\underline{\Or}(#1)}
\newcommand{\higherlim}[3]{{\setbox1=\hbox{\rm lim}
        \setbox2=\hbox to \wd1{\leftarrowfill} \ht2=0pt \dp2=-1pt
        \mathop{\vtop{\baselineskip=5pt\box1\box2}}
        _{#1}}^{#2}#3}

\newcommand{\version}[1]                       
{\begin{center} last edited on #1\\
last compiled on \today\\
name of tex-file: \jobname
\end{center}
}

\newcounter{commentcounter}

\begin{document}

\typeout{----------------------------  eulcathocolim.tex  ----------------------------}

\maketitle


\typeout{-----------------------  Abstract  ------------------------}

\begin{abstract} In a previous article, we introduced notions of finiteness
obstruction, Euler characteristic, and $L^2$-Euler characteristic
for wide classes of categories. In this sequel, we prove the
compatibility of those notions with homotopy colimits of
$\cali$-indexed categories where $\cali$ is any small category
admitting a finite $\cali$-$CW$-model for its $\cali$-classifying
space. Special cases of our Homotopy Colimit Formula include
formulas for products, homotopy pushouts, homotopy orbits, and
transport groupoids. We also apply our formulas to Haefliger
complexes of groups, which extend Bass--Serre graphs of groups to
higher dimensions. In particular, we obtain necessary conditions for
developability of a finite complex of groups from an action of a
finite group on a finite category without loops.
  \\[2mm]


\end{abstract}


 \typeout{--------------------   Section 0: Introduction --------------------------}

\setcounter{section}{-1}
\section{Introduction and Statement of Results}

In our previous paper
\cite{FioreLueckSauerFinObsAndEulCharOfCats(2009)}, we presented a
unified conceptual framework for Euler characteristics of categories
in terms of finiteness obstructions and projective class groups.
Many excellent properties of our invariants stem from the
homological origins of our approach: the theory of modules over
categories and the dimension theory of modules over von Neumann
algebras provide us with an array of tools and techniques. In the
present paper, we additionally draw upon the homotopy theory of
diagrams to prove the compatibility of our invariants with homotopy
colimits.

If $\calc \colon \cali \to \CATS$ is a diagram of categories (or
more generally a pseudo functor into the 2-category of small
categories), then our invariants of the homotopy colimit can be
computed in terms of the invariants of the vertex categories
$\calc(i)$. In particular, our Homotopy Colimit Formula,
Theorem~\ref{the:homotopy_colimit_formula}, states
\begin{equation} \label{equ:intro_hocolim}
\chi\bigl(\hocolim_\cali \calc;R\bigr) = \sum_{n \geq 0} (-1)^n
\cdot \sum_{\lambda \in \Lambda_n} \chi(\calc(i_\lambda);R)
\end{equation}
under certain hypotheses. The set $\Lambda_n$ indexes the
$\cali$-$n$-cells of a finite $\cali$-$CW$-model $E\cali$
for the $\cali$-classifying space of $\cali$, that is, we have a functor $E\cali
\colon \cali^{\op} \to \SPACES$ which is inductively built by gluing
finitely many cells of the form $\mor_\cali(-,i_\lambda)\times D^n$
for $\lambda \in \Lambda_n$, and moreover $E\cali(i)\simeq *$ for all
objects $i$ of $\cali$. Similar formulas hold for the finiteness
obstruction, the functorial Euler characteristic, the functorial
$L^2$-Euler characteristic, and the $L^2$-Euler characteristic.

Motivation for such a formula is provided by the classical
Inclusion-Exclusion Principle: if $A$, $B$, and $A \cap B$ are
finite simplicial complexes, then one has
\[ \chi(A\cup B)=\chi(A)
+\chi(B)-\chi(A\cap B).
\]
However, one cannot expect the Euler characteristic to be compatible
with pushouts, even in the simplest cases. The pushout in $\CATS$ of
the discrete categories \[\{\ast\} \leftarrow \{y,z\} \rightarrow
\{\ast'\}\] is a point, but $\chi(\text{point}) \neq 1+1-2$. On the
other hand, their \emph{homotopy} pushout in $\CATS$ is the category
whose objects and nontrivial morphisms are pictured below.
\[
\xymatrix{y \ar[r] \ar[d] & \ast' \\ \ast & z \ar[l] \ar[u]}
\]
The classifying space of this category has the homotopy type of $S^1$,
so that
\[
\chi(\text{homotopy
pushout})=\chi(\{\ast\})+\chi(\{\ast'\})-\chi(\{y,z\})
\]
is true. In fact, the formula for homotopy pushouts is a special
case of \eqref{equ:intro_hocolim}: the category $\cali=\{1
\leftarrow 0 \rightarrow 2\}$ admits a finite model with
$\Lambda_0=\{1,2\}$ and $\Lambda_1=\{0\}$, as constructed
in~Example~\ref{exa:finite_model_for_pushout}.
See~Example~\ref{exa:homotopy_pushout} for the homotopy pushout
formulas of the other invariants.

The Homotopy Colimit Formula in
Theorem~\ref{the:homotopy_colimit_formula} has many applications
beyond homotopy pushouts. Other special cases are formulas for Euler
characteristics of products, homotopy orbits, and transport
groupoids. Our formulas also have ramifications for the
developability of Haefliger's \emph{complexes of groups} in
geometric group theory. If a group $G$ acts on an
$M_\kappa$-polyhedral complex by isometries preserving cell
structure, and if each $g \in G$ fixes each cell pointwise that $g$
fixes setwise, then the quotient space is also an
$M_\kappa$-polyhedral complex, see Bridson--Haefliger
\cite[page~534]{Bridson-Haefliger(1999)}. Let us call the quotient
$M_\kappa$-polyhedral complex $Q$. To each face $\overline{\sigma}$
of $Q$, one can assign the stabilizer $G_\sigma$ of a chosen
representative cell $\sigma$. This assignment, along with the
various conjugated inclusions of groups obtained from face
inclusions, is called the \emph{complex of groups associated to the
group action}. It is a pseudo functor from the poset of faces of $Q$
into groups. In the finite case, the Euler characteristic and
$L^2$-Euler characteristic of the homotopy colimit can be computed
in terms of the original complex and the order of the group. We
prove this
in~Theorem~\ref{thm:Euler_characteristic_of_hocolim_of_quotient_complex}.
Homotopy colimits of complexes of groups play a special role in
Haefliger's theory, see the discussion after
Definition~\ref{def:complex_of_groups}.

In
Section~\ref{sec:finiteness_obstruction_and_Euler_characteristics},
we review the notions and results
from~\cite{FioreLueckSauerFinObsAndEulCharOfCats(2009)} that we need
in this sequel. Explanations of the finiteness obstruction, the
functorial Euler characteristic, the Euler characteristic, the
functorial $L^2$-Euler characteristic, the $L^2$-Euler
characteristic, and the necessary theorems are all contained in
Section~\ref{sec:finiteness_obstruction_and_Euler_characteristics}
in order to make the present paper self-contained.
Section~\ref{sec:spaces_over_a_category} is dedicated to an
assumption in the Homotopy Colimit Formula, namely the requirement
that a finite $\cali$-$CW$-model exists for the $\cali$-classifying
space of $\cali$. We recall the notion of $\cali$-$CW$-complex,
present various examples, and prove that finite models are preserved
under equivalences of categories. Homotopy colimits of diagrams of
categories are recalled in
Section~\ref{sec:homotopy_colimits_of_categories}. The homotopy
colimit construction in $\CATS$ is the same as the Grothendieck
construction, or the category of elements. Thomason proved that the
homotopy colimit construction has the expected properties.  We prove
our main theorem, the Homotopy Colimit Formula, in
Section~\ref{sec:Homotopy_colimit_formula}, work out various
examples in
Section~\ref{sec:examples_of_the_homotopy_colimit_formula}, and
derive the generalized Inclusion-Exclusion Principle in
Section~\ref{sec:Combinatorial_Applications_of_the_Homotopy_Colimit_Formula}.
We review the groupoid cardinality of Baez--Dolan and the Euler
characteristic of Leinster in
Section~\ref{sec:comparison_with_Leinster}, and compare our Homotopy
Colimit Formula with Leinster's compatibility with Grothendieck
fibrations in terms of weightings. We apply our results to Haefliger
complexes of groups in Section~\ref{sec:complexes_of_groups} to
prove
Theorems~\ref{thm:Euler_characteristic_of_hocolim_of_quotient_complex}
and~\ref{the:extension_of_Haefligers_corollary}, which express Euler
characteristics of complexes of groups associated to group actions
in terms of the initial data.

\vspace{.5cm} \noindent {\bf Acknowledgements.} All three authors
were supported by the Sonderforschungsbereich 878 \--- Groups,
Geometry and Actions \--- and the Leibniz-Preis of Wolfgang
L\"uck.

Thomas M.~Fiore was supported at the University of Chicago by NSF
Grant DMS-0501208. At the Universitat Aut\`{o}noma de Barcelona he
was supported by grant SB2006-0085 of the Spanish Ministerio de
Educaci\'{o}n y Ciencia under the Programa Nacional de ayudas para
la movilidad de profesores de universidad e investigadores
espa$\tilde{\text{n}}$oles y extranjeros. Thomas M. Fiore also
thanks the Centre de Recerca Matem\`{a}tica in Bellaterra
(Barcelona) for its hospitality during the CRM Research Program on
Higher Categories and Homotopy Theory in 2007-2008, where he heard
Tom Leinster speak about Euler characteristics. He also thanks the
Max Planck Insitut f\"ur Mathematik for its hospitality and support
during his stay in Summer 2010.

\tableofcontents


\typeout{--------------   The finiteness obstruction and Euler characteristics ---}

\section{The Finiteness Obstruction and Euler Characteristics}
\label{sec:finiteness_obstruction_and_Euler_characteristics}

We quickly recall the main definitions and results needed from our
first paper \cite{FioreLueckSauerFinObsAndEulCharOfCats(2009)} in
order to make this article as self-contained as possible. See
\cite{FioreLueckSauerFinObsAndEulCharOfCats(2009)} for proofs and
more detail.

Throughout this paper, let $\Gamma$ be a category and $R$ an
associative, commutative ring with identity. The first ingredient we
need is the theory of modules over categories developed by L\"uck
\cite{Lueck(1989)}, and recalled in
\cite{FioreLueckSauerFinObsAndEulCharOfCats(2009)}.  An
\emph{$R\Gamma$-module} is a contravariant functor from $\Gamma$
into the category of left $R$-modules. For example, if $\Gamma$ is a
group $G$ viewed as a one-object category, then an $R\Gamma$-module
is the same as a right module over the group ring $RG$. An
$R\Gamma$-module $P$ is \emph{projective} if it is projective in the
usual sense of homological algebra, that is, for every surjective
$R\Gamma$-morphism $p \colon M \to N$ and every $R\Gamma$-morphism
$f \colon P \to N$ there exists an $R\Gamma$-morphism $\overline{f}
\colon P \to M$ such that $p \circ \overline{f} = f$. An
$R\Gamma$-module $M$ is \emph{finitely generated} if there is a
surjective $R\Gamma$-morphism $B(C) \to M$ from an $R\Gamma$-module
$B(C)$ that is free on a collection $C$ of sets indexed by
$\ob(\Gamma)$ such that $\coprod_{x \in \ob(\Gamma)} C_x$ is finite.
Explicitly, the \emph{free $R\Gamma$-module on the $\ob(\Gamma)$-set
$C$} is
\begin{equation} \label{equ:free_RGamma_module}
B(C) := \bigoplus_{x \in \ob(\Gamma)} \bigoplus_{C_x} R\mor_\Gamma(?,x).
\end{equation}
A contravariant $R\Gamma$-module may be tensored with a covariant
$R\Gamma$-module to obtain an $R$-module: if $M\colon \Gamma^{\op}
\to R\text{-}\MOD$  and $N\colon \Gamma \to R\text{-}\MOD$ are
functors, then the \emph{tensor product} $M \otimes_{R\Gamma} N$ is
the quotient of the $R$-module
\begin{equation*}
 \bigoplus_{x \in \ob(\Gamma)} M(x) \otimes_R N(x)
\end{equation*}
by the $R$-submodule generated by elements of the form
\begin{equation*}
(M(f)m) \otimes n - m \otimes (N(f)n)
\end{equation*}
where $f:x \to y$ is a morphism in $\Gamma$, $m \in M(y)$, and $n
\in N(x)$.

Finite projective resolutions of the constant $R\Gamma$-module
$\underline{R}$ play a special role in our theory of Euler
characteristic for categories. A resolution $P_*$ of an
$R\Gamma$-module $M$ is said to be \emph{finite projective} if it
has finite length and each $P_n$ is finitely generated and
projective. We say that a category $\Gamma$ \emph{is of type
(FP$_R$)} if the constant $R\Gamma$-module $\underline{R} \colon
\Gamma^{\op} \to R\text{-}\MOD$ with value $R$ admits a finite
projective resolution. Categories in which every endomorphism is an
isomorphism, the so-called \emph{EI-categories}, provide important
examples. Finite EI-categories in which $|\aut(x)|$ is invertible in
$R$ for each object $x$ are of type (FP$_R$). Further examples of
categories of type (FP$_R$) include categories $\Gamma$ which admit
a finite $\Gamma$-$CW$-model for the classifying $\Gamma$-space
$E\Gamma$ (see Section \ref{sec:spaces_over_a_category} and Examples
\ref{exa:finite_model_for_I_with_terminal_object},
\ref{exa:finite_model_for_parallel_arrows},
\ref{exa:finite_model_for_pushout}, and
\ref{exa:finite_model_for_q_interior}). In fact, such categories
$\Gamma$ are even \emph{of type (FF$_R$)}: the cellular chains on a
finite $\Gamma$-$CW$-model for $E\Gamma$ provide a finite free
resolution of $\underline{R}$. In general, if a category is of type
(FF$_\IZ$), then it is of type (FF$_R$) for any associative,
commutative ring $R$ with identity.

A home for the finiteness obstruction of a category $\Gamma$ is
provided by the \emph{projective class group} $K_0(R\Gamma)$. The
generators of this abelian group are the isomorphism classes of
finitely generated projective $R\Gamma$-modules and the relations
are given by expressions $[P_0] - [P_1] + [P_2] = 0$ for every exact
sequence $0 \to P_0 \to P_1 \to P_2 \to 0$ of finitely generated
projective $R\Gamma$-modules.

\begin{definition}[Finiteness obstruction of a category] \label{def:finiteness_obstruction_of_a_category}
Let $\Gamma$ be a category of type (FP$_R$) and $P_*$ a finite
projective resolution of the constant $R\Gamma$-module
$\underline{R}$. The \emph{finiteness obstruction of $\Gamma$ with
coefficients in $R$} is
$$o(\Gamma;R) := \sum_{n \geq 0} (-1)^n \cdot [P_n] \; \in K_0(R\Gamma).$$
We also use the notation $[\underline{R}]$, or simply $[R]$, to
denote the finiteness obstruction $o(\Gamma;R)$. The finiteness
obstruction, when it exists, does not depend on the choice $P_*$ of
finite projective resolution of $\underline{R}$.
\end{definition}

The finiteness obstruction is compatible with most everything one
could hope for. If $F \colon \Gamma_1 \to \Gamma_2$ is a right
adjoint, and $\Gamma_1$ is of type (FP$_R$), then $\Gamma_2$ is of
type (FP$_R$) and $F_*o(\Gamma_1;R)=o(\Gamma_2;R)$ (here the group
homomorphism $F_*$ is induced by induction with $F$). Since an
equivalence of categories is a right adjoint (and also a left
adjoint), a particular instance of the previous sentence is: if $F
\colon \Gamma_1 \to \Gamma_2$ is an equivalence of categories, then
$\Gamma_1$ is of type (FP$_R$) if and only if $\Gamma_2$ is, and in
this case $F_*o(\Gamma_1;R)=o(\Gamma_2;R)$. The finiteness
obstruction is also compatible with finite coproducts of categories,
finite products of categories, restriction along admissable
functors, and homotopy colimits, as we prove in
Theorem~\ref{the:homotopy_colimit_formula}. If $G$ is a finitely
presented group of type (FP$_\mathbb{Z}$), then Wall's finiteness
obstruction $o(BG)$ is the same as $o(\widehat{G};\mathbb{Z})$,
which is the finiteness obstruction of $G$ viewed as a one-object
category $\widehat{G}$ with morphisms $G$. The finiteness obstruction in Definition~\ref{def:finiteness_obstruction_of_a_category} is a special case of the finiteness obstruction  of a finitely dominated $R\Gamma$-chain complex $C$, denoted $o(C)\in K_0(R\Gamma)$. The image of $o(C)$ in the reduced $K$-theory $\tilde{K}_0(R\Gamma)$ vanishes if and only if $C$ is $R\Gamma$-homotopy equivalent to a finite free $R\Gamma$-chain complex, see \cite[Chapter 11]{Lueck(1989)}.

We will occasionally work with directly finite categories. A
category is called \emph{directly finite} if for any two objects $x$
and $y$ and morphisms $u \colon x \to y$ and $v \colon y \to x$ the
implication $vu = \id_x \implies uv = \id_y$ holds. If $\Gamma_1$
and $\Gamma_2$ are equivalent categories, then $\Gamma_1$ is
directly finite if and only if $\Gamma_2$ is directly finite.
Examples of directly finite categories include groupoids, and more
generally EI-categories.

A key result in the theory of modules over an EI-category is
L\"uck's splitting of the projective class group of $\Gamma$ into
the projective class groups of the automorphism groups
$\aut_\Gamma(x)$, one \for each isomorphism class of objects. We
next recall the relevant maps and notation. For $x \in \ob(\Gamma)$,
we denote $R\aut_\Gamma(x)$ by $R[x]$ for simplicity. The
\emph{splitting functor at $x \in \ob(\Gamma)$}
\begin{eqnarray} \label{equ:splitting_functor}
  & S_x \colon \MOD\text{-}R\Gamma \to \MOD\text{-}R[x], &
  \label{S_x}
\end{eqnarray}
maps an $R\Gamma$-module $M$ to the quotient of the $R$-module
$M(x)$ by the $R$-submodule generated by all images of $R$-module
homomorphisms $M(u)\colon M(y) \to M(x)$ induced by all
non-invertible morphisms $u\colon x \to y$ in $\Gamma$. The right
$R[x]$-module structure on $M(x)$ induces a right $R[x]$-module
structure on $S_xM$. Note that $S_xM$ is an $R[x]$-module, not an
$R\Gamma$-module.  The functor $S_x$ respects direct sums, sends
epimorphisms to epimorphisms, and sends free modules to free
modules. If $\Gamma$ is directly finite,
then $S_x$ also preserves finitely generated and projective. The
\emph{extension functor at $x \in \ob(\Gamma)$}
\begin{eqnarray} \label{equ:extension_functor}
  & E_x \colon  \MOD\text{-}R[x] \to  \MOD\text{-}R\Gamma&
  \label{E_x}
\end{eqnarray}
maps an $R[x]$-module $N$ to the  $R\Gamma$-module
$N \otimes_{R[x]} R\mor_\Gamma(?,x)$. The functor
$E_x$ respects direct sums, sends epimorphisms to epimorphisms, sends free modules to free modules, and preserves finitely generated and projective. If $\Gamma$ is directly finite, and $P$ is a projective $R[x]$-module, then there is a natural isomorphism $P \cong S_xE_xP$ compatible with direct sums.
\begin{theorem}[Splitting of $K_0(R\Gamma)$ for EI-categories, Theorem~10.34 on page~196 of L\"uck \cite{Lueck(1989)}]
\label{the_splitting_of_K-theory_for_EI_categories} If $\Gamma$ is
an EI-category, then the group homomorphisms
$$\xymatrix{K_0(R\Gamma) \ar@<.4ex>[r]^-S & \ar@<.4ex>[l]^-E }\Split K_0(R\Gamma) :=
\bigoplus_{\overline{x} \in \iso(\Gamma)} K_0(R\aut_\Gamma(x))$$
defined by
\[
S[P]=  \{[S_xP] \mid \overline{x} \in \iso(\Gamma)\}
\]
and
\[
E \{[Q_x] \mid \overline{x} \in \iso(\Gamma)\}=
\sum_{\overline{x} \in \iso(\Gamma)} [E_xQ_x],
\]
are isomorphisms and inverse to one another. They are covariantly natural with respect to functors
between EI-categories.
\end{theorem}

\begin{remark}
If $\Gamma$ is not an EI-category, then the splitting homomorphism
$S\colon K_0(R\Gamma) \to \Split K_0(R\Gamma)$ may not be an
isomorphism. However, $S$ is covariantly natural with respect to
functors between directly finite categories, see
\cite[Lemma~3.15]{FioreLueckSauerFinObsAndEulCharOfCats(2009)}.
\end{remark}

The splitting functors $S_x$ allow us to define the notion of $R\Gamma$-rank $\rk_{R\Gamma}$ for finitely generated $R\Gamma$-modules, which in turn
allows the definition of the functorial Euler characteristic, as we explain next. We assume a fixed notion of a rank $\rk_R(N) \in \IZ$ for finitely generated $R$-modules $N$ such that $\rk_R(R) = 1$ and $\rk_R(N_1) = \rk_R(N_0) + \rk_R(N_2)$ for any sequence
$0 \to N_0 \to N_1 \to N_2 \to 0$ of finitely generated $R$-modules. If $R$ is a commutative principal ideal domain, we use $\rk_R(N) := \dim_F(F \otimes_R N)$,
where $F$ is the quotient field of $R$. Let $U(\Gamma)$ be the free abelian group on the set
of isomorphism classes of objects in $\Gamma$, that is $U(\Gamma)
:= \IZ\iso(\Gamma).$  The augmentation homomorphism
$\epsilon\colon U(\Gamma) \to \IZ$ adds up the components of an element of $U(\Gamma)$.
\begin{definition}[Rank of a finitely generated $R\Gamma$-module]
\label{def:rank_of_fin_gen_prof_RGamma-module} If $M$ is a finitely
generated $R\Gamma$-module $M$, then its \emph{$R\Gamma$-rank} is
$$\rk_{R\Gamma}(M) :=
\bigl\{\rk_R(S_xM \otimes_{R[x]} R) \mid \overline{x} \in
\iso(\Gamma)\bigr\} \quad \in U(\Gamma).$$
\end{definition}

\begin{definition}[The (functorial) Euler characteristic of a category]
  \label{def:functorial_Euler_characteristic_of_a_category}\label{def:Euler_characteristic_of_a_category}
  Suppose that $\Gamma$ is of type (FP$_R$).
  The \emph{functorial Euler  characteristic of $\Gamma$ with coefficients in $R$}
is the image of the finiteness obstruction $o(\Gamma;R) \in K_0(R\Gamma)$ under the
  homomorphism $\rk_{R\Gamma}\colon K_0(R\Gamma) \to U(\Gamma)$, that is
\[
\chi_f(\Gamma;R):=\rk_{R\Gamma} o(\Gamma;R)=\left\{\sum_{n \geq 0
}(-1)^n\rk_R(S_xP_n \otimes_{R[x]} R) \mid \overline{x} \in
\iso(\Gamma) \right\} \quad \in U(\Gamma),
\]
where $P_*$ is any finite projective $R\Gamma$-resolution of the constant $R\Gamma$-module $\underline{R}$.
  The \emph{Euler characteristic of $\Gamma$ with coefficients in $R$}
  is the sum of the components of the functorial Euler characteristic, that is, \[
\chi(\Gamma;R):=\epsilon(\chi_f(\Gamma;R))=\sum_{\overline{x} \in \iso(\Gamma)} \sum_{n \geq 0 }(-1)^n\rk_R(S_xP_n \otimes_{R[x]} R).
\]
\end{definition}

For example, if $\calg$ is a finite groupoid, then $\chi_f(\calg) \in U(\calg)$ is
$(1,1,\dots,1)$, and $\chi(\calg)$ counts the isomorphism classes of objects, or equivalently the connected components, of $\calg$.

\begin{theorem}[Theorem~4.20~of~Fiore--L\"uck--Sauer \cite{FioreLueckSauerFinObsAndEulCharOfCats(2009)}]
\label{the:chi_f_determines_chi} Let $R$ be a Noetherian ring and
$\Gamma$ a directly finite category of type (FP$_R$). Then the Euler
characteristic and topological Euler characteristic of $\Gamma$
agree.  That is, $H_n(B\Gamma;R)$ is a finitely generated $R$-module
for every $n \geq 0$, there exists a natural number $d$ with
$H_n(B\Gamma;R) = 0$ for all $n > d$, and
$$\chi(\Gamma;R) = \chi(B\Gamma;R)= \sum_{n \geq 0} (-1)^n \cdot \rk_{R}(H_n(B\Gamma;R)) \in \IZ.$$
Here $\chi(\Gamma;R)$ is defined in
Definition~\ref{def:Euler_characteristic_of_a_category} and
$B\Gamma$ denotes the geometric realization of the nerve of
$\Gamma$.
\end{theorem}

The functorial Euler characteristic and Euler characteristic have
many desirable properties. They are invariant under equivalence of
categories and are compatible with finite products and finite
coproducts. As we prove in
Theorem~\ref{the:homotopy_colimit_formula}, they are also compatible
with homotopy colimits.

The $L^2$-Euler characteristic, which is in some sense the better
invariant, can be defined similarly by taking $R=\IC$ and using the
$L^2$-rank $\rk^{(2)}_\Gamma$ rather than the $R\Gamma$-rank. For
this we need group von Neumann algebras and their dimension theory
from L\"uck \cite{Lueck(1998a)} and \cite{Lueck(1998b)}, as recalled
in our first paper
\cite{FioreLueckSauerFinObsAndEulCharOfCats(2009)} for the purpose
of Euler characteristics. If $G$ is a group, its \emph{group von
Neumann algebra}
\begin{equation*}
\caln(G) =  \calb(l^2(G))^G
\end{equation*}
is the algebra of bounded operators on $l^2(G)$ that are equivariant with respect to the right $G$-action. If $G$ is finite, $\caln(G)$ is the group
ring $\IC G$. In any case, the group ring $\IC G$ embeds as a subring of $\caln(G)$ by sending $g \in G$ to the isometric
$G$-equivariant operator $l^2(G) \to l^2(G)$ given by left multiplication with $g$. In particular, we can view $\caln(G)$ as a $\IC
G$-$\caln(G)$-bimodule and tensor $\IC G$-modules on the right with $\caln(G)$. If $G$ is the automorphism group of an object in $\Gamma$, then we write $\caln(x)$ for $\caln\bigl(\aut_\Gamma(x)\bigr)$.

The \emph{von Neumann dimension}, $\dim_{\caln(G)}$, is a function
that assigns to \emph{every} right $\caln(G)$-module $M$ a
non-negative real number of $\infty$. It is the unique such function
which satisfies Hattori-Stallings rank, additivity, cofinality, and
continuity. If $G$ is a finite group, then $\caln(G) = \IC G$ and we
get for a $\IC G$-module $M$ the von Neumann dimension
$$\dim_{\caln(G)}(M) = \frac{1}{|G|} \cdot \dim_{\IC}(M),$$ where
$\dim_{\IC}$ is the dimension of $M$ viewed as a complex vector
space. A category $\Gamma$ is said to be \emph{of type }($L^2$) if
for one (and hence every) projective $\IC\Gamma$-resolution $P_*$ of
the constant $\IC \Gamma$-module $\underline{\IC}$ we have
$$\sum_{\overline{x} \in \iso(\Gamma)} \sum_{n \geq 0} \dim_{\caln(x)} H_n\bigl(S_x P_*\otimes_{\IC [x]}\caln(x)\bigr) < \infty.$$
Note that the projective resolution $P_*$ of $\underline{\IC}$ is
not required to be of finite length, nor finitely generated.
Examples of categories of type ($L^2$) include finite EI-categories,
in particular finite posets and finite groupoids. Infinite
categories can also be of type ($L^2$), for example any (small)
groupoid with finite automorphism groups such that
\begin{equation}\label{equ:finiteness of aut-sum}
  \sum_{\overline{x} \in \iso(\calg)} \frac{1}{|\aut_{\calg}(x)|}<\infty
\end{equation}
holds is of type ($L^2$). The condition of type ($L^2$) is weaker
than (FP$_\IC$), since any directly finite category of type
(FP$_\IC$) is also of type ($L^2$).

\begin{definition}[The (functorial) $L^2$-Euler characteristic of a category]
  \label{def:functorial_L2-Euler_characteristic_of_a_category} \label{def:L2-Euler_characteristic_of_a_category}
  Suppose that $\Gamma$ is of type ($L^2$). Define
  $$U^{(1)}(\Gamma) :=
  \left\{\sum_{\overline{x} \in \iso(\Gamma)} r_{\overline{x}} \cdot \overline{x}
  \;\bigg|\;
  r_{\overline{x}} \in \IR,
  \sum_{\overline{x} \in \iso(\Gamma)} |r_{\overline{x}}| < \infty\right\}\subseteq\prod_{\bar{x}\in\iso(\Gamma)}\IR.$$
  The \emph{functorial $L^2$-Euler characteristic of $\Gamma$} is
\[
\chi_f^{(2)}(\Gamma):=\left\{\sum_{n \geq 0} (-1)^n \dim_{\caln(x)} H_n\bigl(S_x P_*\otimes_{\IC [x]}\caln(x)\bigr)\mid \bar{x} \in
\iso(\Gamma) \right\} \in U^{(1)}(\Gamma),
\]
where $P_*$ is any
projective $\IC \Gamma$-resolution of the constant
$\IC\Gamma$-module
  $\underline{\IC}$. The \emph{$L^2$-Euler characteristic of $\Gamma$} is the sum over $\bar{x} \in
\iso(\Gamma)$ of the components of the functorial Euler characteristic, that is,
\[
\chi^{(2)}(\Gamma):=\sum_{\overline{x} \in \iso(\Gamma)} \sum_{n
\geq 0} (-1)^n \dim_{\caln(x)} H_n\bigl(S_x P_*\otimes_{\IC
[x]}\caln(x)\bigr).
\]
\end{definition}

If $\calg$ is a groupoid such that \eqref{equ:finiteness of aut-sum}
holds, then the functorial $L^2$-Euler characteristic
$\chi_f^{(2)}(\calg) \in \prod_{\overline{x} \in \iso(\calg)} \IR$
has at $\overline{x} \in \iso(\calg)$ the value $1/|\aut_\calg(x)|$.
The $L^2$-Euler characteristic is
\begin{equation} \label{eq:L2_Euler_Characteristic_of_Groupoid}
\chi^{(2)}(\calg) = \sum_{\overline{x} \in \iso(\calg)}
\frac{1}{|\aut_\calg(x)|}.
\end{equation}
See Lemma~\ref{lem:chi(2)_and_chi} for an explicit formula for
$\chi^{(2)}(\Gamma)$ in the case of a finite, skeletal EI-category
$\Gamma$ in which the left $\aut_\Gamma(y)$-action on
$\mor_\Gamma(x,y)$ is free for every two objects $x,y \in
\ob(\Gamma)$.

\begin{definition}[$L^2$-rank of a finitely generated $\IC\Gamma$-module]
\label{def:L2rank_of_fin_gen_prof_RGamma-module}
Let $M$ be a finitely generated $\IC\Gamma$-module $M$.  Its \emph{$L^2$-rank} is
\begin{equation*}
 \rk_{\Gamma}^{(2)}(M)  :=
\bigl\{\dim_{\caln(x)}(S_xM \otimes_{\IC[x]} \caln(x))\mid \bar{x}
\in \iso(\Gamma)\bigr\}
\in U(\Gamma) \otimes_{\IZ} \IR = \bigoplus_{\iso(\Gamma)} \IR.
\end{equation*}
\end{definition}

\begin{theorem}[Relating the finiteness obstruction and the
  $L^2$-Euler characteristic, Theorem~5.22~of~Fiore--L\"uck--Sauer \cite{FioreLueckSauerFinObsAndEulCharOfCats(2009)}]
  \label{the:comparing_o_and_chi(2)}
  Suppose that $\Gamma$ is a directly finite category of type (FP$_{\IC}$).
  Then
  $\Gamma$ is of type ($L^2$) and the image of the finiteness
  obstruction $o(\Gamma;\IC)$ (see
  Definition~\ref{def:finiteness_obstruction_of_a_category})
  under the homomorphism
  $$\rk_{\Gamma}^{(2)} \colon K_0(\IC\Gamma) \to U(\Gamma) \otimes_{\IZ} \IR
  = \bigoplus_{\overline{x} \in \iso(\Gamma)} \IR$$ is the functorial $L^2$-Euler characteristic $\chi_f^{(2)}(\Gamma)$.
\end{theorem}

The $L^2$-Euler characteristic agrees with the groupoid cardinality
of Baez--Dolan \cite{Baez-Dolan(2001)} and the Euler characteristic
of Leinster \cite{Leinster(2008)} in certain cases, see
Lemma~\ref{lem:chi(2)_and_chi} and
Section~\ref{sec:comparison_with_Leinster}. In particular, the
Baez--Dolan groupoid cardinality of a groupoid
satisfying~\eqref{equ:finiteness of aut-sum} is
\eqref{eq:L2_Euler_Characteristic_of_Groupoid}. However, the
Baez--Dolan groupoid cardinality and Leinster's Euler characteristic
$\chi_L(\Gamma)$ only depend on the underlying graph of $\Gamma$,
whereas our invariants truly depend on the category structure. For
instance, $\chi_L$ is $\frac{1}{2}$ for both the two-element monoid
$(\IZ/2,\times)$ and the two-element group $(\IZ/2,+)$, whereas
$\chi^{(2)}$ is 1 respectively $\frac{1}{2}$. The distinction can
already be seen on the level of the finiteness obstructions. The
Euler characteristic $\chi(-)$ and topological Euler characteristic
$\chi(B-)$ can also distinguish categories with the same underlying
directed graph as in the following example. For $S=\{1,2,3,4\}$,
$G_1=\langle(1234)\rangle$, $G_2=\langle (12),(34) \rangle$, and
$k=1,2$, let $\Gamma_k$ be the EI-category with objects $x$ and $y$
and $\mor(x,y):=S$, $\mor(x,x):=\{\id_x\}$, $\mor(y,y):=G_k$, and
$\mor(y,x)=\emptyset$. Composition in $\Gamma_k$ is the composition
in $G_k$ and the left $G_k$-action on $S$, that is, $\Gamma_k$ is
the EI-category associated to the respective $G_k$-$\{1\}$-biset $S$
as in Subsection~6.4 of Fiore--L\"uck--Sauer
\cite{FioreLueckSauerFinObsAndEulCharOfCats(2009)}. Then $\Gamma_1$
and $\Gamma_2$ have the same underlying directed graph but
$\chi(\Gamma_1;\IQ)=\chi(B\Gamma_1;\IQ)=1$ and
$\chi(\Gamma_2;\IQ)=\chi(B\Gamma_2;\IQ)=0$ by Theorem~6.23~(iii) of
Fiore--L\"uck--Sauer
\cite{FioreLueckSauerFinObsAndEulCharOfCats(2009)}. An infinite
example of categories with the same underlying graph but different
Euler characteristics is provided by the groups $\IZ$ and $\IZ*\IZ$,
each of which admits a finite $\Gamma$-$CW$-model for its respective
$\Gamma$-classifying space. The categories $\widehat{\IZ}$ and
$\widehat{\IZ*\IZ}$ have the same underlying directed graph, but we
have $\chi^{(2)}(\widehat{\IZ})=0\neq\chi^{(2)}(\widehat{\IZ*\IZ})$,
and similarly for $\chi$. Typically, the Euler characteristic of a
category $\Gamma_{\text{free}}$ free on a directed graph $(V,E)$ is
the same as the Euler characteristic of the directed graph $(V,E)$.
For the topological Euler characteristic this is clearly true, since
$B\Gamma_{\text{free}}$ is homotopy equivalent to the topological
realization $|(V,E)|$. If $\Gamma_{\text{free}}$ is directly finite
and $R$ is Noetherian, then we also have
$\chi(\Gamma_{\text{free}})=\chi(|(V,E)|)$ by
Theorem~\ref{the:chi_f_determines_chi}. For example for the directed
graph with one vertex and one arrow we have
$\chi(\widehat{\IN})=0=\chi(S^1)$.

The functorial $L^2$-Euler characteristic and the $L^2$-Euler
characteristic have many desirable properties. They are invariant
under equivalence of categories and are compatible with finite
products, finite coproducts, and isofibrations and coverings between
finite groupoids. We prove in
Theorem~\ref{the:homotopy_colimit_formula} the compatibility with
homotopy colimits. In the case of a group $G$, the $L^2$-Euler
characteristic of $\widehat{G}$ coincides with the classical
$L^2$-Euler characteristic of $G$, which is $1/\vert G \vert$ when
$G$ is finite. The $L^2$-Euler characteristic is also closely
related to the geometry and topology of the classifying space for
proper $G$-actions, namely the functorial $L^2$-Euler characteristic
of the proper orbit category $\uor{G}$ is equal to the equivariant
Euler characteristic of the classifying space $\underline{E}G$ for
proper $G$-actions, whenever $\underline{E}G$ admits a finite
$G$-$CW$-model.

The question arises: what are sufficient conditions for the Euler
characteristic and $L^2$-Euler characteristic to coincide with the Euler characteristic of the classifying space?
This is answered in the following Theorem.

\begin{theorem}[Invariants agree for directly finite and type (FF$_\IZ$), Theorem~5.25~of Fiore--L\"uck--Sauer~\cite{FioreLueckSauerFinObsAndEulCharOfCats(2009)}]
\label{the:coincidence} Suppose $\Gamma$ is directly finite and of
type (FF$_\IZ$). Then the functorial $L^2$-Euler characteristic of
Definition~\ref{def:functorial_L2-Euler_characteristic_of_a_category}
coincides with the functorial Euler characteristic of
Definition~\ref{def:functorial_Euler_characteristic_of_a_category}
for any associative, commutative ring $R$ with identity
\[ \chi_f^{(2)}(\Gamma)=\chi_f(\Gamma;R)\in U(\Gamma) \subseteq
U^{(1)}(\Gamma),
\]
and thus $\chi^{(2)}(\Gamma)=\chi(\Gamma;R)$ in
Definition~\ref{def:L2-Euler_characteristic_of_a_category} and
Definition~\ref{def:Euler_characteristic_of_a_category}.

If $R$ is additionally Noetherian, then
\begin{equation} \label{equ:coincidence_for_type_FFZ_and_directly_finite}
\chi(B\Gamma;R)=\chi(\Gamma;R)=\chi^{(2)}(\Gamma).
\end{equation}
Moreover, if $\Gamma$ is merely of type (FF$_\IC$) rather than
(FF$_\IZ$), then equation
\eqref{equ:coincidence_for_type_FFZ_and_directly_finite} holds for
any Noetherian ring $R$ containing $\IC$.
\end{theorem}

Any category $\Gamma$ which admits a finite $\Gamma$-$CW$-model in
the sense of Section~\ref{sec:spaces_over_a_category} is of type
(FF$_R$) for any ring $R$, by an application of the cellular
$R$-chain functor. Thus, Theorem~\ref{the:coincidence} applies to
any directly finite category $\Gamma$ which admits a finite
$\Gamma$-$CW$-model. For example, finite categories without loops
are directly finite and admit finite models
(Lemma~\ref{lem:scwol_directly_finite_EI} and
Theorem~\ref{the:finite_models_for_finite_scwols}), so
equation~\eqref{equ:coincidence_for_type_FFZ_and_directly_finite}
holds for instance for $\{j \rightrightarrows k\}$, $\{k \leftarrow
j \to \ell\}$, and finite posets. The monoid $\IN$ and group $\IZ$,
viewed as one-object categories $\widehat{\IN}$ and $\widehat{\IZ}$,
are also directly finite and admit finite models (see
Example~\ref{exa:finite_model_for_N_and_Z}), so we have
\[
0=\chi(S^1;R)=\chi(B\widehat{\IN};R)=\chi(\widehat{\IN};R)=\chi^{(2)}(\widehat{\IN})
\]
and
\[
0=\chi(S^1;R)=\chi(B\widehat{\IZ};R)=\chi(\widehat{\IZ};R)=\chi^{(2)}(\widehat{\IZ})
\]
($B\widehat{\IN} \to B\widehat{\IZ}\simeq S^1$ is a homotopy
equivalence by Quillen's Theorem A, see
Rabrenovi{\'c}~\cite[Proposition 10]{Rabrenovic(2005)}). The
equations $\chi(\widehat{\IN};R)=0=\chi^{(2)}(\widehat{\IN})$ and
$\chi(\widehat{\IZ};R)=0=\chi^{(2)}(\widehat{\IZ})$ also follow from
Example~\ref{exa:homotopy_hocolimit_and_trivial_functor}, since the
finite models for $\widehat{\IN}$ and $\widehat{\IZ}$ in
Example~\ref{exa:finite_model_for_N_and_Z} each have one
$\cali$-$0$-cell and one $\cali$-$1$-cell.

We may use Theorem~\ref{the:coincidence} to obtain an explicit
formula for Euler characteristics of finite categories without loops
as follows. Let $\Gamma$ be a finite category without loops, and
choose a skeleton $\Gamma'$. Let $c_n(\Gamma')$ denote the number of
paths
\[
i_0 \to i_1 \to i_2 \to \cdots \to i_n
\]
of $n$-many non-identity morphisms in $\Gamma'$. Then $c_n(\Gamma')$
is the number of $n$-cells in the $CW$-complex $B\Gamma'$, and we
have
\begin{equation} \label{equ:Euler_characteristic_for_skeletal_finite_scwols}
\chi(\Gamma;R)=\chi^{(2)}(\Gamma)=\chi(B\Gamma;R)=\chi(B\Gamma';R)=\sum_{n
\geq 0} (-1)^n c_n(\Gamma').
\end{equation}
See~\cite[Corollary~1.5]{Leinster(2008)} for a different derivation
of this formula for Leinster's Euler characteristic $\chi_L(\Gamma)$
in the case $\Gamma$ was already skeletal. See also
Examples~\ref{exa:homotopy_hocolimit_and_trivial_functor}~and~\ref{exa:Euler_characteristics_of_finite_scwols}
where skeletality of $\cali$ is not required.

\begin{remark}[Homotopy Invariance]
If $F:\Gamma_1 \to \Gamma_2$ is a functor such that $BF$ is a homotopy equivalence, and
both $\Gamma_1$ and $\Gamma_2$ are of type (FP$_R$), and if
\[
\chi(\Gamma_1;R)=\chi(B\Gamma_1;R) \;\;\text{ and }\;\;\chi(\Gamma_2;R)=\chi(B\Gamma_2;R),
\]
then clearly $\chi(\Gamma_1;R)=\chi(\Gamma_2;R)$. However, it is
possible for two categories to be homotopy equivalent, one of which
is (FP$_R$) and the other is not, so that one has a notion of Euler
characteristic and the other does not. In
Section~10~of~Fiore--L\"uck--Sauer~\cite{FioreLueckSauerFinObsAndEulCharOfCats(2009)}
such an example is discussed.
\end{remark}

\section{Spaces over a Category} \label{sec:spaces_over_a_category}

An important hypothesis in our Homotopy Colimit Formula involves the
idea of a space over a category, see
Davis--L\"uck~\cite{Davis-Lueck(1998)}. Namely, we assume that the
indexing category $\cali$ for the diagram $\calc$ of categories
admits a finite $\cali$-$CW$-model for its $\cali$-classifying
space. Essentially this means it is possible to functorially assign
a contractible $CW$-complex $E\cali(i)$ to each $i \in \ob(\cali)$,
and moreover, these local $CW$-complexes are constructed globally by
gluing $\cali$-$n$-cells of the form $\mor_\cali(-,i_\lambda)\times
D^n$ onto the already globally constructed $(n-1)$-skeleton
$E\cali_n$. The Homotopy Colimit Formula then expresses the
invariants of the homotopy colimit of $\calc$ in terms of the
invariants of the categories $\calc(i_\lambda)$ at the base objects
$i_\lambda$ for $E\cali$.

The gluing described above takes place in the more general category
of $\cali$-spaces. A \emph{(contravariant) $\cali$-space} is a
contravariant functor from $\cali$ to the category $\SPACES$ of
(compactly generated) topological spaces. As usual, we will always
work in the category of compactly generated spaces
(see~Steenrod~\cite{Steenrod(1967)}). A \emph{map between
$\cali$-spaces} is a natural transformation.  Given an object $i \in
\ob(\cali)$, we obtain an $\cali$-space $\mor_\cali(?,i)$ which
assigns to an object $j$ the discrete space $\mor_\cali(j,i)$.

The next definition is taken
from~Davis--L\"uck~\cite[Definition~3.2]{Davis-Lueck(1998)}, where
an $\cali$-$CW$-complex is called a free $\cali$-$CW$-complex and we
will omit the word free here. The more general notion of
$\cali$-$CW$-complex was defined by Dror~Farjoun~\cite[1.16 and
2.1]{DrorFarjoun(1987)}. See also~Piacenza~\cite{Piacenza(1991)}.

\begin{definition}[$\cali$-$CW$-complex]
  \label{def:calc-CW-complex}
  A \emph{(contravariant) $\cali$-$CW$-complex} $X$ is a contravariant $\cali$-space
  $X$ together with a filtration
  $$\emptyset = X_{-1} \subset X_0 \subset X_1 \subset X_2  \subset\ldots
  \subset X_n \subset \ldots \subset X = \bigcup_{n \geq 0} X_n$$ such
  that $X = \colim_{n \to \infty} X_n$ and for any $n \geq 0$ the
  \emph{$n$-skeleton} $X_n$ is obtained from the $(n-1)$-skeleton
  $X_{n-1}$ by {\it attaching $\cali$-$n$-cells}, i.e., there exists a pushout
  of $\cali$-spaces of the form
  $$
  \comsquare{\coprod_{\lambda \in \Lambda_n} \mor_{\cali}(-,i_\lambda) \times S^{n-1}}
  {}{X_{n-1}}{}{} {\coprod_{\lambda \in \Lambda_n} \mor_{\cali}(-,i_\lambda) \times D^n}
  {}{X_n}
  $$
  where the vertical maps are inclusions, $\Lambda_n$ is an index set, and the
  $i_\lambda$-s are objects of $\cali$. In particular, $X_0=\coprod_{\lambda \in \Lambda_0} \mor_{\cali}(-,i_\lambda)$.

  We refer to the inclusion functor
  $\mor_{\cali}(-,i_\lambda) \times  (D^n-S^{n-1}) \to X$
  as an \emph{$\cali$-$n$-cell based at $i_\lambda$}.

  An $\cali$-$CW$-complex has \emph{dimension $\le n$} if $X = X_n$.
  We call $X$ \emph{finite dimensional} if there exists an
  integer $n$ with $X = X_n$. It is called \emph{finite} if it is finite
  dimensional and $\Lambda_n$ is finite for every $n \geq 0$.

  The definition of a \emph{covariant $\cali$-$CW$-complex} is analogous.
\end{definition}

\begin{definition}[Classifying $\cali$-space]
\label{def:classifying_calc-space} A model for \emph{the classifying
$\cali$-space} $E\cali$ is an $\cali$-$CW$-complex $E\cali$ such that
$E\cali(i)$ is contractible for all objects $i$.
\end{definition}

The universal property of $E\cali$ is that for any
$\cali$-$CW$-complex $X$ there is up to homotopy precisely one map
of $\cali$-spaces from $X$ to $E\cali$. In particular two models for
$E\cali$ are $\cali$-homotopy equivalent
(see~Davis--L\"uck~\cite[Theorem~3.4]{Davis-Lueck(1998)}). A model
for the usual \emph{classifying space} $B\cali$ is given by $E\cali
\otimes_{\cali} \pt$
(see~\cite[Definition~3.10]{Davis-Lueck(1998)}), where $\pt$ is the
constant covariant $\cali$-space with value the one point space and
$\otimes_{\cali}$ denotes the tensor product of a contravariant and
a covariant $\cali$-space as follows
(see~\cite[Definition~1.4]{Davis-Lueck(1998)}).

\begin{definition}[Tensor product of a contravariant and
a covariant $\cali$-space] \label{def:tensor_product_for_I-spaces}
Let $X$ be a contravariant $\cali$-space and $Y$ a covariant
$\cali$-space. The \emph{tensor product of $X$ and $Y$} is
\[ X \otimes_\cali Y = \Biggl( \coprod_{i \in \cali}  X(i) \times Y(i) \Biggr) / \sim
\]
where $(X(\phi)(x),y)\sim(x,Y(\phi)y)$ for all morphisms $\phi:i \to j$ in $\cali$ and points $x \in X(j)$ and $y \in Y(i)$.
\end{definition}

We present some examples of classifying $\cali$-spaces for various
categories $\cali$.

\begin{example} \label{exa:finite_model_for_I_with_terminal_object}
If $\cali$ has a terminal object $t$, then a finite model for the
classifying $\cali$-space $E\cali$ is simply $\mor_\cali(-,t)$.
\end{example}

\begin{example} \label{exa:finite_model_for_parallel_arrows}
Let $\cali=\{j \rightrightarrows k\}$ be the category consisting of
two objects and a single pair of parallel arrows between them. All
other morphisms are identity morphisms. We obtain a finite model $X$
for the classifying $\cali$-space $E\cali$ as follows.  The
$\cali$-$CW$-space $X$ has a single $\cali$-$0$-cell based at $k$ and a single
$\cali$-1-cell based at $j$. The gluing map $\mor_\cali(-,j) \times S^0 \to
\mor_\cali(-,k)$ is induced by the two parallel arrows $j
\rightrightarrows k$. Then $X(j)=D^1\simeq \ast$ and $X(k)=\ast$.
\end{example}

\begin{example}
\label{exa:finite_model_for_pushout} Let $\cali=\{k \leftarrow j \to
\ell \}$ be the category with objects $j$, $k$ and $\ell$, and precisely
one morphism from $j$ to $k$ and one morphism from $j$ to $\ell$. All
other morphisms are identity morphisms.
  A finite model for $E\cali$ is given by the $\cali$-$CW$-complex with
  precisely two $\cali$-$0$-cells $\mor_\cali(?,k)$ and $\mor_\cali(?,\ell)$ and precisely
  one $\cali$-$1$-cell $\mor_\cali(?,j) \times D^1$ whose attaching map $\mor_\cali(?,j)
  \times S^0 \to \mor_\cali(?,k) \amalg\mor_\cali(?,\ell)$ is the disjoint union of
  the canonical maps $\mor_\cali(?,j) \to \mor_\cali(?,k)$ and $\mor_\cali(?,j)\to
  \mor_\cali(?,\ell)$.  The value of this $1$-dimensional $\cali$-$CW$-complex
  at the objects $k$ and $\ell$ is a point and at the object $j$ is
  $D^1$. Hence it is a finite model for $E\cali$.
\end{example}

\begin{example}
\label{exa:finite_model_for_q_interior} Let $\cali$ be the category
with objects the non-empty subsets of $[q]=\{0,1, \dots, q\}$ and a
unique arrow $J \to K$ if and only if $K \subseteq J$. In other
words $\cali$ is the {\it opposite} of the poset of non-empty subsets of
$[q]$. Then $\cali$ admits a finite $\cali$-$CW$-model $X$ for the
classifying $\cali$-space $E\cali$ as follows. The functor $X\colon \cali^{\op} \to \SPACES$ assigns to $L$
the space $|\Delta[L]|$, which is the geometric realization of the simplicial set which maps $[m]$ to the set of weakly order preserving maps $[m] \to L$. The space $|\Delta[L]|$ is homeomorphic to the standard simplex with $\text{card}(L)$ vertices. The $n$-skeleton $X_n$ of $X$ sends each $L$ to the $n$-skeleton of $|\Delta[L]|$.  The $\cali$-cells of $X$ are attached globally in the following way. The 0-skeleton is
$$X_0=\coprod_{J \subseteq [q], |J|=1}\mor_\cali(-,J).$$ For $n\leq q$, we
construct $X_n$ out of $X_{n-1}$ as the pushout
$$\xymatrix{\coprod_J \mor_\cali(-,J) \times |\partial \Delta[n]| \ar[r] \ar[d] & X_{n-1} \ar[d]
\\ \coprod_J \mor_\cali(-,J) \times |\Delta[n]| \ar[r] & X_n.}$$ The disjoint unions are over all
$J \subseteq [q]$ with $|J|=n+1$. The $J$-component of the gluing
map is induced by the $(n-1)$-face inclusion
$$\xymatrix{|\Delta[K]| \ar[r] & \partial |\Delta[J]| \cong \partial |\Delta[n]|}$$
for all $K \subseteq J$ with $|K|=n$. Clearly $X$ is a finite
$\cali$-$CW$-complex. For each object $L$ of $\cali$, we have
$X(L)=|\Delta[L]| \simeq \ast$, so that $X$ is a finite model
for $E\cali$.
\end{example}

\begin{example} \label{exa:finite_model_for_N_and_Z}
Infinite categories may also admit finite models. Let $\cali=\widehat{\IN}$
be the monoid of natural numbers $\IN$ viewed as a one-object category.
A finite model $X$ for the $\widehat{\IN}$-classifying
space has $X_0(*)=\mor_{\widehat{\IN}}(*,*)=\IN$ and $X_1(*)=[0,\infty)$.
This model has a single $\widehat{\IN}$-$0$-cell $\mor_{\widehat{\IN}}(-,*)$ and a single $\widehat{\IN}$-$1$-cell $\mor_{\widehat{\IN}}(-,*)\times D^1$. The gluing map
$\IN \times S^0 \to \IN$ sends $(n,-1)$ and $(n,1)$ to $n$ and $n+1$ respectively.
Similarly, the group of integers $\IZ$ viewed as a one object category admits a finite model $Y$ with one $\widehat{\IZ}$-$0$-cell and one $\widehat{\IZ}$-$1$-cell, so that $Y_0(*)=\IZ$ and $Y_1(*)=\IR$.
\end{example}

\begin{remark}
Suppose a category $\cali$ admits a finite $\cali$-$CW$-model for $E\cali$. Then the cellular $R$-chains of a finite model provide a finite free resolution
of the constant $R\cali$-module $\underline{R}$, so $\cali$ is of type (FF$_R$). If $\cali$ is additionally directly finite and $R$ is Noetherian, then $\chi(B\cali;R)=\chi(\cali;R)=\chi^{(2)}(\cali)$ by Theorem~\ref{the:coincidence}.
\end{remark}

\begin{remark}[Bar construction of classifying $\cali$-space] \label{rem:Ebarcalc}
There exists a functorial construction $E^{\barcon}\cali$ of
$E\cali$ by a kind of bar construction. Namely, the contravariant
functor $E^{\barcon}\cali \colon \cali \to \SPACES$ sends an object
$i$ to the space  $B^{\barcon}(i \downarrow \cali)$, which is the
geometric realization of the nerve of the category of objects under
$i$ (see Davis--L\"uck~\cite[page~230]{Davis-Lueck(1998)} and
also~Bousfield--Kan~\cite[page~327]{Bousfield-Kan(1972)}). An
equivalent definition of the bar construction in terms of the tensor
product in Definition~\ref{def:tensor_product_for_I-spaces} is
\begin{equation} \label{equ:bar_construction_as_tensor}
E^{\barcon}\cali=\{*\} \otimes_\cali B^{\barcon} (?\downarrow \cali
\downarrow ??),
\end{equation}
from which we prove that $E^{\barcon}\cali$ is an
$\cali$-$CW$-complex. The $\cali \times \cali^{\op}$-space
$B^{\barcon} (?\downarrow \cali \downarrow ??)$ is an $\cali \times
\cali^{\op}$-$CW$-complex (see~\cite[page~228]{Davis-Lueck(1998)}).
For each path
\begin{equation*}
i_0 \to i_1 \to i_2 \to \cdots \to i_n
\end{equation*}
of $n$-many non-identity morphisms in $\cali$, $B^{\barcon}
(?\downarrow \cali \downarrow ??)$ has an $n$-cell based at
$(i_0,i_n)$, that is a cell of the form $\mor_\cali(?,i_0) \times
\mor_\cali(i_n,??) \times D^n$. By \cite[Lemma~3.19
(2)]{Davis-Lueck(1998)}, the tensor product $E^{\barcon}\cali$ in
\eqref{equ:bar_construction_as_tensor} is an $\cali$-$CW$-complex:
an $(m+n)$-cell based at $i$ is an $n$-cell of $B^{\barcon}
(?\downarrow \cali \downarrow ??)$ based at $(i,j)$ and an $m$-cell
of the $CW$-complex $*(j)$ (see~\cite[page~229]{Davis-Lueck(1998)}).
More explicitly, for each path of $n$-many non-identity morphisms
\begin{equation} \label{equ:sequence_of_I_morphisms}
i_0 \to i_1 \to i_2 \to \cdots \to i_n
\end{equation}
the $\cali$-$CW$-complex $E^{\barcon}\cali$ has an $n$-cell based at
$i_0$.

Though the bar construction is in general not a finite
$\cali$-$CW$-complex, it is in certain cases. For example, if
$\cali$ has only finitely many morphisms, no nontrivial
isomorphisms, and no nontrivial endomorphisms, then there are only
finitely many paths as in \eqref{equ:sequence_of_I_morphisms}, and
hence only finitely many $\cali$-cells in $E^{\barcon}\cali$.

The bar construction is also compatible with induction. Given a
functor $\alpha \colon \cali \to \cald$, we obtain a map of
$\cald$-spaces
$$E^{\barcon}\alpha \colon \alpha_*E^{\barcon}\cali \to E^{\barcon}\cald,$$
where $\alpha_*$ denotes induction with the functor $\alpha$
(see~\cite[Definition~1.8]{Davis-Lueck(1998)}). If $T \colon \alpha
\to \beta$ is a natural transformation of functors $\cali \to
\cald$, we obtain for any $\cali$-space $X$ a natural transformation
$T_* \colon \alpha_*X \to \beta_*X$ which comes from the map of
$\cali$-$\cald$-spaces $\mor_{\cald}(??,\alpha(?)) \to
\mor_{\cald}(??,\beta(?))$ sending $g \colon {??} \to \alpha(?)$ to
$T(?) \circ g \colon {??} \to \beta(?)$.
\end{remark}

\begin{lemma}[Invariance of finite models under equivalence of categories] \label{lem:finite_models_and_equivalences_of_categories}
Suppose $\cali$ and $\calj$ are equivalent categories. Then
$\cali$ admits a finite $\cali$-$CW$-model for $E\cali$ if and only if $\calj$
admits a finite $\calj$-$CW$-model for $E\calj$. More precisely, if $F\colon \cali \to \calj$
is an equivalence of categories and $Y$ is a finite $\calj$-$CW$-model for $E\calj$, then the restriction
$\res_F Y$ is a finite $\cali$-$CW$-model for $E\cali$.
\end{lemma}
\begin{proof}
For any functor $F\colon \cali \to \calj$, we have an adjunction
\[
\ind_F\colon \cali\text{-}\SPACES \rightleftarrows \calj\text{-}\SPACES \colon \res_F
\]
defined by
\[
\ind_F(X):=X(?)\otimes_\cali \mor_\calj\bigl(??,F(?)\bigr) \;\;\;\;\; \res_F(Y):=Y \circ F(?).
\]
The $\cali$-space $\res_F(Y)$ is naturally homeomorphic to $Y(?)\otimes_\calj \mor_\calj\bigl(F(??),?\bigr)$. But since we are assuming $F$ is an
equivalence of categories, it is a left adjoint in an adjoint equivalence $(F,G)$, and we have natural homeomorphisms of $\cali$-spaces
\begin{eqnarray*}
\res_F(Y) & \cong &  Y(?)\otimes_\calj \mor_\calj\bigl(F(??),?\bigr) \\
& \cong & Y(?)\otimes_\calj \mor_\calj\bigl(??,G(?)\bigr) \\
& \cong & \ind_G(Y).
\end{eqnarray*}
Since $\ind_G$ is a left adjoint, so is $\res_F$, and $\res_F$ therefore preserves pushouts. Note also
\[ \res_F \mor_\calj(?,j) = \mor_\calj\bigl(F(?),j\bigr) \cong  \mor_\cali\bigl(?,G(j)\bigr).
\]

If $Y$ is a finite $\calj$-$CW$-model for $E\calj$ with $n$-skeleton
$$
  \comsquare{\coprod_{\lambda \in \Lambda_n} \mor_{\calj}(-,j_\lambda) \times S^{n-1}}
  {}{Y_{n-1}}{}{} {\coprod_{\lambda \in \Lambda_n} \mor_{\calj}(-,j_\lambda) \times D^n}
  {}{Y_n,}
  $$
then $X:=\res_F Y$ is a finite $\cali$-$CW$-complex with $n$-skeleton
$$
  \comsquare{\coprod_{\lambda \in \Lambda_n} \mor_{\cali}\bigl(-,G(j_\lambda)\bigr) \times S^{n-1}}
  {}{X_{n-1}}{}{} {\coprod_{\lambda \in \Lambda_n} \mor_{\cali}\bigl(-,G(j_\lambda)\bigr) \times D^n}
  {}{X_n.}
  $$
Clearly, $\res_F Y$ is contractible at each object $i$, since $\res_F Y (i)=Y(F(i)) \simeq \ast$.
\end{proof}


\section{Homotopy Colimits of Categories}
\label{sec:homotopy_colimits_of_categories}

\begin{definition}[Homotopy colimit for categories]
\label{def:homotopy_colimit}
  Let $\calc \colon \cali \to \CATS$ be a covariant functor from some
  (small) index
  category $\cali$ to the category of small categories. Its
  \emph{homotopy colimit}
  $$\hocolim_{\cali} \calc$$
  is the following
  category. Objects are pairs $(i,c)$, where $i \in \ob(\cali)$ and $c
  \in \ob\bigl(\calc(i)\bigr)$.  A morphism from $(i,c)$ to $(j,d)$ is
  a pair $(u,f)$, where $u \colon i \to j$ is a morphism in $\cali$
  and $f \colon \calc(u)(c) \to d$ is a morphism in $\calc(j)$.  The
  composition of the morphisms $(u,f) \colon (i,c) \to (j,d)$ and
  $(v,g) \colon (j,d) \to (k,e)$ is the morphism
  $$(v,g) \circ (u,f) = (v \circ u, g \circ \calc(v)(f)) \colon (i,c) \to (k,e).$$
  The identity of $(i,c)$ is given by $(\id_i,\id_c)$.

  This homotopy colimit construction for functors is often called the \emph{Grothendieck construction} or the \emph{category of elements}.

\end{definition}

In which sense is $\hocolim_{\cali} \calc$ a homotopy colimit? First, recall from \cite{Illusie(1972)} that the nerve functor induces an equivalence of categories $\text{Ho}\; \CATS \to \text{Ho}\; \SSETS$, where $\text{Ho}\; \CATS$ denotes the localization of $\CATS$ with respect to nerve weak equivalences and $\text{Ho}\; \SSETS$ denotes the localization of $\SSETS$ with respect to the usual weak equivalences. In \cite{Thomason(1979)}, Thomason proved that $\hocolim_{\cali} \calc$ in $\CATS$ corresponds to the Bousfield--Kan construction in $\SSETS$ under this equivalence of categories. Consequently, $\hocolim_{\cali} \calc$ has a universal property in the form of a bijection
\begin{equation} \label{equ:hocolim_explanation}
\text{Ho}\; \CATS (\hocolim_{\cali} \calc, \Gamma) \cong \text{Ho} \;\CATS^\cali (\calc, \underline{\Gamma}),
\end{equation}
for any category $\Gamma$. Here $\underline{\Gamma}$ indicates the $\cali$-diagram that is constant $\Gamma$.
In \cite{Thomason(1980)}, Thomason proved that $\CATS$ admits a cofibrantly generated model structure in which the weak equivalences are the nerve weak equivalences, so that the associated projective model structure on $\CATS^\cali$ exists. The model-theoretic construction of a homotopy colimit of the $\cali$-diagram $\calc$ in $\CATS$ as a colimit of a cofibrant replacement of $\calc$ in the projective model structure therefore works.  This model-theoretic construction also has the universal property in \eqref{equ:hocolim_explanation}, so is isomorphic to $\hocolim_{\cali} \calc$ in $\text{Ho}\; \CATS$, i.e. weakly equivalent to $\hocolim_{\cali} \calc$ in $\CATS$. A direct proof that $\hocolim_{\cali} \calc$ satisfies the universal property \eqref{equ:hocolim_explanation} is in Grothendieck's letter \cite{Grothendieck(1983)}, see the article of Maltsiniotis\footnote{We thank George Maltsiniotis for clarifying these points about homotopy colimits in $\CATS$.} \cite[Section 3.1]{Maltsiniotis(2005)}.

\begin{remark} \label{rem:homotopy_colimit_of_pseudo_functor}
If $\calc$ is merely a pseudo functor, then it of course still has a
homotopy colimit. A \emph{pseudo functor} $\calc \colon \cali \to
\CATS$ is like an ordinary functor, but only preserves composition
and unit up to specified coherent natural isomorphisms
$\calc_{v,u}\colon \calc(v) \circ \calc(u) \Rightarrow \calc(v \circ
u)$ and $\calc_i\colon 1_{\calc(i)} \Rightarrow \calc(\id_i)$.
Moreover, $\calc_{v,u}$ is required to be natural in $v$ and $u$.
The objects and morphisms of the \emph{homotopy colimit}
$\hocolim_{\cali} \calc$ are defined as in the strict case of
Definition~\ref{def:homotopy_colimit}. The composition in
$\hocolim_{\cali} \calc$ is defined by the modified rule
$$(v,g) \circ (u,f) = (v \circ u, g \circ (\calc(v)(f))\circ \calc_{v,u}^{-1}(c))$$
while the identity of the object $(i,c)$ is given by
$$(\id_i,\calc_i^{-1}(c)).$$
The homotopy colimit of a pseudo functor $\calc \colon \cali \to
\CATS$ is an ordinary 1-category with strictly associative and
strictly unital composition.
\end{remark}

\begin{remark} \label{rem:homotopy_colimit_is_2-functor}
For a fixed category $\cali$, the homotopy colimit construction
$\hocolim_{\cali}-$ is a strict 2-functor from the strict 2-category
of pseudo functors $\cali \to \CATS$, pseudo natural
transformations, and modifications into the strict 2-category
$\CATS$.
\end{remark}

\begin{example}[Homotopy colimit of a constant functor] \label{exa:hocolim_constant}
If $\calc \colon \cali \to \CATS$ is a constant functor, say
constantly a category also called $\calc$, then $\hocolim_{\cali}
\calc= \cali \times \calc$.
\end{example}

\begin{example}[Homotopy colimit for $\cali$ with a terminal object]
\label{exa:homotopy_hocolimit_and_when_I_has_terminal_object}
Suppose $\cali$ has a terminal object $t$ and $\calc \colon \cali
\to \CATS$ is a strict covariant functor. Then $\hocolim_\cali
\calc$ is homotopy equivalent to $\calc(t)$ as follows. This is
analogous to the familiar fact that $\calc(t)$ is a colimit of
$\calc$. The components of the universal cocone
\begin{equation}\label{equ:cocone_for_I_with_terminal_object}
\pi \colon \calc \Rightarrow \Delta_{\calc(t)}
\end{equation}
are $\calc(i \to t)$. Applying $\hocolim_\cali-$ to
\eqref{equ:cocone_for_I_with_terminal_object} and composing with the projection
gives us a functor $F$ $$\xymatrix@R=1pc@C=4pc{\hocolim_\cali \calc
  \ar[r]_-{\hocolim_\cali \pi} \ar@/^1.5pc/[rr]^F & \cali \times \calc(t)
  \ar[r]_-{\pr_{\calc(t)}} & \calc(t) \\ (i,c) \ar@{|->}[rr] & & \calc(i \to
  t)(c).}$$ The functor $G \colon \calc(t) \to \hocolim_\cali \calc$,
$G(c)=(t,c)$ is a homotopy inverse, since $F\circ G =\id_{\calc(t)}$ and we have
a natural transformation $\id_{\hocolim_\cali \calc} \Rightarrow G\circ F$ with
components
$$\xymatrix{(i \to t, \id_{\calc(i \to t)}) \colon (i,c) \ar[r] & (t,\calc(i \to t) c).}$$
\end{example}

Let $\calh$ denote the homotopy colimit of the $\cali$-diagram of
categories $\calc$. We now construct an $\cali$-diagram of
$\calh$-spaces $E^\calh$ with the property that its tensor product
with $E \cali$ is $\calh$-homotopy equivalent to a classifying
$\calh$-space for $\calh$. This $\cali$-diagram of $\calh$-spaces
$E^\calh$ will play an important role in the inductive proof of the
Homotopy Colimit Formula Theorem~\ref{the:homotopy_colimit_formula}.

\begin{construction}[Construction of $E^\calh$]
Let $\calc \colon \cali \to \CATS$ be a strict covariant functor,
and abbreviate $\calh = \hocolim_{\cali} \calc$. Define a functor
\begin{eqnarray}
  & E^{\calh} \colon \cali \to \calh\text{-}\SPACES &
  \label{Ecalh}
\end{eqnarray}
as follows. Given an object $i \in \cali$, we have the functor
\begin{eqnarray}
  & \alpha(i) \colon \calc(i) \to \calh &
  \label{alpha(i)}
\end{eqnarray}
sending an object $c$ to the object $(i,c)$ and a morphism $f \colon
c \to d$ to the morphism $(\id_i,f)$. We define
$$E^{\calh}(i) = \alpha(i)_* E^{\barcon}\bigl(\calc(i)\bigr).$$
Consider a morphism $u \colon i \to j$ in $\cali$. It induces a
natural transformation $T(u) \colon \alpha(i) \to \alpha(j) \circ
\calc(u)$ from the functor $\alpha(i) \colon \calc(i) \to \calh$ to
the functor $\alpha(j) \circ \calc(u) \colon \calc(i) \to \calh$ by assigning
to an object $c$ in $\calc(i)$ the morphism
$$(u,\id_{\calc(u)(c)}) \colon \alpha(i)(c) = (i,c)
\to \alpha(j)\circ \calc(u)(c) = (j,\calc(u)(c)).$$ From
Remark~\ref{rem:Ebarcalc} we obtain a map of $\calh$-spaces
$$T(u)_* \colon\alpha(i)_*E^{\barcon}\bigl(\calc(i)\bigr)
\to \alpha(j)_*\calc(u)_*E^{\barcon}\bigl(\calc(i)\bigr)$$ and a
map of $\calc(j)$-spaces
$$E^{\barcon}\bigl(\calc(u)\bigr) \colon \calc(u)_*E^{\barcon}\bigl(\calc(i)\bigr) \to E^{\barcon}\bigl(\calc(j)\bigr).$$
Finally, for the morphism $u$ in $\cali$, we define $E^{\calh}(u)
\colon E^{\calh}(i) \to E^{\calh}(j)$ by the composite
$$\alpha(i)_*E^{\barcon}\bigl(\calc(i)\bigr)
\xrightarrow{T(u)_*} \alpha(j)_*\calc(u)_*E^{\barcon}\bigl(\calc(i)\bigr)
\xrightarrow{\alpha(j)_*(E^{\barcon}\bigl(\calc(u)\bigr)}
\alpha(j)_*E^{\barcon}\bigl(\calc(j)\bigr).$$
\end{construction}

Define the
homotopy colimit of the covariant functor $E^{\calh}$ of~\eqref
{Ecalh} to be the contravariant $\calh$-space
\begin{eqnarray}
  \hocolim_{\cali} E^{\calh} & := & (i,c) \mapsto E\cali \otimes_{\cali} \left(E^{\calh}(i,c)\right).
  \label{hocolim_cali_ebarcon}
\end{eqnarray}

\begin{lemma} \label{lem:model_for_E(hocolim_cali_calc)}
Consider any model $E\cali$ for the classifying $\cali$-space of the category
$\cali$. Then the  contravariant $\calh$-space $E\cali \otimes_{\cali} E^{\calh}$
of~\eqref{hocolim_cali_ebarcon} is $\calh$-homotopy equivalent to
the classifying $\calh$-space $E\calh$ of the category $\calh :=\hocolim_{\cali} \calc$.
\end{lemma}
\begin{proof}
We first show that for any object $(i,c)$ in $\calh$ the space
$E\cali \otimes_{\cali} \left( E^{\calh}(i,c) \right)$ is contractible. The covariant
functor $ E^{\calh}(i,c) \colon \cali \to \SPACES$ sends an object $j$
to
\begin{eqnarray*}
\alpha(j)_*\left(E^{\barcon} \calc (j) \right)(i,c) & = & \alpha(j)_*\left(E^{\barcon}\calc (j)\right) (?) \otimes_{\calh} \mor_{\calh}\bigl((i,c),?)\bigr)
\\
& = &
\left(E^{\barcon}\calc (j)\right)(?) \otimes_{\calc (j)} \mor_{\calh}\bigl((i,c),(j,?)\bigr)
\\
& = &
\left(E^{\barcon}\calc (j)\right)(?) \otimes_{\calc (j)}
\left(\coprod_{u \in \mor_{\cali}(i,j)} \mor_{\calc(j)}\bigl(\calc(u)(c),?\bigr)\right)
\\
& = &
\coprod_{u \in \mor_{\cali}(i,j)}
\left(E^{\barcon}\calc (j)\right)(?) \otimes_{\calc (j)}\mor_{\calc(j)}\bigl(\calc(u)(c),?\bigr)
\\
& = &
\coprod_{u \in \mor_{\cali}(i,j)} \left(E^{\barcon}\calc(j)\right)\bigl(\calc(u)(c)\bigr).
\end{eqnarray*}
Since  $\left(E^{\barcon}\calc(j)\right)\bigl(\calc(u)(c)\bigr)$ is contractible, the
projection
$$\coprod_{u \in \mor_{\cali}(i,j)} \left(E^{\barcon}\calc(j)\right)\bigl(\calc(u)(c)\bigr)
\to  \mor_{\cali}(i,j)$$ is a homotopy equivalence. Hence the
collection of these projections for $j \in \ob(\cali)$ induces a map
of $\cali$-spaces
$$\pr \colon E^{\calh}(i,c) \to \mor_{\cali}(i,?)$$
whose evaluation at each object $j$ in $\ob(\cali)$ is a homotopy
equivalence. We conclude
from~Davis--L\"uck~\cite[Theorem~3.11]{Davis-Lueck(1998)} that
$$E\cali \otimes_{\cali} \pr \colon E\cali \otimes_{\cali} E^{\calh}(i,c)
\xrightarrow{\simeq} E\cali \otimes_{\cali} \mor_{\cali}(i,?).$$ is
a homotopy equivalence. Since $E\cali \otimes_{\cali}
\mor_{\cali}(i,?) = E\cali(i)$ is contractible, this implies that
for any object $(i,c)$ in $\calh$ the space $E\cali \otimes_{\cali}
\left( E^{\calh}(i,c) \right)$ is contractible, as we initially claimed.

It remains to show that $E\cali \otimes_{\cali} E^{\calh}$ has the
$\calh$-homotopy type of an $\calh$-$CW$-complex. It is actually an
$\calh$-$CW$-complex. The following argument, that $E\cali \otimes_{\cali} E^{\calh}$ has the homotopy type of
an $\calh$-$CW$-complex, will be used again later.\footnote{This is a well-known standard argument, which
we present only so that the reader easily sees that it works in the setting of $\calh$-spaces.}

We have a filtration of $E\cali$
$$\emptyset = E\cali_{-1} \subseteq E\cali_{0} \subseteq E\cali_{1}
\subseteq \ldots \subseteq E\cali_{n} \subseteq \ldots \subseteq
E\cali = \bigcup_{n \geq 0}  E\cali_{n}$$
such that
$$E\cali =\colim_{n \to \infty}  E\cali_{n}$$
and for every
$n \geq 0$ there exists a pushout of $\cali$-spaces
\begin{eqnarray}
&
  \comsquare{\coprod_{\lambda \in \Lambda_n} \mor_{\cali}(-,i_\lambda) \times S^{n-1}}
  {}
  {E\cali_{n-1}}
  {}{}
  {\coprod_{\lambda \in \Lambda_n} \mor_{\cali}(-,i_\lambda) \times D^n}
  {}
  {E\cali_n.}
&
\label{cali_CW_decomposition_of_Ecali}
\end{eqnarray}
Since $- \otimes_{\cali} Z$ has a right adjoint~\cite[Lemma~1.9]{Davis-Lueck(1998)}
we get
$$E\cali \otimes_{\cali} E^{\calh}
=  \colim_{n \to \infty} E\cali_n \otimes_{\cali} E^{\calh}$$
as a colimit of $\calh$-spaces. After an application of $-\otimes_\cali E^{\calh}$ to \eqref{cali_CW_decomposition_of_Ecali}, we obtain pushouts of $\calh$-spaces
\begin{eqnarray}
&
  \comsquare{\coprod_{\lambda \in \Lambda_n} E^{\calh}(i_\lambda) \times S^{n-1}}
  {f_{n-1}}
  {E\cali_{n-1} \otimes_{\cali} E^{\calh}}
  {}{}
  {\coprod_{\lambda \in \Lambda_n} E^{\calh}(i_\lambda) \times D^n}
  {}{E\cali_{n} \otimes_{\cali} E^{\calh}}
&
\label{pushout_of_calh-spaces}
\end{eqnarray}
where the left vertical arrow and hence the right vertical arrow are
cofibrations of $\calh$-spaces. By induction we may assume that
$E\cali_{n-1} \otimes_{\cali} E^{\calh}$ has the homotopy type of an
$\calh$-$CW$-complex. Since the vertical maps are cofibrations, by
replacing it with a homotopy equivalent $\calh$-$CW$-complex we do
not change the homotopy type of the pushout (the usual proof for
spaces goes through for $\calh$-spaces). Hence we may assume that
$E\cali_{n-1} \otimes_{\cali} E^{\calh}$ is a $\calh$-$CW$-complex.
We may also assume that $f_{n-1}$ is cellular: since the vertical
maps are cofibrations, by replacing $f_{n-1}$ by a homotopic
cellular map, which exists
by~Davis--L\"uck~\cite[cf.~Theorem~3.7]{Davis-Lueck(1998)}, we also
do not change the homotopy type of the pushout.
See~Selick~\cite[Theorem~7.1.8]{Selick(1997)} for a proof of this
statement for spaces which translates verbatim to the setting of
$\calh$-spaces. If $f_{n-1}$ is cellular,
diagram~\eqref{pushout_of_calh-spaces} is a cellular pushout. Hence
we completed the induction step, showing that $E\cali_{n}
\otimes_{\cali} E^{\calh}$ has the homotopy type of an
$\calh$-$CW$-complex.

It remains to show that $E\cali \otimes_{\cali} E^{\calh}$ has the homotopy
type of a $\calh$-$CW$-complex: choose $\calh$-$CW$-complexes $Z_n$ and
$\calh$-homotopy equivalences
$g_n:Z_n\to E\cali_{n} \otimes_{\cali} E^{\calh}$.
By iteratively
replacing $Z_n$ by the mapping cylinder of
\[
    Z_{n-1}\xrightarrow{g_{n-1}} E\cali_{n-1} \otimes_{\cali} E^{\calh}\to E\cali_{n} \otimes_{\cali} E^{\calh}\xrightarrow{\bar g_n}Z_n,
\]
where $\bar g_n$ is a homotopy inverse of $g_n$, one finds
a new sequence of
homotopy equivalences $g_n':Z_n\to E\cali_{n} \otimes_{\cali} E^{\calh}$
(with the modified $\calh$-$CW$-complexes $Z_n$)
such that $g_n'\vert_{Z_{n-1}}=g_{n-1}'$.
\end{proof}


\typeout{---------------  Section: Homotopy Colimit Formula --------------------}
\section{Homotopy Colimit Formula for Finiteness Obstructions and Euler Characteristics}
\label{sec:Homotopy_colimit_formula}

In this section we prove the main theorem of this paper: the
Homotopy Colimit Formula. It expresses the finiteness obstruction,
the Euler characteristic, and the $L^2$-Euler characteristic of the
homotopy colimit of a diagram in $\CATS$ in terms of the respective
invariants for the diagram entries at the base objects for cells in
a finite model for the $\cali$-classifying space of $\cali$.
Analogous formulas for the functorial counterparts of the Euler
characteristic and $L^2$-Euler characteristic are included. The
Homotopy Colimit Formula is initially stated and proved for strict
functors $\calc\colon \cali \to \CATS$, but we prove that it also
holds for pseudo functors $\cald\colon \cali \to \CATS$ in
Corollary~\ref{cor:homotopy_colimit_formula_for_pseudo_functors}.
The full generality of pseudo functors is needed for the
applications to complexes of groups in
Section~\ref{sec:complexes_of_groups}.

\subsection{Homotopy Colimit Formula}

\begin{theorem}[Homotopy Colimit Formula]
  \label{the:homotopy_colimit_formula}  Let
  $\cali$ be a small category such that there exists a finite $\cali$-$CW$-model for
  its classifying $\cali$-space. Fix such a  finite $\cali$-$CW$-model $E\cali$.
  Denote by $\Lambda_n$ the finite set of $n$-cells $\lambda =
  \mor_\cali(?,i_\lambda) \times D^n$ of $E\cali$. Let $\calc \colon \cali \to \CATS$ be a
  covariant functor.  Abbreviate $\calh = \hocolim_{\cali}
  \calc$. Then:

  \begin{enumerate}

  \item \label{the:homotopy_colimit_formula:directly_finite}
    If $\cali$ is directly finite, and $\calc(i)$ is directly finite for every object $i \in
    \ob(\cali)$, then the category $\calh$ is directly finite;

  \item \label{the:homotopy_colimit_formula:(EI)} If
    $\cali$ is an EI-category,
    $\calc(i)$ is an EI-category for every object $i \in
    \ob(\cali)$,
    and for every automorphism $u \colon i \xrightarrow{\cong} i$ the
    map $\iso(\calc(i)) \to \iso(\calc(i)), \; \overline{x} \mapsto
    \overline{\calc(u)(x)}$ is the identity, then the category $\calh$ is an EI-category;

  \item \label{the:homotopy_colimit_formula:(FP)} If for
    every object $i$ the category $\calc(i)$ is of type (FP$_R$), then
    the category $\hocolim_{\cali} \calc$ is of type (FP$_R$);

  \item \label{the:homotopy_colimit_formula:(FF)} If for
    every object $i$ the category $\calc(i)$ is of type (FF$_R$), then
    the category $\hocolim_{\cali} \calc$ is of type (FF$_R$);

  \item \label{the:homotopy_colimit_formula:o} If for every
    object $i$ the category $\calc(i)$ is of type (FP$_R$),  then we obtain
    for the finiteness obstruction
    $$o(\calh;R) =  \sum_{n \geq 0} (-1)^n \cdot
    \sum_{\lambda \in \Lambda_n} \alpha(i_\lambda)_*(o(\calc(i_\lambda);R)),$$
    where $\alpha(i_\lambda)_* \colon K_0(R\calc(i_\lambda)) \to K_0(R\calh)$ is the
    homomorphism induced by the canonical functor
    $\alpha(i_\lambda) \colon \calc(i_\lambda) \to \calh$ defined in~\eqref{alpha(i)};

\item \label{the:homotopy_colimit_formula:chi} Suppose that
    $\cali$ is directly finite and
    $\calc(i)$ is directly finite for every object $i \in
    \ob(\cali)$. If for every
    object $i$ the category $\calc(i)$ is additionally of type (FP$_R$) then we obtain for the functorial Euler characteristic
  $$\chi_f(\calh;R) =
  \sum_{n \geq 0} (-1)^n \cdot \sum_{\lambda \in \Lambda_n} \alpha(i_\lambda)_*(\chi_f(\calc(i_\lambda);R)),$$
  where $\alpha(i_\lambda)_* \colon U(\calc(i_\lambda)) \to U(\calh)$ is the
  homomorphism induced by the canonical functor
  $\alpha(i_\lambda) \colon \calc(i_\lambda) \to \calh$ defined
  in~\eqref{alpha(i)}. Summing up, we also have
$$\chi\bigl(\calh;R\bigr) = \sum_{n \geq 0} (-1)^n \cdot \sum_{\lambda \in \Lambda_n} \chi(\calc(i_\lambda);R).$$  If $R$ is Noetherian, in addition to the direct finiteness and (FP$_R$) hypotheses, we obtain for the Euler characteristics of the classifying spaces
  $$\chi\bigl(B\calh;R\bigr) = \sum_{n \geq 0} (-1)^n \cdot \sum_{\lambda \in \Lambda_n} \chi(B\calc(i_\lambda);R);$$

\item \label{the:homotopy_colimit_formula:chi(2)} Suppose that $\cali$ is
  directly finite and $\calc(i)$ is directly finite for every object $i \in \ob(\cali)$.  If for every object $i$ the category
  $\calc(i)$ is additionally of type ($L^2$), then $\calh$ is of type ($L^2$) and we obtain for
  the functorial $L^2$-Euler characteristic
  $$\chi_f^{(2)}(\calh) =
  \sum_{n \geq 0} (-1)^n \cdot \sum_{\lambda \in \Lambda_n}
  \alpha(i_\lambda)_*\bigl(\chi_f^{(2)}(\calc(i_\lambda))\bigr),$$ where $\alpha(i_\lambda)_*
  \colon U^{(1)}(\calc(i_\lambda)) \to U^{(1)}(\calh)$ is the homomorphism induced by
  the canonical functor $\alpha(i_\lambda) \colon \calc(i_\lambda) \to \calh$ defined
  in~\eqref{alpha(i)}, and we obtain for the $L^2$-Euler characteristic
  $$\chi^{(2)}(\calh) =
  \sum_{n \geq 0} (-1)^n \cdot \sum_{\lambda \in \Lambda_n} \chi^{(2)}(\calc(i_\lambda)).$$

\end{enumerate}
\end{theorem}
\begin{proof}%
\ref{the:homotopy_colimit_formula:directly_finite} Consider
  morphisms $(u,f) \colon(i,c) \to (j,d)$ and $(v,g) \colon(j,d) \to
  (i,c)$ in $\calh$ with $(v,g) \circ (u,f) = \id_{(i,c)}$. This
  implies $vu = \id_i$ and $g \circ \calc(v)(f) = \id_c$. Since
  $\cali$ and $\calc(i)$ are by assumption directly finite, we
  conclude $uv = \id_j$ and $\calc(v)(f) \circ g =
  \id_{\calc(v)(d)}$. Hence
  \begin{multline*}
    (u,f) \circ (v,g) = \bigl(uv,f \circ \calc(u)(g)\bigr) = \bigl(uv,
    \calc(uv)(f) \circ \calc(u)(g)\bigr) =
    \bigl(uv,\calc(u)(\calc(v)(f) \circ g)\bigr)
    \\
    = \bigl(uv,\calc(u)(\id_{\calc(v)(d)})\bigr) =
    \bigl(\id_j,\id_{\calc(u)\bigl(\calc(v)(d)\bigr)}\bigr) =(\id_j,\id_d).
  \end{multline*}%
\ref{the:homotopy_colimit_formula:(EI)}
  Consider an endomorphism
  $(u,f) \colon (i,c) \to (i,c)$ in $\calh$.  Since $\cali$ is an EI-category,
  $u\colon i \to i$ is an automorphism. Since $\overline{\calc(u)(c)} = \overline{c}$
  by assumption, we can choose an isomorphism $g \colon c
  \xrightarrow{\cong} \calc(u)(c)$. Hence $fg$ is an
  endomorphism in $\calc(i)$. Since $\calc(i)$ is an EI-category,
  and $g$ is an isomorphism, $f$ is also an isomorphism. Since $u$ and $f$ are isomorphisms, $(u,f)$
  is an isomorphism.  \\[1mm]%
\ref{the:homotopy_colimit_formula:(FP)}
  and~\ref{the:homotopy_colimit_formula:o}.  We say that an
  $R\calh$-chain complex $C_*$ is of type (FP$_R$) if it admits a
  \emph{finite projective approximation}, i.e., there is a finite length
  chain complex $P_*$ of finitely generated, projective $R\calh$-modules together with an
  $R\calh$-chain map $f_*\colon P_* \to C_*$ such that $H_n(f_*(i,c))$
  is bijective for all $n \geq 0$ and $(i,c) \in \ob(\calh)$. If $C_*$
  is of type (FP$_R$), define its finiteness obstruction
  $$o(C_*) := \sum_{n \geq 0} (-1)^n \cdot [P_n] \quad \in K_0(R\calh)$$
  for any choice $P_*$ of finite projective approximation.  This is
  independent of the choice of $P_*$ and the basic properties of it
  were studied by L\"uck~\cite[Chapter~11]{Lueck(1989)}.  If
  $0[\underline{R}]$ is the $R\calh$-chain complex concentrated in
  dimension zero and given there by the constant $R\calh$-module
  $\underline{R}$, then $\calh$ is of type (FP$_R$) if and only if
  $0[\underline{R}]$ is of type (FP$_R$) and in this case
  $$o(\calh;R) = o(0[\underline{R}]) \in K_0(R\calh).$$

  Consider a finite $\cali$-$CW$-complex $X$.  We want to show by
  induction over the dimension of $X$ that the $R\calh$-chain complex
  $C_*(X \otimes_\cali E^{\calh})$ is of type (FP$_R$) and satisfies
  $$o\bigl(C_*(X \otimes_\cali E^{\calh})\bigr)
  = \sum_{n \geq 0} (-1)^n \cdot \sum_{\lambda \in \Lambda_n}
  \alpha(i_\lambda)_*(o(\calc(i_\lambda);R)),$$
  where $\Lambda_n$ denotes the set of
  $n$-cells of $X$ and $i_\lambda$ is the object at which the $n$-cell
  $\lambda = \mor_\cali(?,i_\lambda) \times D^n$ of $X$ is based.

  The induction beginning, where $X$ is the empty set, is obviously
  true.  The induction step is done as follows. Let $d$ be the
  dimension of $X$.  Then $X_d$ is obtained from $X_{d-1}$ by a
  pushout of $\cali$-spaces
  $$
  \xycomsquareminus{\coprod_{\lambda \in \Lambda_d} \mor_{\calc}(-,i_\lambda) \times
  S^{d-1}} {} {X_{d-1}} {}{} {\coprod_{\lambda \in \Lambda_d} \mor_{\calc}(-,i_\lambda)
  \times D^d} {} {X = X_d.}
  $$
  Applying $- \otimes_{\cali} E^{\calh}$ to it yields, because
  $E^{\calh}(i) = \alpha(i)_*E^{\barcon}\bigl(\calc(i)\bigr)$, a pushout of
  $\calh$-spaces with a cofibration as left vertical arrow
  $$\xymatrix{
  \coprod_{\lambda \in \Lambda_d} \alpha(i_\lambda)_*E^{\barcon}\bigl(\calc(i_\lambda)\bigr) \times S^{d-1} \ar[d]
  \ar[r] & X_{d-1} \otimes_{\cali} E^{\calh} \ar[d]
  \\
  \coprod_{\lambda \in \Lambda_d} \alpha(i_\lambda)_*E^{\barcon}\bigl(\calc(i_\lambda)\bigr) \times D^{d} \ar[r] & X
  \otimes_{\cali} E^{\calh}.}
  $$
  In the sequel we can assume without loss of generality that
  $X_{d-1} \otimes_{\cali} E^{\calh}$ and $X \otimes_{\cali} E^{\calh}$ are
  $\calh$-$CW$-complexes and the diagram above is a pushout of
  $\calh$-$CW$-complexes, since this can be arranged by replacing them
  by homotopy equivalent $\calh$-$CW$-complexes (see the proof of
Lemma~\ref{lem:model_for_E(hocolim_cali_calc)}).  We obtain an exact
  sequence of $R\calh$-chain complexes
  $$0 \to C_*(X_{d-1} \otimes_{\cali} E^{\calh}) \to C_*(X \otimes_{\cali} E^{\calh}) \to
  \bigoplus_{\lambda \in \Lambda_d} \Sigma^d C_*\bigl(\alpha(i_\lambda)_*E^{\barcon}\calc(i_\lambda)\bigr)
  \to 0.$$ Consider $\lambda \in \Lambda_d$.  Since $\calc(i_\lambda)$ is of type (FP$_R$),
  we can find a finite projective $R\calc(i_\lambda)$-chain complex $P_*$
  whose homology is concentrated in dimension zero and given there by
  the constant $R\calc(i_\lambda)$-module $\underline{R}$. Since
  $C_*(E^{\barcon}\calc(i_\lambda))$ is a projective $R\calc(i_\lambda)$-chain
  complex with the same homology, there is an $R\calc(i_\lambda)$-chain
  homotopy equivalence $f_* \colon P_* \xrightarrow{\simeq}
  C_*\bigl(E^{\barcon}\calc(i_\lambda)\bigr)$ (see~L\"uck~\cite[Lemma~11.3 on page~213]{Lueck(1989)} and
  $$o(\calc(i_\lambda);R) = o(P_*) = \sum_{n \geq 0} (-1)^n \cdot [P_n] \in K_0(R\calc(i_\lambda)).$$
  Obviously
  $$\alpha(i_\lambda)_*f_* \colon \alpha(i_\lambda)_*P_* \xrightarrow{\simeq}
  \alpha(i_\lambda)_*C_*\bigl(E^{\barcon}\bigl(\calc(i_\lambda)\bigr)\bigr) =
  C_*\bigl(\alpha(i_\lambda)_*E^{\barcon}\calc(i_\lambda)\bigr)$$ is an $R\calh$-chain
  homotopy equivalence. Hence
  $C_*(\alpha(i_\lambda)_*E^{\barcon}\calc(i_\lambda))$ and, by the induction
  hypothesis, $C_*(X_{d-1} \otimes_\cali E^{\calh})$ are $R\calh$-chain
  complexes of type (FP$_R$). We conclude
  from~L\"uck~\cite[Lemma~11.3 on page~213]{Lueck(1989)} that
  $C_*(X \otimes_\cali E^{\calh})$ is of type (FP$_R$) and
  $$o\bigl(C_*(X \otimes_\cali E^{\calh})\bigr)
  = o\bigl(C_*(X_{d-1} \otimes_\cali E^{\calh})\bigr) + \sum_{\lambda \in \Lambda_d}
  o\bigl(\Sigma^d \alpha(i_\lambda)_*C_*(E^{\barcon}\calc(i_\lambda))\bigr).$$
  This implies together with the induction hypothesis applied to $X_{d-1}$
  \begin{eqnarray*}
  \lefteqn{o\bigl(C_*(X \otimes_\cali E^{\calh})\bigr)}
  & &
  \\
  & = &
  \sum_{n = 0}^{d-1} (-1)^n \cdot \sum_{\lambda \in \Lambda_n} \alpha(i_\lambda)_*(o(\calc(i_\lambda);R))
  +
  \sum_{\lambda \in \Lambda_d} (-1)^d \cdot \alpha(i_\lambda)_*(o(\calc(i_\lambda);R))
  \\
  & = &
  \sum_{n = 0}^d (-1)^n \cdot \sum_{\lambda \in \Lambda_n} \alpha(i_\lambda)_*(o(\calc(i_\lambda);R)).
  \end{eqnarray*}
  This finishes the induction step.

  Assertions~\ref{the:homotopy_colimit_formula:(FP)}
  and~\ref{the:homotopy_colimit_formula:o} follow by taking $X =
  E\cali$.
\\[1mm]\ref{the:homotopy_colimit_formula:(FF)}
This proof is analogous to that of
assertion~\ref{the:homotopy_colimit_formula:(FP)}.
\\[1mm]\ref{the:homotopy_colimit_formula:chi}
By~\ref{the:homotopy_colimit_formula:directly_finite}
and~\ref{the:homotopy_colimit_formula:(FP)}, the category $\calh$ is
directly finite and of type (FP$_R$). Then an application of
$\rk_{R\calh}$ to the formula for $o(\calh;R)$
in~\ref{the:homotopy_colimit_formula:o} yields the formula for
$\chi_f(\calh;R)$ in~\ref{the:homotopy_colimit_formula:chi} by the
naturality of $\rk_{R-}$ with respect to the functors
$\alpha(i_\lambda)$ between directly finite categories, see
Fiore--L\"uck--Sauer~\cite[Lemma~4.9]{FioreLueckSauerFinObsAndEulCharOfCats(2009)}.

An application of the augmentation homomorphism $\epsilon \colon
U(\calh) \to \IZ$ to the formula for $\chi_f(\calh;R)$ yields the
formula for $\chi(\calh;R)$. We also use the naturality of the
augmentation homomorphism, that is, the commutativity of diagram
(4.5)~in~\cite{FioreLueckSauerFinObsAndEulCharOfCats(2009)} for
$F=\alpha(i_\lambda)$.

If $R$ is additionally Noetherian, then Theorem~\ref{the:chi_f_determines_chi} applies, and the Euler characteristics of the categories agree with the Euler characteristics of the classifying spaces. \\[1mm]%
\ref{the:homotopy_colimit_formula:chi(2)}
  The proofs for the functorial $L^2$-Euler characteristic and the $L^2$-Euler characteristic are
  somewhat more complicated since the property ($L^2$) is more general than
(FP$_R$), and the $L^2$-Euler
  characteristic comes from the finiteness obstruction
  only in the case (FP$_R$). The proofs are variations of the proofs for
  assertions~\ref{the:homotopy_colimit_formula:(FP)}
  and~\ref{the:homotopy_colimit_formula:o}. Instead of using~L\"uck~\cite[Lemma~11.3 on page~213]{Lueck(1989)},
  we now use the basic properties of $L^2$-Euler characteristics for chain complexes of
  modules over group von Neumann algebras
  \cite[Lemma~5.7]{FioreLueckSauerFinObsAndEulCharOfCats(2009)}.
  For example, we use
  \cite[Lemma~5.7~(iv)]{FioreLueckSauerFinObsAndEulCharOfCats(2009)},
  which says for any injective group homomorphism $i \colon H \to G$
  and $\caln(H)$-chain complex $C_*$, we have $\chi^{(2)}(C_*) = \chi^{(2)}(\ind_{i_*}
  C_*)$, provided the sum of the $L^2$-Betti numbers of $C_*$ is
  finite. The injectivity hypothesis is easily verified: for
  every object $i \in \ob(\cali)$ and object $x \in \calc(i)$ the
  functor $\alpha(i) \colon\calc(i) \to \calh$ clearly induces an injection
  $\aut_{\calc(i)}(x) \to \aut_{\calh}(i,x)$. This finishes the proof of
  Theorem~\ref{the:homotopy_colimit_formula}.
\end{proof}

\begin{corollary} \label{cor:homotopy_colimit_formula_for_pseudo_functors}
Theorem~\ref{the:homotopy_colimit_formula} on homotopy colimits
holds for pseudo functors $\cald \colon \cali \to \CATS$.
\end{corollary}
\begin{proof}
We first remark that the pseudo functor $\cald \colon \cali \to
\CATS$ is equivalent to a strict functor $\calc \colon \cali \to
\CATS$ in the following sense. As usual, we denote by
$\Hom(\cali,\CATS)$ the strict 2-category of pseudo functors $\cali
\to \CATS$, pseudo natural transformations between them, and
modifications. The pseudo functor $\cald$ is equivalent to a
strict functor $\calc$ \emph{as objects of the 2-category}
$\Hom(\cali,\CATS)$. For example, we may take $\calc$ to be the
strict functor
$$i \mapsto \mor_{\Hom(\cali,\CATS)}(\cali(i,-),\cald).$$

The equivalence between $\calc$ and $\cald$ in $\Hom(\cali,\CATS)$
has two useful consequences. Since
$$\hocolim_{\cali}\colon\Hom(\cali,\CATS) \to \CATS$$ is a strict
2-functor, it sends any equivalence between $\calc$ and $\cald$ to
an equivalence in $\CATS$ between the categories
$\hocolim_{\cali}\calc$ and $\hocolim_{\cali}\cald$. Another
consequence of the equivalence between $\calc$ and $\cald$ is that
for every $i \in \cali$, the categories $\calc(i)$ and $\cald(i)$
are equivalent. With these observations we reduce
Corollary~\ref{cor:homotopy_colimit_formula_for_pseudo_functors} to
Theorem~\ref{the:homotopy_colimit_formula}.
\\[1mm]~%
\ref{the:homotopy_colimit_formula:directly_finite} Suppose
$\cald(i)$ is directly finite for every $i\in\ob(\cali)$ and $\cali$
is directly finite. Since direct finiteness is preserved under
equivalence of categories by Fiore--L\"uck--Sauer
\cite[Lemma~3.2]{FioreLueckSauerFinObsAndEulCharOfCats(2009)}, and
$\calc(i)$ is equivalent to $\cald(i)$, we see that $\calc(i)$ is
directly finite for every $i\in\ob(\cali)$. Hence
$\hocolim_\cali\calc$ is directly finite by
Theorem~\ref{the:homotopy_colimit_formula}~\ref{the:homotopy_colimit_formula:directly_finite}.
Since $\hocolim_\cali\cald$ is equivalent to $\hocolim_\cali\calc$,
it is also directly finite, again by
\cite[Lemma~3.2]{FioreLueckSauerFinObsAndEulCharOfCats(2009)}.
\\[1mm]~%
\ref{the:homotopy_colimit_formula:(EI)} Suppose that $\cali$ is an
EI-category, $\cald(i)$ is an EI-category for every $i \in
\ob(\cali)$, and for every automorphism $u \colon i
\xrightarrow{\cong} i$ the map $\iso(\cald(i)) \to \iso(\cald(i)),
\; \overline{y} \mapsto \overline{\cald(u)(y)}$ is the identity.
Since EI is preserved under equivalence of categories
\cite[Lemma~3.11]{FioreLueckSauerFinObsAndEulCharOfCats(2009)},
and $\calc(i)$ is equivalent to $\cald(i)$, we see $\cald(i)$ is an
EI-category.

We claim that for each automorphism $u$, the functor $\calc(u)$ also
induces the identity on isomorphism classes of objects of
$\calc(i)$. Let $\phi\colon \cald \to \calc$ be a pseudo
equivalence, that is, an equivalence in the 2-category
$\Hom(\cali,\CATS)$. For $x \in \calc(i)$, there is a $y \in
\cald(i)$ and an isomorphism $x \cong \phi_i(y)$. We have
isomorphisms
$$\calc(u)(x)\cong\calc(u) \phi_i(y) \cong \phi_i \cald(u)(y)\cong \phi_i(y) \cong
x,$$ and $\calc(u)$ induces the identity on isomorphism classes of
objects of $\calc(i)$. Then $\hocolim_\cali\calc$ is an EI-category
by
Theorem~\ref{the:homotopy_colimit_formula}~\ref{the:homotopy_colimit_formula:(EI)},
and so is $\hocolim_\cali\cald$, again by
\cite[Lemma~3.11]{FioreLueckSauerFinObsAndEulCharOfCats(2009)}.
\\[1mm]~%
\ref{the:homotopy_colimit_formula:(FP)}
and~\ref{the:homotopy_colimit_formula:(FF)} similarly follow from
Theorem~\ref{the:homotopy_colimit_formula}~\ref{the:homotopy_colimit_formula:(FP)}
and~\ref{the:homotopy_colimit_formula:(FF)}, since property
(FP$_R$), property (FF$_R$), and the finiteness obstruction are all
invariant under equivalence of categories
\cite[Theorem~2.8]{FioreLueckSauerFinObsAndEulCharOfCats(2009)}.
\\[1mm]~%
\ref{the:homotopy_colimit_formula:o} Suppose $\cald(i)$ is of type
(FP$_R$) for every $i\in\ob(\cali)$. Then every $\calc(i)$ is also
of type (FP$_R$), since property (FP$_R$) is invariant under
equivalence of categories
\cite[Theorem~2.8]{FioreLueckSauerFinObsAndEulCharOfCats(2009)}. As
in~\eqref{alpha(i)}, we have for each $i \in \cali$ the functor
$$\alpha^\cald(i) \colon \cald(i) \to \hocolim_\cali \cald $$
which sends an object $d$ to the object $(i,d)$ and a morphism $g
\colon d \to d'$ to the morphism $(\id_i,g \circ \cald^{-1}_i(d))$.
From a pseudo equivalence $\psi\colon \calc \to \cald$ we obtain a
strictly commutative diagram
\begin{equation} \label{equ:naturality_for_alpha_to_hocolim}
\xymatrix@C=3pc{\calc(i) \ar[d]_{\psi_i} \ar[r]^-{\alpha^\calc(i)} & \hocolim_\cali
\calc \ar[d]^{\hocolim_\cali \psi}
\\ \cald(i) \ar[r]_-{\alpha^\cald(i)} & \hocolim_\cali \cald}
\end{equation}
for each $i \in \ob(\cali)$. Since the finiteness obstruction is
invariant under equivalence of categories
\cite[Theorem~2.8]{FioreLueckSauerFinObsAndEulCharOfCats(2009)},
we may use
Theorem~\ref{the:homotopy_colimit_formula}~\ref{the:homotopy_colimit_formula:o}
for $\calc$ to obtain
$$\aligned
o(\hocolim_\cali\cald;R) &=
(\hocolim_\cali\psi)_\ast(o(\hocolim_\cali\calc;R)) \\
&=(\hocolim_\cali\psi)_\ast\left(  \sum_{n \geq 0} (-1)^n \cdot
\sum_{\lambda \in \Lambda_n} \alpha^\calc(i_\lambda)_*(o(\calc(i_\lambda);R))\right) \\
&= \sum_{n \geq 0} (-1)^n \cdot \sum_{\lambda \in \Lambda_n}
(\hocolim_\cali\psi)_\ast\circ \alpha^\calc(i_\lambda)_*(o(\calc(i_\lambda);R)) \\
&=\sum_{n \geq 0} (-1)^n \cdot \sum_{\lambda \in \Lambda_n}
\alpha^\cald(i_\lambda)_*\circ (\psi_{i_\lambda})_\ast(o(\calc(i_\lambda);R)) \\
&=\sum_{n \geq 0} (-1)^n \cdot \sum_{\lambda \in \Lambda_n} \alpha^\cald(i_\lambda)_*(
o(\cald(i_\lambda);R)).
\endaligned$$
\\[1mm]~%
\ref{the:homotopy_colimit_formula:chi} follows
from~\ref{the:homotopy_colimit_formula:directly_finite},~%
\ref{the:homotopy_colimit_formula:(FP)},
and~\ref{the:homotopy_colimit_formula:o} in the same way that
Theorem~\ref{the:homotopy_colimit_formula}~\ref{the:homotopy_colimit_formula:chi}
follows from
Theorem~\ref{the:homotopy_colimit_formula}~\ref{the:homotopy_colimit_formula:(EI)},~%
\ref{the:homotopy_colimit_formula:(FP)}, and~\ref{the:homotopy_colimit_formula:o}.
\\[1mm]~%
\ref{the:homotopy_colimit_formula:chi(2)} Suppose that $\cali$ is
directly finite and $\cald(i)$ is directly finite for every object $i
\in \ob(\cali)$. Suppose also for every object $i \in \cali$ the category
$\cald(i)$ is of type ($L^2$). By the proof of
Corollary~\ref{cor:homotopy_colimit_formula_for_pseudo_functors}~\ref{the:homotopy_colimit_formula:directly_finite}
above, the values of the strict functor $\calc$ are directly finite
categories. If $\Gamma_1$ and $\Gamma_2$ are equivalent
categories, then $\Gamma_1$ is both directly finite and of type
($L^2$) if and only if $\Gamma_2$ is both directly finite and of
type ($L^2$)
\cite[Lemma~5.15~(i)]{FioreLueckSauerFinObsAndEulCharOfCats(2009)}.
Since each $\cald(i)$ is directly finite, of type ($L^2$), and
equivalent to $\calc(i)$, we see that each $\calc(i)$ is also
directly finite and of type ($L^2$). So we may now apply
Theorem~\ref{the:homotopy_colimit_formula}~\ref{the:homotopy_colimit_formula:directly_finite}
and \ref{the:homotopy_colimit_formula:chi(2)} to $\calc$ and
conclude that $\hocolim_\cali \calc$ is directly finite and of type
($L^2$).  Again using the preservation of the direct finiteness and
($L^2$) under equivalence
\cite[Lemma~5.15~(i)]{FioreLueckSauerFinObsAndEulCharOfCats(2009)},
and the equivalence of $\hocolim_\cali \calc$ with $\hocolim_\cali
\cald$, we see $\hocolim_\cali \cald$ is both directly finite and of
type ($L^2$).

To prove the formulas for $\chi_f^{(2)}$ and $\chi^{(2)}$, we use
\cite[Lemma~5.15~(ii)]{FioreLueckSauerFinObsAndEulCharOfCats(2009)},
which says: if $F \colon \Gamma_1 \to \Gamma_2$ is an equivalence of
categories, and both $\Gamma_1$ and $\Gamma_2$ are both directly
finite and of type ($L^2$), then
$U^{(1)}(F)\chi^{(2)}_f(\Gamma_1)=\chi^{(2)}_f(\Gamma_2)$ and
$\chi^{(2)}(\Gamma_1) = \chi^{(2)}(\Gamma_2).$ We apply this to the
equivalences $\psi_i$ and $\hocolim_\cali \psi$, and use the
commutativity of diagram
\eqref{equ:naturality_for_alpha_to_hocolim}. For readability, we
write $(\hocolim_\cali \psi)_*$ for $U(\hocolim_\cali \psi)$ and
$\alpha(i_\lambda)_*$ for $U^{(1)}(\alpha(i_\lambda))$, et cetera.
$$\aligned
\chi^{(2)}_f(\hocolim_\cali \cald)&=(\hocolim_\cali \psi)_* \chi^{(2)}_f(\hocolim_\cali \calc) \\
&=(\hocolim_\cali \psi)_*\sum_{n \geq 0} (-1)^n \cdot \sum_{\lambda \in \Lambda_n}
  \alpha(i_\lambda)_*\bigl(\chi_f^{(2)}(\calc(i_\lambda))\bigr) \\
&=\sum_{n \geq 0} (-1)^n \cdot \sum_{\lambda \in \Lambda_n}
(\hocolim_\cali\psi)_\ast\circ \alpha^\calc(i_\lambda)_*\bigl(\chi_f^{(2)}(\calc(i_\lambda))\bigr) \\
&=\sum_{n \geq 0} (-1)^n \cdot \sum_{\lambda \in \Lambda_n}
\alpha^\cald(i_\lambda)_*\circ (\psi_{i_\lambda})_\ast\bigl(\chi_f^{(2)}(\calc(i_\lambda))\bigr)\\
&=\sum_{n \geq 0} (-1)^n \cdot \sum_{\lambda\in \Lambda_n} \alpha^\cald(i_\lambda)_* \bigl(\chi_f^{(2)}(\cald(i_\lambda))\bigr).\\
\endaligned$$
The formula for $\chi^{(2)}$ follows by summing up the components of the functorial $L^2$-Euler characteristics.
\end{proof}

\subsection{The Case of an Indexing Category of Type (FP$_R$)}
\label{subsec:The_case_of_an_indexing_categegory_of_type_(FP)}

The Homotopy Colimit Formula of Theorem~\ref{the:homotopy_colimit_formula} can be extended to the
case, where $\cali$ is of type (FP$_R$) and not necessarily of type (FF$_R$)
as follows (recall that the existence of a finite $\cali$-$CW$-model for $E\cali$ implies $\cali$ is of type (FF$_R$), since cellular chains then provide a finite free resolution of $\underline{R}$. ). The evaluation of the covariant functor
$$E^{\calh} \colon \cali \to \calh\text{-}\SPACES$$
of \eqref{Ecalh} at every object $i \in \cali$ is an $\calh$-$CW$-complex. Composing
it with the cellular chain complex functor yields a covariant functor
$$C_*(E^{\calh}) \colon \cali \to R\calh\text{-}\CHAINCOMPLEXES$$
whose evaluation at every object in $\cali$ is a free $R\calh$-chain
complexes. Since by assumption $\calc(i)$ is of type (FP$_R$),
$C_*(E^{\calh})(i)$ is $R\calh$-chain homotopy equivalent to a
finite projective  $R\calh$-chain complex for every object $i \in
\cali$. Since $R\mor_{\cali}(?,i) \otimes_{R\cali} C_*(E^{\calh})$
is $R\calh$-isomorphic to $C_*(E^{\calh})$, we conclude for every
finitely generated projective $R\Gamma$-module $P$ that $P
\otimes_{R\cali} C_*(E^{\calh})$ is $R\calh$-chain homotopy
equivalent to finite projective  $R\calh$-chain complex and in
particular possesses a finiteness obstruction $o\bigl(P
\otimes_{R\cali} C_*(E^{\calh}\bigr) \in K_0(R\calh)$
(see~L\"uck~\cite[Theorem~11.2 on page~212]{Lueck(1989)}). Because
of~L\"uck~\cite[Theorem~11.2 on page~212]{Lueck(1989)} we obtain a
homomorphism
$$\alpha_{\calc} \colon K_0(R\cali) \to K_0(R\calh),
\quad [P] \mapsto o\bigl(P \otimes_{R\cali} C_*(E^{\calh})\bigr).$$
The chain complex version of the proof of
Lemma~\ref{lem:model_for_E(hocolim_cali_calc)} shows that the
$R\calh$-chain complex $C_*(\cali) \otimes_{R\cali} C_*(E^{\calh})$
is a projective $R\calh$-resolution of the constant $R\Gamma$-module
$\underline{R}$. Choose a finite projective $R\cali$-chain complex
$P_*$ and an $R\cali$-chain homotopy equivalence $f_* \colon P_*
\xrightarrow{\simeq}  C_*(\cali)$. Then $f_*\otimes_{R\cali} \id
\colon P_* \otimes_{R\cali} C_*(E^{\calh}) \to C_*(\cali)
\otimes_{R\cali} C_*(E^{\calh})$ is an $R\Gamma$-chain homotopy
equivalence of $R\Gamma$-chain complexes and $P_* \otimes_{R\cali}
C_*(E^{\calh})$ is is $R\calh$-chain homotopy equivalent to finite
projective $R\calh$-chain complex by L\"uck~\cite[Theorem~11.2 on
page~212]{Lueck(1989)}. This implies
$$o(\Gamma;R) = o\bigl(P_* \otimes_{R\cali} C_*(E^{\calh})\bigr).$$
We conclude from~\cite[Theorem~11.2 on page~212]{Lueck(1989)}
$$
o\bigl(P_* \otimes_{R\cali} C_*(E^{\calh})\bigr) =
\sum_{n \geq p} (-1)^n \cdot o\bigl(P_n \otimes_{R\cali} C_*(E^{\calh})\bigr)
$$
Since $o(\cali;R)$ is $\sum_{n \geq p} (-1)^n \cdot [P_n]$, this implies

\begin{theorem}[The Homotopy Colimit Formula for an indexing category of type (FP$_R$)]
\label{the:the_homotopy_colimit_formula_for_an_indexing_category_of_type_(FP)}
We obtain under the conditions above
$$\alpha_{\calc} \bigl(o(\cali;R)\bigr)
=
o(\calh;R).
$$
\end{theorem}


\begin{remark}
See Section \ref{sec:comparison_with_Leinster} for a comparison with Leinster's Euler characteristic and his results.
\end{remark}


\typeout{---------------  Section: Examples  of the Homotopy Colimit Formula  --------------------}

\section{Examples of the Homotopy Colimit Formula}
\label{sec:examples_of_the_homotopy_colimit_formula}

We now present several examples of the Homotopy Colimit Formula
Theorem~\ref{the:homotopy_colimit_formula}. These include the cases:
$\cali$ with a terminal object, the constant functor, the trivial
functor, homotopy pushouts, homotopy orbits, and the transport
groupoid. For the transport groupoid in the finite case, see also
Example~\ref{exa:transport_groupoid_finite_case}.

\begin{example}[Homotopy Colimit Formula for $\cali$ with a terminal
  object] \label{exa:hocolim_formula_for_I_with_terminal_object} Suppose that
  $\cali$ has a terminal object $t$ and $\calc \colon \cali \to \CATS$ is a
  functor. Then $\mor_\cali(-,t)$ is a finite $\cali$-CW model for $E\cali$. If
  every category $\calc(i)$ is of type (FP$_R$), then
  $o(\calh;R)=\alpha(t)_*(o(\calc(t);R)$.  If $\cali$ and $\calc$ additionally
  satisfy the hypotheses of
  Theorem~\ref{the:homotopy_colimit_formula}~\ref{the:homotopy_colimit_formula:chi},
  then $\chi_f(\calh;R)=\chi_f(\calc(t);R)$ and
  $\chi(\calh;R)=\chi(\calc(t);R)$, as we anticipated in
  Example~\ref{exa:homotopy_hocolimit_and_when_I_has_terminal_object}.
  Similar statements hold for $\chi^{(2)}_f$ and $\chi^{(2)}$ in the
  $L^2$ case.
\end{example}

\begin{example}[Homotopy Colimit Formula for a constant functor]
\label{exa:homotopy_hocolimit_formula_of_constant_functor} Consider
the situation of Theorem~\ref{the:homotopy_colimit_formula} in the
special case where the covariant functor $\calc \colon \cali \to
\CATS$ is constant $\calc \in \CATS$. Suppose that $\cali$ admits a
finite $\cali$-$CW$-model for $E\cali$. Then we may draw various
conclusions about the homotopy colimit $\calh=\cali \times \calc$.
If $\cali$ and $\calc$ are of type (FP$_R$), then so is $\cali
\times \calc$. If $\cali$ and $\calc$ are of type (FF$_R$), then so
is $\cali \times \calc$. The statements in
Theorem~\ref{the:homotopy_colimit_formula} provide us with formulas
in terms of $\calc$ for $o(\cali\times \calc;R)$,
$\chi_f(\cali\times \calc;R)$, $\chi(\cali\times \calc;R)$,
$\chi_f^{(2)}(\cali\times \calc)$, and $\chi^{(2)}(\cali\times
\calc)$. We recall that the invariants $o$, $\chi_f$, $\chi$,
$\chi_f^{(2)}$, and $\chi^{(2)}$ are multiplicative, see
Fiore--L\"uck--Sauer~\cite[Theorems~2.17,~4.22,~and~5.17]{FioreLueckSauerFinObsAndEulCharOfCats(2009)}.
\end{example}

\begin{example}[Homotopy Colimit Formula for the trivial functor]
  \label{exa:homotopy_hocolimit_and_trivial_functor}
  Consider the situation of Theorem~\ref{the:homotopy_colimit_formula}
  in the special case where the covariant functor $\calc \colon \cali
  \to \CATS$ is constantly the terminal category, which consists of a single object and its identity morphism. Then $\hocolim_{\cali} \calc$
  agrees with $\cali$, as we see from Example~\ref{exa:hocolim_constant}.
Obviously $\calc(i)$ is of type (FF$_R$), its
  finiteness obstruction is $[R] \in K_0(R) = K_0(R\calc(i))$ and both
  its Euler characteristic and $L^2$-Euler characteristic equals
  $1$. We obtain from Theorem~\ref{the:homotopy_colimit_formula}
$$
\begin{array}{lcll}
  o(\cali;R)
  & = &
  \sum_{n \geq 0} (-1)^n \cdot \sum_{\lambda \in \Lambda_n} [R\mor_{\cali}(?,i_\lambda)] & \in K_0(R\cali);
  \\
  \chi_f(\cali;R)
  & = &
  \sum_{n \geq 0} (-1)^n \cdot \sum_{\lambda \in \Lambda_n} \overline{i_\lambda} & \in U(\Gamma);
  \\
  \chi(\cali;R)
  & = &
  \sum_{n \geq 0} (-1)^n \cdot |\Lambda_n| & \in \IZ;
  \\
  \chi_f^{(2)}(\cali)
  & = &
  \sum_{n \geq 0} (-1)^n \cdot \sum_{\lambda \in \Lambda_n} \overline{i_\lambda} & \in U^{(1)}(\cali);
  \\
  \chi^{(2)}(\cali)
  & = &
  \sum_{n \geq 0} (-1)^n \cdot |\Lambda_n| & \in \IR.
\end{array}
$$
\end{example}

\begin{example}[Homotopy pushout formula]
  \label{exa:homotopy_pushout}
  Let $\cali$ be the category with objects $j$, $k$ and $\ell$ such that
  there is precisely one morphism from $j$ to $k$ and from $j$ to $\ell$
  and all other morphisms are identity morphisms.
  $$\cali=\{\xymatrix{k & j \ar[l]_{g} \ar[r]^{h} & \ell}\}$$
  By Example~\ref{exa:finite_model_for_pushout}, the category $\cali$ admits a
  finite model for the classifying $\cali$-space $E\cali$.

  A covariant functor $\calc \colon \cali \to \CATS$ is the same as specifying
  three categories $\calc(j)$, $\calc(k)$ and $\calc(\ell)$ and two
  functors $\calc(g) \colon \calc(j) \to \calc(k)$ and $\calc(h)
  \colon \calc(j) \to \calc(\ell)$. Let $\calh = \hocolim_{\cali} \calc$
  be the homotopy colimit. Let $\alpha(i) \colon \calc(i) \to \calh$
  be the canonical functor for $i = j,k,\ell$. Then we obtain a square
  of functors which commutes up to natural transformations
  $$
  \xymatrix{\calc(j) \ar[r]^{\calc(g)} \ar[d]_{\calc(h)} \ar[dr]^{\alpha(j)}
  & \calc(k) \ar[d]^{\alpha(k)}
  \\
  \calc(\ell) \ar[r]_{\alpha(\ell)} & \calh.}
  $$
  It induces  diagrams which do {\bf not} commute in general
  $$
  \xymatrix{K_0(R\calc(j)) \ar[r]^{\calc(g)_*} \ar[d]_{\calc(h)_*} \ar[dr]^{\alpha(j)_*}
  & K_0(R\calc(k)) \ar[d]^{\alpha(k)_*}
  \\
  K_0(R\calc(\ell)) \ar[r]_{\alpha(\ell)_*} & K_0(\calh)}
  $$
  and
  $$
  \xymatrix{U(\calc(j)) \ar[r]^{\calc(g)_*} \ar[d]_{\calc(h)_*} \ar[dr]^{\alpha(j)_*}
  & U(R\calc(k)) \ar[d]^{\alpha(k)_*}
  \\
  U(R\calc(\ell)) \ar[r]_{\alpha(\ell)_*} & U(\calh).}
  $$
  Suppose that $\calc(i)$ is of type (FP$_R$) for $i = j,k,\ell$. We conclude from
  Theorem~\ref{the:homotopy_colimit_formula}~\ref{the:homotopy_colimit_formula:(FP)}
  that $\calh$ is of type (FP$_R$) and
  $$\begin{array}{lcll}
  o(\calh;R)
  & = &
  \alpha(k)_*\bigl(o(\calc(k);R)) + \alpha(\ell)_*\bigl(o(\calc(\ell);R))
  - \alpha(j)_*\bigl(o(\calc(j;R));
  & \in K_0(R\calh);
  \\
  \chi_f(\calh;R)
  & = &
  \alpha(k)_*\bigl(\chi_f(\calc(k);R)\bigr) + \alpha(\ell)_*\bigl(\chi_f(\calc(\ell);R)\bigr) -
  \alpha(j)_*\bigl(\chi_f(\calc(j);R)\bigr);
  & \in U(\calh);
  \\
  \chi(\calh;R)
  & = &
  \chi(\calc(k);R) + \chi(\calc(\ell);R) - \chi(\calc(j);R);
  & \in \IZ;
  \\
  \chi_f^{(2)}(\calh)
  & = &
  \alpha(k)_*\bigl(\chi^{(2)}_f(\calc(k)\bigr) + \alpha(\ell)_*\bigl(\chi^{(2)}_f(\calc(\ell))\bigr) -
  \alpha(j)_*\bigl(\chi^{(2)}_f(\calc(j))\bigr);
  & \in U^{(1)}(\calh);
  \\
  \chi^{(2)}(\calh)
  & = &
  \chi^{(2)}(\calc(k)) + \chi^{(2)}(\calc(\ell)) - \chi^{(2)}(\calc(j));
  & \in \IR.
  \end{array}
  $$
\end{example}

\begin{example}[Homotopy orbit formula]
  \label{exa:homotopy_colimit_and_group_actions}
  Suppose that a group $G$ acts on a category $\calc$ from the
  left.  This can be viewed as a covariant functor $\widehat{G} \to
  \CATS$ whose source is the groupoid $\widehat{G}$
  with one object and $G$ as its automorphism
  group. Let $\calh = \hocolim_{\widehat{G}} \calc$ be its homotopy
  colimit, also called the \emph{homotopy orbit}.
  Notice that $\calh$ and $\calc$ have the same set of
  objects.

  Suppose there is a finite model for $BG$ of the group $G$, or
  equivalently, a finite model for the $\widehat{G}$-classifying space $E\widehat{G}$
  of the category $\widehat{G}$.  Let $\chi(BG) \in \IZ$ be its Euler characteristic.
  Let $\alpha \colon \calc \to \calh$
  be the canonical inclusion. Suppose that $\calc$ is of type
  (FP$_R$). Then we conclude from
  Theorem~\ref{the:homotopy_colimit_formula}~\ref{the:homotopy_colimit_formula:(FP)} that $\calh$ is of type
  (FP$_R$) and we have
  $$\begin{array}{lcll}
  o(\calh;R) & = & \chi(BG) \cdot \alpha_*\bigl(o(\calc;R)\bigr) & \in K_0(R\calh);
  \\
  \chi_f(\calh;R) & = & \chi(BG) \cdot \alpha_*\bigl(\chi_f(\calc;R)\bigr)
  & \in U(\calh);
  \\
  \chi(\calh;R) & = & \chi(BG) \cdot \chi(\calc;R) & \in \IZ;
  \\
  \chi_f^{(2)}(\calh;R) & = & \chi(BG) \cdot \alpha_*\bigl(\chi_f^{(2)}(\calc;R)\bigr)
  & \in U^{(1)}(\calh);
  \\
  \chi^{(2)}(\calh;R) & = & \chi(BG) \cdot \chi^{(2)}(\calc;R) & \in \IR.
  \end{array}
  $$
\end{example}

\begin{example}[Transport groupoid]
  \label{exa:transport_groupoid}
  Let $G$ be a group and let $S$ be a left $G$-set. Its
  \emph{transport groupoid} $\calg^G(S)$ has $S$ as its set of objects. The set
  of morphisms from $s_1$ to $s_2$ is $\{g \in G \mid gs_1 = s_2\}$. The
  composition is given by the multiplication in $G$.  Denote by
  $\underline{S}$ the category whose set of objects is $S$ and which
  has no morphisms besides the identity morphisms. The group $G$ acts
  from the left on $\underline{S}$.  One easily checks that $\calg^G(S)$ is the
  homotopy orbit of $\underline{S}$ defined in
  Example~\ref{exa:homotopy_colimit_and_group_actions}.

Recall from
Fiore--L\"uck--Sauer~\cite[Lemma~6.15~(iv)]{FioreLueckSauerFinObsAndEulCharOfCats(2009)}:
 if $\Gamma$ is a quasi-finite EI-category and for any morphism $f
\colon x \to y$ in $\Gamma$, the order of the finite group $\{g \in
\aut(x) \mid f \circ g = f\}$ is invertible in $R$, then $\Gamma$ is
of type (FP$_R$) if and only if $\iso(\Gamma)$ is finite and for
every object $x \in \ob(\Gamma)$ the trivial $R[x]$-module $R$ is of
type (FP$_R$). Thus, category $\underline{S}$ is of type (FP$_R$) if
and only if $S$ is finite. Suppose that $\underline{S}$ is of type
(FP$_R$) and there is a finite model for $BG$.  Obviously
$o(\underline{S};R)$ is given in $K_0(R\underline{S}) = \oplus_S
K_0(R)$ by the collection $\{[R] \in K_0(R) \mid s \in S\}$.

  Suppose for simplicity that $G$ acts transitively on $S$. Fix an
  element $s \in S$. Let $G_s$ be its isotropy group. Since $S$ is
  finite, $G_s$ is a subgroup of $G$ of finite index, namely
  $[G:G_s] = |S|$.  The transport groupoid $\calg^G(S)$ is connected and the
  automorphism group of $s$ is $G_s$. Hence evaluation at $s$ induces
  an isomorphism
  $$\ev \colon K_0(R\calg^G(S)) \xrightarrow{\cong} K_0(R[G_s]).$$
  The composition
  $$K_0(R\underline{S}) \xrightarrow{\alpha_*} K_0(R\calg^G(S))
  \xrightarrow{\cong} K_0(R[G_s])$$
  sends $o(\underline{S};R)$ to $|S| \cdot [RG_s]$, where
  $\alpha \colon \underline{S} \to \calg^G(S)$ is the obvious inclusion. Hence
  Example~\ref{exa:homotopy_colimit_and_group_actions} implies
  $$\ev\bigl(o(\calg^G(S);R)\bigr)
  = \chi(BG) \cdot |S| \cdot [RG_S] \quad \in K_0(RG_s).$$
  Since $BG$ has a finite model, $BG_s$ as a finite covering of $BG$ has a finite
  model. The cellular $RG_s$-chain complex of $EG_s$ yields a finite
  free resolution of the trivial $RG_s$-module $R$. This implies
  $$\ev\bigl(o(\calg^G(S);R)\bigr)
  = \chi(BG_s) \cdot [RG_s] \quad \in K_0(RG_s).$$
  Hence we obtain the equality in $K_0(RG_s)$
  $$\chi(BG_s) \cdot [RG_s] = \chi(BG) \cdot |S| \cdot [RG_S]
  = \chi(BG) \cdot [G:G_s] \cdot [RG_s].$$
  This is equivalent to the equality of integers
  $$\chi(BG_s)  = \chi(BG) \cdot [G:G_s].$$
  This equation is compatible with the well-know fact that for a $d$-sheeted
  covering $\overline{X} \to X$ of a finite $CW$-complex $X$ the total
  space $\overline{X}$ is again a finite $CW$-complex and we have
  $\chi(\overline{X}) = d \cdot \chi(X)$.

  For the transport groupoid in the finite case, see also
Example~\ref{exa:transport_groupoid_finite_case}.
\end{example}


\typeout{---- Section 11: Combinatorial Illustrations of the Homotopy Colimit Formula ---}
\section{Combinatorial Illustrations of the Homotopy Colimit Formula}
\label{sec:Combinatorial_Applications_of_the_Homotopy_Colimit_Formula}

The classical Inclusion-Exclusion Principle follows from the
Homotopy Colimit Formula Theorem~\ref{the:homotopy_colimit_formula}.
We can also easily calculate the cardinality of a coequalizer in
$\SETS$ in certain cases. These are different proofs of Examples
3.4.d and 3.4.b of Leinster's paper \cite{Leinster(2008)}.

\begin{example}[Inclusion-Exclusion
Principle]\label{the:inclusion_exclusion} Let $X$ be a set and $S_0,
\dots, S_q$ finite subsets of $X$. Then
$$|S_0 \cup S_1 \cup \cdots \cup S_q|=
\sum_{\emptyset \neq J \subseteq [q]}(-1)^{|J|-1} \cdot
\left|\bigcap_{j \in J} S_j\right|.$$
\end{example}
\begin{proof}
Let $\cali$ be the category in Example~\ref{exa:finite_model_for_q_interior}
and consider the finite
$\cali$-$CW$-model for its classifying $\cali$-space constructed
there. We define a functor $\calc: \cali \to \SETS$ by
$\calc(J):=\bigcap_{j\in J} S_j$. The functor
$$\hocolim_{\cali} \calc \xymatrix{\ar[r] &} \colim_{\cali} \calc = S_0 \cup S_1 \cup \cdots \cup S_q$$
is an equivalence of categories, since it is surjective on objects
and fully faithful. We have
$$\aligned |S_0 \cup S_1 \cup \cdots \cup S_q| &=\chi(S_0 \cup S_1 \cup \cdots \cup S_q)
\\ &= \chi(\hocolim_{\cali} \calc)
\\ &= \sum_{n \geq 0} (-1)^n \cdot \sum_{\lambda \in \Lambda_n} \chi (\calc(i_\lambda))
\\ &= \sum_{n \geq 0} (-1)^n \cdot \sum_{J \subseteq [q] \;\text{and}\;
|J|=n+1} \chi(\calc(J))
\\ &= \sum_{n \geq 0} (-1)^n \left( \sum_{J \subseteq [q] \;\text{and}\;
|J|=n+1} \left|\bigcap_{j \in J} S_j \right| \right)
\\&= \sum_{\emptyset \neq J \subseteq [q]} \left( (-1)^{|J|-1} \left|\bigcap_{j \in J} S_j \right|
\right).
\endaligned$$
\end{proof}

\begin{example}[Cardinality of a Coequalizer]
\label{the:cardinality_of_coequalizer}
Let $\cali$ be the category
$$\xymatrix{a \ar@<.5ex>[r]^f \ar@<-.5ex>[r]_g & b}$$
and $\calc:\cali \to \SETS$ a functor such that:
\begin{enumerate}
\item
the maps $\calc f$ and $\calc g$ are injective,
\item
the images of the maps $\calc f$ and $\calc g$ are disjoint, and
\item
the sets $\calc a$ and $\calc b$ are finite.
\end{enumerate}
Then the coequalizer $\colim \calc$ has cardinality $|\calc b|-|\calc a|$.
\end{example}
\begin{proof}
The assumptions that $\calc f$ and $\calc g$ are injective and have
disjoint images imply that the functor
$$\hocolim_{\cali} \calc \xymatrix{\ar[r] &} \colim_{\cali} \calc$$
is fully faithful. Clearly it is also surjective on objects, and
hence an equivalence of categories. The category $\cali$ has a
finite $\cali$-$CW$-model for its classifying $\cali$-space,
constructed explicitly in Example~\ref{exa:finite_model_for_parallel_arrows}. By
Theorem~\ref{the:homotopy_colimit_formula}, we have
$$\aligned
 \chi(\colim_{\cali} \calc)&=\chi(\hocolim_{\cali}\calc)
\\ &= \sum_{n \geq 0} (-1)^n \cdot \sum_{\lambda \in \Lambda_n} \chi (\calc(i_\lambda))
\\ &= \chi(\calc b) - \chi(\calc a)
\\ &= |\calc b| - |\calc a|.
\endaligned$$
\end{proof}

\typeout{------------- Section: Comparison with Results of
Baez--Dolan and Leinster  ---------}
\section{Comparison with Results of Baez--Dolan and Leinster} \label{sec:comparison_with_Leinster}

We recall Baez--Dolan's groupoid cardinality \cite{Baez-Dolan(2001)}
and Leinster's Euler characteristic for certain finite categories
\cite{Leinster(2008)}, compare our Homotopy Colimit Formula with his
result on compatibility with Grothendieck fibrations, prove an
analogue for indexing categories $\cali$ that admit finite
$\cali$-$CW$-models for their classifying $\cali$-spaces, and
finally mention a Homotopy Colimit Formula for Leinster's invariant
in a restricted case.

\subsection{Review of Leinster's Euler Characteristic} Let $\Gamma$ be a
category with finitely many objects and finitely many morphisms. A
\emph{weighting} on $\Gamma$ is a function $q^{\bullet} \colon
\ob(\Gamma) \to \IQ$ such that for all objects $x \in \ob(\Gamma)$,
we have \[\sum_{y \in \ob(\Gamma)} |\mor_\Gamma(x,y)| \cdot q^y =
1.\] A \emph{coweighting} $q_{\bullet}$ on $\Gamma$ is a weighting
on $\Gamma^{\op}$. If a finite category admits both a weighting
$q^{\bullet}$ and a coweighting $q_{\bullet}$, then $\sum_{y \in
\ob(\Gamma)} q^y = \sum_{x \in \ob(\Gamma)} q_x$. For a discusion of
which matrices have the form $\left(|\mor_\Gamma(x,y)| \right)_{x,y
\in \ob(\Gamma)}$ for some finite category $\Gamma$, see Allouch
\cite{Allouch2008} and \cite{Allouch2010}.

As proved in \cite{FioreLueckSauerFinObsAndEulCharOfCats(2009)},
free resolutions of the constant $R\Gamma$-module $\underline{R}$
give rise to weightings on $\Gamma$.

\begin{theorem}[Weighting from a free resolution, Theorem~7.6~of~Fiore--L\"uck--Sauer~\cite{FioreLueckSauerFinObsAndEulCharOfCats(2009)}]
\label{the:weighting_from_free_resolution} Let $\Gamma$ be a finite
category.
Suppose that the constant $R \Gamma$-module
$\underline{R}$ admits a finite free resolution $P_*$. If $P_n$ is
free on the finite $\ob(\Gamma)$-set $C_n$, that is
\begin{equation}\label{the:weighting_from_free_resolution:equation}
P_n=B(C_n)=\bigoplus_{y \in \ob(\Gamma)} \bigoplus_{C^y_n}
R\mor_\Gamma(?,y),
\end{equation}
then the function $q^{\bullet} \colon \ob(\Gamma) \to \IQ$ defined
by
$$q^y:=\sum_{n \geq 0}(-1)^n \cdot |C^y_n|$$
is a weighting on $\Gamma$.
\end{theorem}

\begin{corollary}[Construction of a weighting from a
finite $\cali$-$CW$-model for the classifying $\cali$-space,
Corollary~7.8~of~Fiore--L\"uck--Sauer~\cite{FioreLueckSauerFinObsAndEulCharOfCats(2009)}]
\label{cor:weighting_from_finite_model} Let $\cali$ be a finite
category. Suppose that $\cali$ admits a finite $\cali$-$CW$-model
$X$ for the classifying $\cali$-space. Then the function
$q^{\bullet} \colon \ob(\cali) \to \IQ$ defined by
$$q^y:=\sum_{n \geq 0}(-1)^n(\text{number of $n$-cells of $X$ based at $y$})$$
is a weighting on $\cali$.
\end{corollary}

As explained in
Section~7.5~of~\cite{FioreLueckSauerFinObsAndEulCharOfCats(2009)},
we use this Corollary to obtain several of Leinster's weightings in
\cite{Leinster(2008)} from $\cali$-$CW$-models for
$\cali$-classifying spaces. If $\cali$ has a terminal object, then
we obtain from the finite model in Example
\ref{exa:finite_model_for_I_with_terminal_object} the weighting
which is 1 on the terminal object and 0 otherwise. The category
$\cali=\{j \rightrightarrows k\}$ in Example
\ref{exa:finite_model_for_parallel_arrows} has weighting
$(q^j,q^k)=(-1,1)$. The category $\cali=\{k \leftarrow j \to \ell\}$
in Example \ref{exa:finite_model_for_pushout} has weighting
$(q^j,q^k,q^\ell)=(-1,1,1)$. Lastly, the category in Example
\ref{exa:finite_model_for_q_interior} has weighting
$q^J:=(-1)^{|J|-1}$.

Weightings and coweightings play a key role in Leinster's notion of
Euler characteristic. See also
Berger--Leinster~\cite{Berger+Leinster(2008)}.

\begin{definition}[Definition 2.2 of Leinster~\cite{Leinster(2008)}] \label{def:Leinsters_Euler_characteristic}
A finite category $\Gamma$ \emph{has an Euler characteristic in the
sense of Leinster} if it admits both a weighting and a coweighting. In this case, its
\emph{Euler characteristic in the sense of Leinster} is defined as
$$\chi_L(\Gamma) := \sum_{y \in \ob(\Gamma)} q^y = \sum_{x \in \ob(\Gamma)} q_x$$
for any choice of weighting $q^{\bullet}$ or coweighting
$q_{\bullet}$.
\end{definition}

The Euler characteristic of Leinster agrees with the \emph{groupoid
cardinality} of Baez--Dolan~\cite{Baez-Dolan(2001)} in the case of a
finite groupoid $\calg$, namely they are both $$\sum_{\overline{x}
\in \iso(\calg)} \frac{1}{|\aut_{\calg}(x)|}.$$ The Euler
characteristic of Leinster agrees with our $L^2$-Euler
characteristic in some cases, as in the following Lemma.

\begin{lemma}[Lemma~7.3~of~Fiore--L\"uck--Sauer~\cite{FioreLueckSauerFinObsAndEulCharOfCats(2009)}] \label{lem:chi(2)_and_chi}
Let $\Gamma$ be a finite EI-category which is skeletal, i.e.,
if two objects are isomorphic, then they are equal. Suppose that the left
$\aut_\Gamma(y)$-action on $\mor_\Gamma(x,y)$ is free for every two objects $x,y
\in \ob(\Gamma)$.

Then $\Gamma$ is of type (FP$_\IC$) and of type ($L^2$), and has an
Euler characteristic in the sense of Leinster. Furthermore, the
$L^2$-Euler characteristic $\chi^{(2)}(\Gamma)$ of
Definition~\ref{def:L2-Euler_characteristic_of_a_category} coincides
with Leinster's Euler characteristic $\chi_L(\Gamma)$ of
Definition~\ref{def:Leinsters_Euler_characteristic}:
$$\chi^{(2)}(\Gamma) = \chi_L(\Gamma).$$
Moreover, these are both equal to
\[
\sum_{l \ge 0} (-1)^l \cdot \sum_{x_0,x_l \in \ob(\Gamma)} \sum
\frac{1}{|\aut(x_{l})| \cdot |\aut(x_{l-1})|  \cdot \cdots \cdot
|\aut(x_{0})|},
\]
where the inner sum is over all paths $x_0 \to x_1 \to \cdots \to
x_l$ from $x_0$ to $x_l$ such that $x_0, \ldots, x_l$ are all
distinct
\cite[Example~6.33]{FioreLueckSauerFinObsAndEulCharOfCats(2009)}.
\end{lemma}

This concludes the review of Leinster's and Baez--Dolan's invariants
and how they relate to our $L^2$-Euler characteristic. Next we turn
to a comparison of homotopy colimit results.

\subsection{Comparison with Leinster's Proposition 2.8} Leinster's result on homotopy colimits, rephrased in our notation to make the comparison more apparent, is below.
\begin{theorem}[Proposition 2.8 of Leinster~\cite{Leinster(2008)}]
Let $\cali$ be a category with finitely many objects and finitely many morphisms, and $\calc:\cali \to \CATS$ a pseudo functor. Assume
that $\hocolim_\cali \calc$ has finitely many objects and finitely many morphisms. Let $q^{\bullet}$ be a weighting on $\cali$ and suppose
that $\hocolim_\cali \calc$ and all $\calc(i)$ have Euler characteristics. Then
\[\chi_L(\hocolim_\cali \calc) =\sum_{i \in \ob(\cali)} q^i \chi_L(\calc(i)).
\]
\end{theorem}

For example, if $\cali=\{k \leftarrow j \to \ell\}$, then $\cali$ admits the weighting $(q^j,q^k,q^\ell)=(-1,1,1)$ as discussed above. If
$\calc:\cali \to \CATS$ is a pseudo functor, and the homotopy pushout has finitely many objects and finitely many morphisms, and $\hocolim_\cali \calc$ and all $\calc(i)$ have Euler characteristics, then Leinster's result says that the homotopy pushout has the Euler characteristic $\chi_L(\calc(k))+\chi_L(\calc(\ell))-\chi_L(\calc(j))$.

Leinster's Proposition 2.8 tells us how the Euler characteristic is
compatible with Grothendieck fibrations. We can remove the
hypothesis of finite from that Proposition, at the expense of
requiring a finite model, as in the following theorem for our invariants.

\begin{theorem} \label{the:hocolim_weighting}
Let $\cali$ be a finite category. Suppose that $\cali$ admits a
finite $\cali$-$CW$-model $X$ for the classifying $\cali$-space of
$\cali$. Let $q^{\bullet} \colon \ob(\cali) \to \IQ$ be the
$\cali$-Euler characteristic of $X$, namely
$$q^i:=\sum_{n \geq 0}(-1)^n(\text{number of $n$-cells of $X$ based at $i$}).$$
Let $\calc:\cali \to \CATS$ be a functor such that for every object
$i$ the category $\calc(i)$ is of type (FP$_R$). Suppose that
$\cali$ is directly finite and $\calc(i)$ is directly finite for all $i \in \ob(\cali)$. Then
$$\chi(\hocolim_{\cali}\calc;R)=\sum_{i\in\ob(\cali)}q^i\chi(\calc(i);R).$$
If each $\calc(i)$ is of type ($L^2$) rather than (FP$_R$), we have
$$\chi^{(2)}(\hocolim_{\cali}\calc)=\sum_{i\in\ob(\cali)}q^i\chi^{(2)}(\calc(i)).$$
\end{theorem}
\begin{proof}
By Theorem~\ref{the:homotopy_colimit_formula}~\ref{the:homotopy_colimit_formula:chi}, we have
$$\aligned
\chi(\hocolim_{\cali}\calc; R)&=\sum_{n \geq 0} (-1)^n \cdot \sum_{\lambda
\in \Lambda_n} \chi (\calc(i_\lambda); R)
\\ &=\sum_{n \geq 0} (-1)^n \cdot \sum_{i\in\ob(\cali)}(\text{number of $n$-cells of $X$ based at $i$})\chi (\calc(i);R)
\\ &=\sum_{i\in\ob(\cali)}\sum_{n \geq 0} (-1)^n(\text{number of $n$-cells of $X$ based at $i$})\chi (\calc(i);R)
\\ &=\sum_{i\in\ob(\cali)}q^i\chi(\calc(i);R).
\endaligned$$
The statement for $\chi^{(2)}$ is proved similarly from Theorem~\ref{the:homotopy_colimit_formula}~\ref{the:homotopy_colimit_formula:chi(2)}.
\end{proof}

\begin{remark}
Whenever $\chi(\colim_{\cali}\calc;R)=\chi(\hocolim_{\cali}\calc;R)$,
Theorem~\ref{the:homotopy_colimit_formula} and
Theorem~\ref{the:hocolim_weighting} can be used to calculate the Euler
characteristic of a colimit. Indeed, the hypotheses of Examples~\ref{the:inclusion_exclusion}
and~\ref{the:cardinality_of_coequalizer} guaranteed the equivalence of
the colimit and the homotopy colimit, and this equivalence was a
crucial ingredient in those proofs. For example, under Leinster's
hypothesis of familial representability on $\mathcal{C}$, each
connected component of $\hocolim_{\cali}\calc$ has an initial
object, so
$$\chi(\hocolim_{\cali}\calc;R)=\chi(\colim_{\cali}\calc;R)$$ (recall
that $\colim_{\cali}\calc$ is the set of connected components of
$\hocolim_{\cali}\calc$ whenever $\calc$ takes values in $\SETS$).
This is the role of familial representability in his Examples 3.4.
\end{remark}

As a corollary to our Homotopy Colimit Formula for the $L^2$-Euler characteristic,
we have a Homotopy Colimit Formula for Leinster's Euler characteristic when they agree.

\begin{corollary}[Homotopy Colimit Formula for Leinster's Euler characteristic]
  Let $\cali$ be a skeletal, finite, EI-category such that the left
  $\aut_\cali(y)$-action on $\mor_\cali(x,y)$ is free for every two objects $x,y
\in \ob(\cali)$. Assume there
  exists a finite $\cali$-$CW$-model for the $\cali$-classifying space of $\cali$. Let $\calc \colon
  \cali \to \CATS$ be a covariant functor such that for each $i \in \ob(\cali)$, the category $\calc(i)$ is a skeletal,
  finite, EI and the left
  $\aut_{\calc(i)}(d)$-action on $\mor_{\calc(i)}(c,d)$ is free for every two objects $c,d
\in \ob(\calc(i))$. Assume for every object $i \in \ob(\cali)$, for
  each automorphism $u \colon i \to i$ in $\cali$, and each
  $\overline{x} \in \iso(\calc(i))$ we have $\overline{\calc(u)(x)} =
  \overline{x}$.

  Then $\calh:=\hocolim_{i \in I} \calc$ is again a skeletal, finite,
  EI-category such that the left $\aut_{\calh}(h)$-action on $\mor_{\calh}(g,h)$ is free for every two objects $g,h \in \ob(\hocolim_{i \in I} \calc)$, and
$$\chi_L(\calh) =
  \sum_{n \geq 0} (-1)^n \cdot \sum_{\lambda \in \Lambda_n} \chi_L(\calc(i_\lambda); R).$$
\end{corollary}
\begin{proof}
The category $\calh$ is an EI-category by Theorem~\ref{the:homotopy_colimit_formula}~\ref{the:homotopy_colimit_formula:(EI)}. Skeletality and
finiteness of $\calh$ follow directly from the skeletality and finiteness of $\cali$ and $\calc(i)$, and the definition of $\calh$. The hypotheses on $\calc(i)$ imply that $\chi^{(2)}(\calc(i))=\chi_L(\calc(i))$ by Theorem~\ref{lem:chi(2)_and_chi}, and similarly $\chi^{(2)}(\calh)=\chi_L(\calh)$. Finally, Theorem~\ref{the:homotopy_colimit_formula}~\ref{the:homotopy_colimit_formula:chi(2)}, which is the Homotopy Colimit Formula for the $L^2$-Euler characteristic $\chi^{(2)}$, implies the formula is also true for Leinster's Euler
characteristic $\chi_L$ in the special situation of the Corollary.
\end{proof}


\typeout{---- Section 15: Complexes of Groups----------------------------------- -------}
\section{Scwols and Complexes of Groups} \label{sec:complexes_of_groups}

As an illustration of the Homotopy Colimit Formula, we consider
Euler characteristics of small categories without loops
(\emph{scwols}) and complexes of groups in the sense of
Haefliger~\cite{Haefliger(1991)}, \cite{Haefliger(1992)} and
Bridson--Haefliger \cite{Bridson-Haefliger(1999)}. One-dimensional
complexes of groups are the classical Bass--Serre graphs of groups
\cite{Serre(1980)}. For finite scwols, the Euler characteristic,
$L^2$-Euler characteristic, and Euler characteristic of the
classifying space all coincide, essentially because finite scwols
admit finite models for their classifying spaces. The Euler
characteristic of a finite scwol is particularly easy to find: one
simply chooses a skeleton, counts the paths of non-identity
morphisms of length $n$, and then computes the alternating sum of
these numbers.

Scwols and complexes of groups are combinatorial models for
polyhedral complexes and group actions on them. The poset of faces
of a polyhedral complex is a scwol. Suppose a group $G$ acts on an
$M_\kappa$-polyhedral complex by isometries preserving cell
structure, and suppose each group element $g \in G$ fixes each cell
pointwise that $g$ fixes setwise. In this case, the quotient is also
an $M_\kappa$-polyhedral complex, say $Q$, and we obtain a pseudo
functor from its scwol of faces into groups.  Namely, to a face
$\overline{\sigma}$ of $Q$, one associates the stabilizer $G_\sigma$
for a selected representative $\sigma$ of $\overline{\sigma}$.
Inclusions of subfaces of $Q$ then correspond to inclusions of
stabilizer subgroups up to conjugation. This pseudo functor is the
complex of groups associated to the group action.

However, it is sometimes easier to work directly with the
combinatorial model rather than with the original
$M_\kappa$-polyhedral complex, and consider instead appropriate
group actions on the associated scwol, as in
Definition~\ref{def:group_action_on_a_scwol}. Then the quotient
category of a scwol is again a scwol, and the associated pseudo
functor on the quotient scwol is called the \emph{complex of groups
associated to the group action}. Any group-valued pseudo functor on
a scwol that arises in this way is called \emph{developable}.

Our main results in this section concern the Euler characteristics
of homotopy colimits of complexes of groups associated to group
actions in the sense of
Definition~\ref{def:group_action_on_a_scwol}.
Theorem~\ref{thm:Euler_characteristic_of_hocolim_of_quotient_complex},
concludes that the Euler characteristic and $L^2$-Euler
characteristic of the homotopy colimit are $\chi(\calx/G)$ and
$\chi^{(2)}(\calx)/\vert G \vert$ respectively, $G$ and $\calx$ are
finite. These formulas provide necessary conditions for
developability. That is, if $F$ is a pseudo functor from a scwol
$\caly$ to groups, one may ask if there are a scwol $\calx$ and a
group $G$ such that $\caly$ is isomorphic to $\calx/G$ and $F$ is
the associated complex of groups. To obtain conditions on
$\chi(\calx)$, $\chi^{(2)}(\calx)$, and $\vert G \vert$, one forms
the homotopy colimit of $F$, calculates its Euler characteristic
and $L^2$-Euler characteristic, and then compares with the formulas
of
Theorem~\ref{thm:Euler_characteristic_of_hocolim_of_quotient_complex}.
A simple case is illustrated in
Example~\ref{exa:necessary_conditions_for_developability:single_arrow}.
Another application of the formulas is the computation of the Euler
characteristic and $L^2$-Euler characteristic for the transport
groupoid of a finite left $G$-set, as in
Example~\ref{exa:transport_groupoid_finite_case}. We finish with
Theorem~\ref{the:extension_of_Haefligers_corollary}, which extends
Haefliger's formula for the Euler characteristic of the classifying
space of the homotopy colimit of a complex of groups in terms of
Euler characteristics of lower links and groups.

One novel aspect of our approach is that we do not require scwols to
be skeletal. We prove in
Theorem~\ref{the:reduction_to_skeletal_case} that any scwol with a
$G$-action in the sense of
Definition~\ref{def:group_action_on_a_scwol} can be replaced by a
skeletal scwol with a $G$-action, and this process preserves
quotients, stabilizers, complexes of groups, and homotopy colimits.
Moreover, if the initial $G$-action was free on the object set, then
so is the $G$-action on the object set of the skeletal replacement.

We begin by recalling the notions in Chapter III.$\calc$ of
Bridson--Haefliger~\cite{Bridson-Haefliger(1999)}, rephrased in the
conceptual framework of 2-category theory.

\begin{notation}[2-Category of
groups]\label{not:2-category_of_groups} We denote by $\GROUPS$ the
2-category of groups. Objects are groups and morphisms are group
homomorphisms. The 2-cells are given by conjugation: a 2-cell
$(g,a)$
$$\xymatrix{H \ar@/^1.5pc/[rr]_{\quad}^{a}="1"
\ar@/_1.5pc/[rr]_{a'}="2" & & G
\ar@{}"1";"2"|(.2){\,}="7" \ar@{}"1";"2"|(.8){\,}="8" \ar@{=>}"7"
;"8"^{(g,a)} }$$
is an element $g \in G$ such that $ga(h)g^{-1}=a'(h)$ for all $h \in H$. The vertical composition is
$(g_2,a_2)\odot (g_1,a_1)=(g_2g_1,a_1)$ and the horizontal composition of
$$\xymatrix@C=3pc{ H \ar@/^2pc/[rr]_{\quad}^{a}="1" \ar@/_2pc/[rr]_{a'}
="2" & & \ar@{}"1";"2"|(.2){\,}="7" \ar@{}"1";"2"|(.8){\,}="8"
\ar@{=>}"7" ;"8"^{(g,a)} G
\ar@/^2pc/[rr]_{\quad}^{b}="1"
\ar@/_2pc/[rr]_{b'}="2" & & \ar@{}"1";"2"|(.2){\,}="7"
\ar@{}"1";"2"|(.8){\,}="8" \ar@{=>}"7" ;"8"^{(k,b)}
K}$$
is $(kb(g),ba)$.
\end{notation}

\begin{definition}[Scwol]
A \emph{scwol}\footnote{Bridson--Haefliger additionally require
scwols to be skeletal \cite[page~574]{Bridson-Haefliger(1999)}.
However, we do not require scwols to be skeletal, since we prove in
Theorem~\ref{the:reduction_to_skeletal_case} that general statements
about scwols can be reduced to the skeletal case.} is a {\bf s}mall
{\bf c}ategory {\bf w}ith{\bf o}ut {\bf l}oops, that is, a small
category in which every endomorphism is trivial.
\end{definition}

\begin{example}
The categories $\{j \rightrightarrows k\}$ and $\cali=\{k \leftarrow
j \to \ell\}$ of Examples~\ref{exa:finite_model_for_parallel_arrows}
and \ref{exa:finite_model_for_pushout} are scwols.  Every partially
ordered set is a scwol, for example, the set of simplices of a
simplicial complex, ordered by the face relation, is a scwol. The
poset of non-empty subsets of $[q]$, and its opposite category in
Example~\ref{exa:finite_model_for_q_interior}, are scwols.  The
opposite category of a scwol is also a scwol.
\end{example}

\begin{lemma}  \label{lem:scwol_directly_finite_EI}
Every scwol is an EI-category and consequently also directly finite.
\end{lemma}
\begin{proof}
Every endomorphism in a scwol is trivial, and therefore an
automorphism, so every scwol is an EI-category. By
Fiore--L\"uck--Sauer
\cite[Lemma~3.13]{FioreLueckSauerFinObsAndEulCharOfCats(2009)},
every EI-category is also directly finite.

For a direct proof of direct finiteness: if $u \colon x \to y$ and
$v \colon y \to x$ are morphisms in a scwol, then $vu $ and $uv$ are
automorphisms, and hence both $vu = \id_x$ and $uv = \id_y$ hold
automatically.
\end{proof}

\begin{theorem}[Finite scwols admit finite models]\label{the:finite_models_for_finite_scwols}
Suppose $\cali$ is a finite scwol. Then $\cali$ admits a finite
$\cali$-$CW$-model for its $\cali$-classifying space in the sense of
Definition~\ref{def:classifying_calc-space}.
\end{theorem}
\begin{proof}
By Lemma~\ref{lem:finite_models_and_equivalences_of_categories}, we may assume that $\cali$ is skeletal.

Since $\cali$ has only finitely many morphisms, no nontrivial
isomorphisms, and no nontrivial endomorphisms, there are only
finitely many paths of non-identity morphisms. Thus the bar
construction of $E^{\barcon}\cali$ Remark~\ref{rem:Ebarcalc} has
only finitely many $\cali$-cells.
\end{proof}

\begin{corollary} \label{cor:finite_scwols_FF}
Any finite scwol $\cali$ is of types (FF$_R$) and (FP$_R$) for every
associative, commutative ring $R$ with identity.  Moreover, any
finite scwol is also of type ($L^2$).
\end{corollary}
\begin{proof}
The cellular $R$-chains of the finite model in
Theorem~\ref{the:finite_models_for_finite_scwols} provide a finite,
free resolution of the constant module $\underline{R}$.
By~Theorem~\ref{the:comparing_o_and_chi(2)},
any directly finite category of type (FP$_\IC$) is of type ($L^2)$.
Scwols are directly finite by
Lemma~\ref{lem:scwol_directly_finite_EI}.
\end{proof}

\begin{example}[Invariants coincide for finite scwols] \label{exa:Euler_characteristics_of_finite_scwols}
Let $\cali$ be any finite scwol. Then by
Corollary~\ref{cor:finite_scwols_FF} it is of type (FF$_R$) for any
associative, commutative ring with identity, and by
Theorems~\ref{the:chi_f_determines_chi}~and~\ref{the:coincidence},
we have
\begin{equation*}
\chi(\cali;R)=\chi(B\cali;R)=\chi^{(2)}(\cali).
\end{equation*}
If $\Gamma$ is any skeleton of $\cali$, then by \eqref{equ:Euler_characteristic_for_skeletal_finite_scwols},
\begin{equation} \label{equ:Euler_characteristic_given_by_paths}
\chi(\Gamma;R)=\sum_{n \geq 0}
(-1)^n c_n(\Gamma),
\end{equation}
where $c_n(\Gamma)$ is the number of paths of $n$-many non-identity
morphisms in $\Gamma$. But
by~Fiore--L\"uck--Sauer~\cite[Theorem~2.8~and~Corollary~4.19]{FioreLueckSauerFinObsAndEulCharOfCats(2009)},
type (FF$_R$) and the Euler characteristic are invariant under
equivalence of categories between directly finite categories, so
$\chi(\cali;R)=\chi(\Gamma;R)$ and all three invariants
$\chi(\cali;R)$, $\chi(B\cali;R)$, $\chi^{(2)}(\cali)$ are given
by~\eqref{equ:Euler_characteristic_given_by_paths}.
\end{example}

We now arrive at the main notion of this section: a complex of
groups. We will apply our Homotopy Colimit Formula to complexes of
groups.

\begin{definition}[Complex of groups] \label{def:complex_of_groups}
Let $\caly$ be a scwol. A \emph{complex of groups over $\caly$} is a
pseudo functor $F\colon \caly \to \GROUPS$ such that $F(a)$
is injective for every morphism $a$ in $\caly$. For each object
$\sigma$ of $\caly$, the group $F(\sigma)$ is called the \emph{local
group at $\sigma$}.
\end{definition}

In 2.5 and 2.4 of \cite{Haefliger(1991)} and \cite{Haefliger(1992)}
respectively, Haefliger denotes by $CG(X)$ the homotopy colimit of a
complex of groups $G(X)\colon C(X) \to \GROUPS$. Bridson--Haefliger
use the notation $CG(\caly)$ in
\cite[III.$\calc$.2.8]{Bridson-Haefliger(1999)}. The fundamental
group of a complex of groups $G(X)$ in the sense
of~\cite[Definition~3.5 on p.~548]{Bridson-Haefliger(1999)} equals
the fundamental group of the geometric realization of
$CG(X)$~\cite[Appendix A.12 on p.~578 and Remark A.14 on
p.~579]{Bridson-Haefliger(1999)}. Categories which are homotopy
colimits of complexes of groups are characterized by Haefliger on
page~283 of \cite{Haefliger(1992)}.  From the homotopy colimit
$CG(X)$, Haefliger reconstructs the category $C(X)$ and the complex
of groups $G(X)$ up to a coboundary on pages~282-283 of
\cite{Haefliger(1992)}. Every aspherical realization
\cite[Definition~3.3.4]{Haefliger(1992)} of $G(X)$ has the homotopy
type of the geometric realization of the homotopy colimit, denoted
$BG(X)$ \cite[page~296]{Haefliger(1992)}. The homotopy colimit also
plays a role in the homology and cohomology of complexes of groups
\cite[Section~4]{Haefliger(1992)}; a left $G(X)$-module is a functor
$CG(X)\to \ABELIANGROUPS$.

We return to our recollection of complexes of groups and examples
that arise from group actions.

\begin{definition}[Morphism from a complex of groups to a group]
A \emph{morphism from a complex of groups $F$ to a group $G$} is a pseudo natural
transformation $F \Rightarrow \Delta_G$, where $\Delta_G$ indicates the constant 2-functor $\caly \to \GROUPS$ with value $G$.
\end{definition}

A typical example of a complex of groups equipped with a morphism to
a group $G$ arises from an action of a group $G$ on a scwol, as we now explain.

\begin{definition}[Group action on a scwol, 1.11 of Bridson--Haefliger~\cite{Bridson-Haefliger(1999)}]
\label{def:group_action_on_a_scwol} An \emph{action of a group $G$
on a scwol $\calx$} is a group homomorphism from $G$ into the group
of strictly invertible functors $\calx \to \calx$ such that
\begin{enumerate}
\item \label{def:group_action_on_a_scwol:(i)}
For every nontrivial morphism $a$ of $\calx$ and every $g \in G$, we
have $gs(a) \neq t(a)$,
\item \label{def:group_action_on_a_scwol:(ii)}
For every nontrivial morphism $a$ of $\calx$ and every $g \in G$, if
$gs(a)=s(a)$, then $ga=a$.
\end{enumerate}
\end{definition}

\begin{example} \label{exa:Z2_action_on_circle}
The group $G=\IZ_2$ acts in the sense of
Definition~\ref{def:group_action_on_a_scwol} on the scwol $\calx$
pictured below.
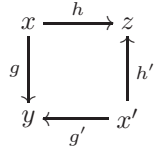
\[
\xymatrix{x \ar[r]^{h} \ar[d]_{g} & z \\ y & x' \ar[l]^{g'}
\ar[u]_{h'}}
\]
The group $\IZ_2$ permutes respectively $x$ and $x'$, $g$ and $g'$,
and $h$ and $h'$. The objects $y$ and $z$ are fixed. This action of
$\IZ_2$ on $\calx$ is a combinatorial model for a reflection action
on $S^1$.
\end{example}

\begin{example} \label{exa:Z2_Z_action_on_line}
Consider the scwol $\calx$ pictured below. The group $G=\{\pm 1
\}\ltimes \IZ$ acts on $\calx$ in the sense of
Definition~\ref{def:group_action_on_a_scwol} where $-1 \cdot m:=-m$
and $n\cdot m:=m +2n$.
\[
\xymatrix@C=2pc{  \cdots \ar[r] & -2 &  -1 \ar[r]  \ar[l] & 0   & 1
\ar[r] \ar[l] & 2 & \cdots \ar[l] }
\]
This action of $\{\pm 1 \}\ltimes \IZ$ on $\calx$ is a combinatorial
model for the reflection and translation action on $\mathbb{R}$.
\end{example}

\begin{lemma}[Consequences of group action conditions] \label{lem:consequences_of_group_action}
If a group $G$ acts on a scwol $\calx$ in the sense of
Definition~\ref{def:group_action_on_a_scwol}, then the following
statements hold.
\begin{enumerate}
\item \label{lem:consequences_of_group_action:(i)}
If $\sigma$ is an object of $\calx$ and $g, h \in G$, then $g \sigma \cong h \sigma$ implies $g \sigma = h \sigma$.
\item \label{lem:consequences_of_group_action:(ii)}
If $a$ is a morphism in $\calx$ and $g,h\in G$, then $gs(a)=hs(a)$
implies $ga=ha$.
\item \label{lem:consequences_of_group_action:(iii)}
If $\sigma\cong\tau$, then the stabilizers $G_\sigma$ and $G_\tau$
are equal.
\end{enumerate}
\end{lemma}
\begin{proof}
For statement~\ref{lem:consequences_of_group_action:(i)}, $g \sigma \cong h \sigma$ implies $\sigma \cong (g^{-1}h) \sigma$, so $\sigma = (g^{-1}h) \sigma$ by Definition~\ref{def:group_action_on_a_scwol} part \ref{def:group_action_on_a_scwol:(i)}, and $g\sigma =h \sigma$.

For statement~\ref{lem:consequences_of_group_action:(ii)}, $gs(a)=hs(a)$ implies $(h^{-1} g) s(a) = s(a)$ and $(h^{-1} g) a = a$ by Definition~\ref{def:group_action_on_a_scwol} part \ref{def:group_action_on_a_scwol:(ii)}, and finally $ga=ha$.

For statement~\ref{lem:consequences_of_group_action:(iii)}, suppose
$\sigma\cong\tau$ and $g\sigma=\sigma$. We have
\[
\tau \cong \sigma = g\sigma \cong g\tau.
\]
Then $\tau=g\tau$ by \ref{lem:consequences_of_group_action:(i)}, and
$G_\sigma \subseteq G_\tau$. The proof is symmetric, so we also have
$G_\tau \subseteq G_\sigma$.
\end{proof}

\begin{definition}[Quotient of a scwol by a group action] \label{def:quotient_scwol}
If a scwol $\calx$ is equipped with a $G$-action as above,
then the \emph{quotient scwol} $\calx / G$ has objects and morphisms
$$\ob (\calx /G ) := (\ob (\calx)) /G  $$
$$\mor (\calx /G ) := (\mor (\calx)) /G.$$
Composition and identities are induced by those of $\calx$.
\end{definition}

\begin{remark}[III.$\calc$.1.13 of Bridson--Haefliger~\cite{Bridson-Haefliger(1999)}]\label{rem:projection_functor_for_group_action}
  The projection functor $p\colon \calx \to \calx / G$ induces a bijection
  \begin{equation} \label{equ:projection_covering_1} \xymatrix{\{a \in
      \mor(\calx) \vert sa=x\} \ar[r] & \{b \in \mor(\calx / G)\vert sb=p(x) \}}
  \end{equation}
  for each $x \in \calx$. If $G / \calx$ is connected and the action of $G$ on
  $\ob(\calx)$ is free, then $p$ is a \emph{covering of scwols}. That is, in
  addition to the bijection \eqref{equ:projection_covering_1}, $p$ induces a
  bijection
  \begin{equation} \label{equ:projection_covering_2} \xymatrix{\{a \in
      \mor(\calx) \vert ta=x\} \ar[r] & \{b \in \mor(\calx / G)\vert tb=p(x) \}}
  \end{equation}
  for each $x \in \calx$.
\end{remark}

\begin{lemma}[Quotients of skeletal scwols are skeletal] \label{lem:quotient_scwol_of_skeletal_scwol_is_skeletal}
If $\calx$ is a skeletal scwol, and a group $G$ acts on $\calx$ in
the sense of Definition~\ref{def:group_action_on_a_scwol}, then the
quotient scwol $\calx / G$ is also skeletal.
\end{lemma}
\begin{proof}
Suppose $\overline{\sigma}$ is isomorphic to $\overline{\tau}$ in
$\calx / G$. We show $\overline{\sigma}$ is actually equal to
$\overline{\tau}$. If $\overline{a}\colon \overline{\sigma} \to
\overline{\tau}$ is an isomorphism with inverse $\overline{b}$, then
there are lifts $a\colon \sigma \to \tau$ and $b \colon \tau \to
\sigma'$ in $\calx$, and an element $g \in G$ such that
$g(ba)=\id_\sigma$. Since $g$ fixes the source of $ba$, the group
element $g$ fixes also $ba$, so $ba=\id_\sigma$ and
$\sigma'=\sigma$. Since $ab$ is an endomorphism of $\tau$, it is
therefore $\id_\tau$. By the skeletality of $\calx$, we have $\sigma
= \tau$, and also $\overline{\sigma}=\overline{\tau}$.
\end{proof}

\begin{lemma}[Quotient of path set is set of paths in quotient] \label{lem:cells_of_quotient_scwol_are_quotient_of_scwol_cells}
Suppose $\calx$ is a scwol equipped with an action of a group $G$ in
the sense of Definition~\ref{def:group_action_on_a_scwol}. Let
$\Lambda_n(\calx)$ respectively $\Lambda_n(\calx/G)$ denote the set
of paths of $n$-many non-identity composable morphisms in $\calx$
respectively $\calx/G$. Give $\Lambda_n(\calx)$ the induced
$G$-action. Then the function
\[
\Lambda_n(\calx) \to \Lambda_n(\calx/G)
\]
\[
 (a_1,\dots,a_n) \mapsto (\overline{a}_1, \dots, \overline{a}_n)
\]
induces a bijection $\Lambda_n(\calx)/G\to \Lambda_n(\calx/G)$.
\end{lemma}
\begin{proof}
Remark~\ref{rem:projection_functor_for_group_action} implies that a
path $(a_1,\dots,a_n)$ in $\calx$ consists entirely of non-identity
morphisms if and only if the projection $(\overline{a}_1, \dots,
\overline{a}_n)$ in $\calx/G$ consists entirely of non-identity
morphisms, so from now on we work only with non-identity morphisms.
Note
\[
(g_1a_1,g_2a_2,\dots,g_na_n)=(g_1a_1,g_1a_2,\dots,g_1a_n)
\]
by
Definition~\ref{def:group_action_on_a_scwol}~\ref{def:group_action_on_a_scwol:(ii)}.
For injectivity, we have $(\overline{a}_1, \dots,
\overline{a}_n)=(\overline{b}_1, \dots, \overline{b}_n)$ if and only
if for some $g_i \in G$
\[
(g_1a_1,g_2a_2,\dots,g_na_n)=(b_1,\dots,b_n),
\]
which happens if and only if for some $g \in G$
\[
(ga_1,ga_2,\dots,ga_n)=(b_1,\dots,b_n),
\]
(take $g=g_1$). For the surjectivity, we can lift any path
$(\overline{a}_1, \dots, \overline{a}_n)$ by first lifting
$\overline{a}_1$ to $a_1$, then $\overline{a}_2$ to $a_2$, and so on
using Remark~\ref{rem:projection_functor_for_group_action}.
\end{proof}

\begin{definition}[Complex of groups from a group
action on a scwol, 2.9 of
Bridson--Haefliger~\cite{Bridson-Haefliger(1999)}]\label{def:complex_of_groups_from_a_group_action}
Let $G$ be a group and $\calx$ a scwol upon which $G$ acts in the
sense of Definition~\ref{def:group_action_on_a_scwol}. Let $p \colon
\calx \to \calx/G$ denote the quotient map.

Haefliger and Bridson--Haefliger define a pseudo functor $F\colon
\calx/G \to \GROUPS$ as follows. In the procedure choices are made,
but different choices lead to isomorphic complexes of groups.  For
each object $\overline{\sigma}$ of $\calx/G$, choose an object
$\sigma$ of $\calx$ such that $p(\sigma)=\overline{\sigma}$ (our
overline convention is the opposite of that in
\cite{Bridson-Haefliger(1999)}). Then $F(\overline{\sigma})$ is
defined to be $G_\sigma$, the isotropy group of $\sigma$ under the
$G$-action.

If $\overline{a}\colon \overline{\sigma} \to \overline{\tau}$ is a
morphism in $\calx/G$, then there exists a unique morphism $a$ in
$\calx$ such that $p(a)=\overline{a}$ and $sa=\sigma$, as in
\eqref{equ:projection_covering_1}. For $\overline{a}$ we choose an
element $h_{\overline{a}} \in G$ such that $h_{\overline{a}} \cdot
ta$ is the object $\tau$ of $\calx$ chosen above so
that $p(\tau)=\overline{\tau}$. An injective group
homomorphism $F(\overline{a})\colon G_\sigma \to G_\tau$ is defined by
$$F(\overline{a})(g):=h_{\overline{a}} g h_{\overline{a}}^{-1}.$$

Suppose $\overline{a}$ and $\overline{b}$ are composable morphisms
of $\calx/G$. We define a 2-cell in $\GROUPS$
$$F_{\overline{b},\overline{a}}\colon F(\overline{b}) \circ F(\overline{a}) \Rightarrow
F(\overline{b} \circ \overline{a})$$ to be
$(h_{\overline{b}\overline{a}}h_{\overline{a}}^{-1}
h_{\overline{b}}^{-1} ,F(\overline{b}) \circ F(\overline{a}))$ as in
Notation~\ref{not:2-category_of_groups}.

The pseudo functor $F\colon \calx/G \to \GROUPS$ is called the {\it
complex of groups associated to the group action of $G$ on the scwol
$\calx$}. This complex of groups comes equipped with a morphism to the group $G$, that is, a pseudo natural
transformation $F \Rightarrow \Delta_G$. The inclusion of each isotropy group $F(\overline{\sigma})=G_\sigma$ into $G$ provides the components of the
pseudo natural transformation.
\end{definition}

\begin{example} \label{exa:Z2_complex_of_groups}
The quotient scwols for the actions in
Examples~\ref{exa:Z2_action_on_circle}~and~\ref{exa:Z2_Z_action_on_line}
are both $\{k \leftarrow j \rightarrow \ell \}$, and the associated
complexes of groups are both
\[
\xymatrix{\IZ_2 & \{0\} \ar[l] \ar[r] &\IZ_2.}
\]
\end{example}

\begin{remark} \label{rem:finite_group_finite_scwol_imply_hocolim_hypotheses_satisfied}
If a group $G$ acts on a scwol in the sense of
Definition~\ref{def:group_action_on_a_scwol}, each object stabilizer
is finite, and the quotient scwol is finite, then the associated
complex of groups $F\colon \calx/G \to \GROUPS$ satisfies all of the
hypotheses of the Homotopy Colimit Formula in
Theorem~\ref{the:homotopy_colimit_formula}~\ref{the:homotopy_colimit_formula:chi(2)}
and in Corollary~\ref{cor:homotopy_colimit_formula_for_pseudo_functors}~\ref{the:homotopy_colimit_formula:chi(2)}. If, in addition, $R$ is a ring such that the order $\vert H\vert $
of each object stabilizer $H\subset G$ is invertible in $R$, then
$F\colon \calx/G \to \GROUPS$ also satisfies all of the
hypotheses of the Homotopy Colimit Formula in
Theorem~\ref{the:homotopy_colimit_formula}~\ref{the:homotopy_colimit_formula:chi}
and in Corollary~\ref{cor:homotopy_colimit_formula_for_pseudo_functors}~\ref{the:homotopy_colimit_formula:chi}.
See
Examples~\ref{exa:Z2_action_on_circle},~\ref{exa:Z2_Z_action_on_line},~and~\ref{exa:Z2_complex_of_groups}.
\end{remark}

Even without finiteness assumptions, it is possible to replace scwols with skeletal scwols and preserve much of the accompanying structure, as Theorem~\ref{the:reduction_to_skeletal_case} explains.

\begin{theorem}[Reduction to skeletal case] \label{the:reduction_to_skeletal_case}
Let $G$ be a group acting on a scwol $\calx$ in the sense of Definition~\ref{def:group_action_on_a_scwol}. Let $\Gamma$ be any skeleton of $\calx$, $i\colon \Gamma \to \calx$ the inclusion, and $r \colon \calx \to \Gamma$ a functor equipped with a natural isomorphism $ir \cong \id_\calx$, and  satisfying $r i = \id_\Gamma$. Then there is a $G$-action on the scwol $\Gamma$ in the sense of Definition~\ref{def:group_action_on_a_scwol} such that following hold.
\begin{enumerate}
\item \label{the:reduction_to_skeletal_case:r_equivariant}
The functor $r$ is $G$-equivariant.
\item \label{the:reduction_to_skeletal_case:r_induces_equivalence_of_quotients}
The induced functor $\overline{r}$ on quotient categories is an equivalence of categories compatible with the quotient maps, that is, the diagram below commutes.
\begin{equation} \label{equ:the:reduction_to_skeletal_case:rbar_compatible_with_projections}
\xymatrix{\calx \ar[r]^r \ar[d]_{p^\calx} & \Gamma \ar[d]^{p^{\Gamma}} \\
\calx/G \ar[r]_{\overline{r}} & \Gamma/G}
\end{equation}
\item \label{the:reduction_to_skeletal_case:inclusion_preserves_stabilizers}
The inclusion $i\colon \Gamma \to \calx$ preserves stabilizers, that is $G_{i\gamma}=G_{\gamma}$ for all $\gamma \in \ob(\Gamma)$. Note that the inclusion may not be $G$-equivariant.
\item \label{the:reduction_to_skeletal_case:complexes_of_groups_agree}
Choices can be made in the definitions of $F^\calx$ and $F^\Gamma$ (the complexes of groups associated to the $G$-actions on $\calx$ and $\Gamma$ in Definition~\ref{def:complex_of_groups_from_a_group_action}), so that the diagram below strictly commutes.
\begin{equation} \label{equ:the:reduction_to_skeletal_case:complexes_of_groups_agree}
\xymatrix{\calx/G \ar[rr]^{\overline{r}} \ar[dr]_{F^\calx} & &
\Gamma/G \ar[dl]^{F^\Gamma} \\ & \GROUPS & }
\end{equation}
\item \label{the:reduction_to_skeletal_case:homotopy_colimits_equivalent}
The functor $(\overline{r},\id)$ is an equivalence of categories
\[
\xymatrix{(\overline{r},\id)\colon \hocolim_{\calx/G} F^\calx \ar[r] & \hocolim_{\Gamma/G} F^\Gamma.}
\]
\item \label{the:reduction_to_skeletal_case:free_implies_free}
If $G$ acts freely on $\ob(\calx)$, then $G$ acts freely on $\ob(\Gamma)$.
\end{enumerate}
\end{theorem}
\begin{proof}
To define the group action, let $\Aut(\calx)$ and $\Aut(\Gamma)$ denote the
strictly invertible endofunctors on $\calx$ and $\Gamma$
respectively, and consider the monoid homomorphism
\begin{equation} \label{eqn:action_on_skeleton}
\varphi \colon \Aut(\calx) \to \End(\Gamma), \;\; F \mapsto r \circ F \circ i.
\end{equation}
This is strictly multiplicatively because the natural isomorphism of functors
\begin{eqnarray*}
r \circ G \circ F \circ i & = & r \circ G \circ \id_{\calx} \circ F \circ i \\
& \cong & (r \circ G \circ i) \circ (r  \circ F \circ i),
\end{eqnarray*}
and skeletality of $\Gamma$
imply $\varphi(GF)$ agrees with $\varphi(G)\varphi(F)$ on objects of $\Gamma$, so each component $\varphi(GF)(\gamma)\cong\varphi(G)\varphi(F)(\gamma)$ is an endomorphism in the scwol $\Gamma$, and is therefore trivial. By naturality, $\varphi(GF)$ and $\varphi(G)\varphi(F)$ agree on morphisms also.
Consequently, $\varphi$ takes values in $\Aut(\Gamma)$ and is a homomorphism
$\varphi \colon \Aut(\calx) \to \Aut(\Gamma)$.

We define a $G$-action on $\Gamma$ as the composite of the action $G
\to \Aut(\calx)$ with $\varphi$ in \eqref{eqn:action_on_skeleton}. We indicate the action of $g$ on $\Gamma$ by $\varphi(g)\gamma$ and the action of $g$ on $\calx$ by $gx$. For simplicity, we suppress $i$ from the notation when indicating the $G$-action in $\calx$ on objects and morphisms of $\Gamma$, so
for example, if $a$ is morphism in $\Gamma$, then $gs(a)$ actually means $gis(a)$ throughout.

To verify Definition~\ref{def:group_action_on_a_scwol}~\ref{def:group_action_on_a_scwol:(i)}
for $\Gamma$, suppose $a$ is a nontrivial morphism in $\Gamma$ and $\varphi(g)s(a)=t(a)$, that is $rgs(a)=t(a)$. Then $gs(a)
\cong t(a)$ in $\calx$, but $gs(a) \neq t(a)$ (for if $gs(a)=t(a)$, then $a$ must be trivial by Definition~\ref{def:group_action_on_a_scwol}~\ref{def:group_action_on_a_scwol:(i)} for $\calx$). Let $b\colon t(a) \to gs(a)$ be an isomorphism in $\calx$ and consider the composite $ba \colon s(a) \to t(a)
\to gs(a)$. Then $gs(ba)=gs(a)=t(ba)$, so $ba$ must be trivial by
Definition~\ref{def:group_action_on_a_scwol}~\ref{def:group_action_on_a_scwol:(i)}
for $\calx$. Consequently $a=b^{-1}$ is a nontrivial \emph{iso}morphism in
$\Gamma$, and we have a contradiction to either skeletality or the
no loops requirement. Thus $\varphi(g)s(a)\neq t(a)$, and
Definition~\ref{def:group_action_on_a_scwol}~\ref{def:group_action_on_a_scwol:(i)}
holds for $\Gamma$. The verification of
Definition~\ref{def:group_action_on_a_scwol}~\ref{def:group_action_on_a_scwol:(ii)}
is shorter: if $a$ is a nontrivial morphism in $\Gamma$ and $\varphi(g)s(a)=s(a)$, that is $rgs(a)=s(a)$, then $gs(a)\cong s(a)$, and $gs(a)=s(a)$ by
Lemma~\ref{lem:consequences_of_group_action}~\ref{lem:consequences_of_group_action:(i)}
for $\calx$. Finally, $ga=a$ by Definition~\ref{def:group_action_on_a_scwol}~\ref{def:group_action_on_a_scwol:(ii)} for $\calx$, $rga=a$ as $a$ is in $\Gamma$, and $\varphi(g)a=a$. The action of $G$ on
$\Gamma$ satisfies Definition~\ref{def:group_action_on_a_scwol} and
we may form the quotient scwol $\Gamma/G$ as in Defition~\ref{def:quotient_scwol}, which is skeletal by Lemma~\ref{lem:quotient_scwol_of_skeletal_scwol_is_skeletal}.
\\[1mm]\ref{the:reduction_to_skeletal_case:r_equivariant} For the $G$-equivariance of $r$, let $f\colon x \to y$ be a morphism in $\calx$ and consider the naturality diagram.
\[
\xymatrix@C=6pc{rgirx \ar[r]^{rgirf=\varphi(g)r(f)} \ar[d]_\cong & rgiry \ar[d]^\cong \\ rgx \ar[r]_{rgf} & rgy}
\]
The vertical morphisms must be identities by skeletality of $\Gamma$ and the no loops condition, so $\varphi(g)r(f)=r(gf)$. Equivariance on objects then follows by taking $f=\id_x$.
\\[1mm]\ref{the:reduction_to_skeletal_case:r_induces_equivalence_of_quotients}
Diagram \eqref{equ:the:reduction_to_skeletal_case:rbar_compatible_with_projections} commutes by definition of $\overline{r}$. The functor $\overline{r}$ is surjective on objects because $p^\Gamma r$ and $p^\calx$ are. The functor $\overline{r}$ is fully faithful since the equivariant bijection $r(x,y)\colon \mor_\calx(x,y) \to \mor_\Gamma(r(x),r(y))$ induces the equivariant bijection $\overline{r}(p^\calx x, p^\calx y)$.
\\[1mm]\ref{the:reduction_to_skeletal_case:inclusion_preserves_stabilizers}
Let $\gamma \in \ob(\Gamma)$, and suppose  $gi\gamma=i\gamma$. Then
\begin{eqnarray*}
\varphi(g)\gamma & \overset{\text{def}}{=} & r(g i\gamma)  \\
& = & r (i \gamma) \\
& = & \gamma
\end{eqnarray*}
and $G_{i\gamma} \subseteq G_\gamma$. Now suppose $\varphi(g)\gamma =\gamma$.
Then $r(gi\gamma)=\gamma$ by definition, and $g i \gamma
\cong i\gamma$ in $\calx$, which says $g\cdot i \gamma = i\gamma$ by
Lemma~\ref{lem:consequences_of_group_action}~\ref{lem:consequences_of_group_action:(i)},
and $G_\gamma \subseteq G_{i\gamma} $.
\\[1mm]\ref{the:reduction_to_skeletal_case:complexes_of_groups_agree}
We claim that choices can be made in the definitions of the
associated complexes of groups $F^\calx$ and $F^\Gamma$
(see~Definition~\ref{def:complex_of_groups_from_a_group_action}) so
that diagram
\eqref{equ:the:reduction_to_skeletal_case:complexes_of_groups_agree}
strictly commutes. First choose a skeleton $\calq$ of the quotient
$\calx/G$, define $F^\calx$ on object in the skeleton $\calq$, and
then extend to all objects in $\calx/G$. For every $\overline{q} \in
\ob(\calq)$, select a $q \in \ob(\calx)$ such that
$p^\calx(q)=\overline{q}$ and define $F^\calx(\overline{q})=G_q$. We
remain with the choice of the selected preimage $q$ of
$\overline{q}$ throughout. If $\overline{\sigma} \in \ob(\calx/G)$
and $\overline{a}\colon \overline{q} \cong \overline{\sigma}$ is an
isomorphism in $\calx/G$, then also define
$F^\calx(\overline{\sigma})=G_q$. This is allowed, since
$\overline{a}\colon \overline{q}\cong \overline{\sigma}$ implies
existence of morphisms $a\colon q \to g_\sigma \sigma$ and $b \colon
\sigma \to g_qq$ in $\calx$, and the composite
\[
\xymatrix{q \ar[r]^{a} & g_\sigma\sigma \ar[r]^{g_\sigma b} &
g_\sigma g_q q}
\]
is trivial by
Definition~\ref{def:group_action_on_a_scwol}~\ref{def:group_action_on_a_scwol:(i)}.
The opposite composite is also trivial, as it is a loop, and we have
$q \cong g_\sigma \sigma$ in $\calx$. Then by
Lemma~\ref{lem:consequences_of_group_action}~\ref{lem:consequences_of_group_action:(iii)},
$G_q = G_{g_\sigma \sigma}$ and we may define
$F^\calx(\overline{\sigma})=G_q$ because $p^{\calx}(g_\sigma
\sigma)=\overline{\sigma}$. In particular, the selected preimage of
$\overline{\sigma}$ in $\calx$ is $g_\sigma \sigma$ and we select
$h_{\overline{a}}=e_G$ for $\overline{a}\colon \overline{q}\cong
\overline{\sigma}$ in
Definition~\ref{def:complex_of_groups_from_a_group_action}, so
$F^\calx(\overline{a})=\id_{G_q}$.  We remark that the isomorphism $\overline{a}$ is
the only morphism $\overline{q} \to \overline{\sigma}$ because there
are no loops in $\calx/G$, so the element $g_\sigma \sigma$ is
uniquely defined as the target of the unique morphism $a$ with
source $q$ and $p^\calx$-image $\overline{a}$.

We next define $F^\Gamma$ on objects of $\Gamma/G$ using the
equivalence $\overline{r}$ and the definition of $F^\calx$ on
objects of $\calq$. For $\overline{q} \in \ob(\calq)$, we also
define $F^\Gamma(\overline{r}(\overline{q}))=G_q$. This is allowed:
for $\overline{r}(\overline{q})=\overline{r(q)}$ we choose $r(q)$ as
the selected preimage in $\ob (\Gamma)$, and $ir(q) \cong q$ in
$\calx$, so $G_{r(q)}=G_{ir(q)}=G_q$ by
\ref{the:reduction_to_skeletal_case:inclusion_preserves_stabilizers}
and
Lemma~\ref{lem:consequences_of_group_action}~\ref{lem:consequences_of_group_action:(iii)}.
Every $\overline{\gamma} \in \ob (\Gamma/G)$ is of the form
$\overline{r}(\overline{q})$ for a unique $\overline{q} \in \calq$,
so $F^\Gamma$ is now defined on all objects of $\Gamma/G$, and we
have $F^\Gamma \circ \overline{r}=F^\calx$ on all objects of
$\calx/G$.

We must now define $F^\calx$ and $F^\Gamma$ on morphisms so that
$F^\Gamma \circ \overline{r}=F^\calx$ for morphisms also. The idea
is to first define $F^\calx$ on morphisms in the skeleton $\calq$,
then extend to all of $\calx/G$, and then define $F^\Gamma$ on morphisms of $\Gamma/G$. If
$\overline{a}\colon \overline{q}_1 \to \overline{q}_2$ is a morphism
in $\calq$, then there is a unique morphism $a$ in $\calx$ with
source $q_1$ and $p^\calx(a)=\overline{a}$. Select any
$h_{\overline{a}}$ such that $h_{\overline{a}}ta=q_2$. Then we
define an injective group homomorphism $F(\overline{a})\colon
G_{q_1} \to G_{q_2}$  by
\[
F(\overline{a})(g):=h_{\overline{a}} g h_{\overline{a}}^{-1}.
\]
If $\overline{b}\colon \overline{\sigma}_1 \to \overline{\sigma}_2$
is any morphism in $\calx/G$, then there exists a unique
$\overline{a}$ in $\calq$ and a unique commutative diagram with
vertical isomorphisms as below.
\[
\xymatrix{\overline{q}_1 \ar[r]^{\overline{a}} \ar[d]_{\cong} &
\overline{q}_2 \ar[d]^{\cong} \\ \overline{\sigma}_1
\ar[r]_{\overline{b}} & \overline{\sigma}_2 }
\]
Then we choose $h_{\overline{b}}$ to be $h_{\overline{a}}$, and we
consequently have $F(\overline{a})=F(\overline{b})$. If
$\overline{c}\colon \overline{r}(\overline{q}_1) \to
\overline{r}(\overline{q}_2)$ is a morphism in $\Gamma/G$, then
there is a unique $\overline{a}\colon \overline{q}_1 \to
\overline{q}_2$ in $\calq$ with
$\overline{r}(\overline{a})=\overline{c}$ and we choose
$h_{\overline{c}}$ to be $h_{\overline{a}}$. Manifestly, we have
$F^\Gamma \circ \overline{r}=F^\calx$. The coherences of $F^\calx$
and $F^\Gamma$ are also compatible, since they are determined by the
$h_{\overline{a}}$'s.
\\[1mm]\ref{the:reduction_to_skeletal_case:homotopy_colimits_equivalent}
From~\ref{the:reduction_to_skeletal_case:r_induces_equivalence_of_quotients} we know $\overline{r}$ is a surjective-on-objects equivalence of categories and from~\ref{the:reduction_to_skeletal_case:complexes_of_groups_agree} we have $F^\calx=F^\Gamma \circ \overline{r}$. From this, one sees
\[
\xymatrix{(\overline{r},\id)\colon \hocolim_{\calx/G} F^\calx = \hocolim_{\calx/G} F^\Gamma \circ \overline{r} \ar[r] & \hocolim_{\Gamma/G} F^\Gamma}
\]
is an equivalence of categories.
\\[1mm]\ref{the:reduction_to_skeletal_case:free_implies_free}
If the action of $G$ on $\ob(\calx)$ is free, then for each $\gamma \in \ob(\Gamma)$, the group
$G_\gamma=G_{i\gamma}$ (see \ref{the:reduction_to_skeletal_case:inclusion_preserves_stabilizers}) is trivial, and $G$ acts freely on $\ob(\Gamma)$.
\end{proof}

\begin{remark}
In Theorem~\ref{the:reduction_to_skeletal_case}, it is even possible
to select a skeleton so that the inclusion is $G$-equivariant,
though we will not need this. See Section~\ref{sec:appendix}.
\end{remark}

In~\cite[Theorems~5.30~and~5.37]{FioreLueckSauerFinObsAndEulCharOfCats(2009)},
we proved the compatibility of the $L^2$-Euler characteristic with
coverings and isofibrations of finite connected groupoids.
Theorem~\ref{the:compatibility_of_chi_with_free_group_action_on_finite_scowls}
is an analogue for scwols (see
Remark~\ref{rem:projection_functor_for_group_action}).

\begin{theorem}[Compatibility with free actions on finite scwols] \label{the:compatibility_of_chi_with_free_group_action_on_finite_scowls}
Let $G$ be a finite group acting on a finite scwol $\calx$. If $G$
acts freely on $\ob(\calx)$, then
\[
\chi(\calx/G;R)=\frac{\chi(\calx;R)}{\vert G \vert}  \;\; \text{ and }
\; \; \chi^{(2)}(\calx/G)=\frac{\chi^{(2)}(\calx)}{\vert G \vert}.
\]
Recall $\chi(-;R)$ and $\chi^{(2)}(-)$ agree for finite scwols by Example~\ref{exa:Euler_characteristics_of_finite_scwols}.
\end{theorem}
\begin{proof}
By
Theorem~\ref{the:reduction_to_skeletal_case}~\ref{the:reduction_to_skeletal_case:r_equivariant},
\ref{the:reduction_to_skeletal_case:r_induces_equivalence_of_quotients},
and \ref{the:reduction_to_skeletal_case:free_implies_free}, we may
assume $\calx$ is skeletal.

A consequence of
Definition~\ref{def:group_action_on_a_scwol}~\ref{def:group_action_on_a_scwol:(ii)}
(independent of skeletality) is that an element $g \in G$ fixes a
path $a=(a_1,\dots,a_n)$ in $\calx$ if and only if $g$ fixes $sa_1$,
so $G_{sa_1}=G_a$. Then $G$ acts freely on $\Lambda_n(\calx)$, since
it acts freely on $\ob(\calx)$.

The scwol $\calx/G$ is skeletal by Lemma~\ref{lem:quotient_scwol_of_skeletal_scwol_is_skeletal}, and by Example~\ref{exa:Euler_characteristics_of_finite_scwols} and Lemma~\ref{lem:cells_of_quotient_scwol_are_quotient_of_scwol_cells} we have
\begin{eqnarray*}
\chi^{(2)}(\calx/G) & = & \sum_{n \geq 0} (-1)^n c_n(\calx/G) \\
& = & \sum_{n \geq 0} (-1)^n \vert\Lambda_n(\calx/G) \vert \\
& = & \sum_{n \geq 0} (-1)^n \vert\Lambda_n(\calx)/G \vert \\
& = & \sum_{n \geq 0} (-1)^n \frac{\vert\Lambda_n(\calx)\vert}{\vert G\vert} \\
& = & \frac{1}{\vert G\vert}\sum_{n \geq 0} (-1)^n \vert\Lambda_n(\calx)\vert \\
& = & \frac{1}{\vert G\vert}\sum_{n \geq 0} (-1)^n c_n(\calx) \\
& = & \frac{\chi^{(2)}(\calx)}{\vert G \vert}.
\end{eqnarray*}
\end{proof}

A complex of groups is called \emph{developable} if it is isomorphic
to a complex of groups associated to a group action. A classical
theorem of Bass--Serre says that every complex of groups on a scwol
with maximal path length 1 is developable. The following gives a
necessary condition of developability of a complex of groups from a
scwol and group of specified Euler characteristics.

\begin{theorem}[Euler characteristics of associated complexes of groups]
\label{thm:Euler_characteristic_of_hocolim_of_quotient_complex}
Let $G$ be a finite group that acts on a finite scwol $\calx$ in the sense of
Definition~\ref{def:group_action_on_a_scwol}. Let $F\colon \calx/G
\to \GROUPS$ be the associated complex of groups. Then
\[
\chi^{(2)}(\hocolim_{\calx/G} F )= \frac{\chi^{(2)}(\calx)}{\vert
G\vert} = \frac{\chi(\calx;\IC)}{\vert G\vert} =
\frac{\chi(B\calx;\IC)}{\vert G\vert}.
\]
If $R$ is a ring such that the orders of subgroups $H\subset G$ are invertible in
$R$, then we also have
\[
\chi(\hocolim_{\calx/G} F;R)=\chi(\calx/G;R).
\]
\end{theorem}
\begin{proof}
By
Theorem~\ref{the:reduction_to_skeletal_case}~\ref{the:reduction_to_skeletal_case:r_equivariant},
\ref{the:reduction_to_skeletal_case:r_induces_equivalence_of_quotients},
\ref{the:reduction_to_skeletal_case:complexes_of_groups_agree}, and
\ref{the:reduction_to_skeletal_case:homotopy_colimits_equivalent},
we may assume $\calx$ is skeletal. Then $\calx/G$ is also skeletal
by~Lemma~\ref{lem:quotient_scwol_of_skeletal_scwol_is_skeletal}.

Let $\Lambda_n(\calx)$ respectively $\Lambda_n(\calx/G)$ denote the
set of paths of $n$-many non-identity composable morphisms in
$\calx$ respectively $\calx/G$. Then by
Lemma~\ref{lem:cells_of_quotient_scwol_are_quotient_of_scwol_cells},
the sets $\Lambda_n(\calx)/G$ and $\Lambda_n(\calx/G)$ are in
bijective correspondence.

We will also use the fact that an element $g \in G$ fixes a path
$a=(a_1,\dots,a_n)$ in $\calx$ if and only if $g$ fixes $sa_1$, so
$G_{sa_1}=G_a$. This is a consequence of
Definition~\ref{def:group_action_on_a_scwol}~\ref{def:group_action_on_a_scwol:(ii)}.

By Theorem~\ref{the:finite_models_for_finite_scwols},
$E^{\barcon}\calx$ and $E^{\barcon}( \calx / G )$ are finite models
for the skeletal scwols $\calx$ and $\calx/G$, and the $n$-cells are
indexed by $\Lambda_n(\calx)$ and $\Lambda_n(\calx/G)$,
respectively. For each path $(a_1,\dots,a_n)$ in $\calx$, there is
an $n$-cell in $E^{\barcon}\calx$ based at $sa_1$. A similar
statement holds for $\calx/G$ and $E^{\barcon}( \calx / G )$.

Now we may apply the Homotopy Colimit Formula to the associated
complex of groups $F:\calx / G \to \GROUPS$ by
Remark~\ref{rem:finite_group_finite_scwol_imply_hocolim_hypotheses_satisfied}.
For the Euler characteristic, we have
\begin{eqnarray*}
\chi(\hocolim_{\calx/G} F;R) & = & \sum_{n \geq 0} (-1)^n \cdot \left(
\sum_{\overline{a} \in \Lambda_n(\calx/G)} \chi(F(s\overline{a}_1);R)
\right) \\ & = & \sum_{n \geq 0} (-1)^n \cdot \left(
\sum_{\overline{a} \in \Lambda_n(\calx/G)} 1\right)
\\ & = & \sum_{n \geq 0} (-1)^n \vert \Lambda_n(\calx/G) \vert
\\ & = & \sum_{n \geq 0} (-1)^n c_n(\calx/G)
\\ & = & \chi(\calx/G;R).
\end{eqnarray*}
For the $L^2$-Euler characteristic on the other hand, we have
\begin{eqnarray*}
\chi^{(2)}(\hocolim_{\calx/G} F) &=& \sum_{n \geq 0} (-1)^n \cdot \left( \sum_{\overline{a} \in \Lambda_n(\calx/G)} \chi^{(2)}(F(s\overline{a}_1)) \right) \\
& = & \sum_{n \geq 0} (-1)^n \cdot \left( \sum_{\overline{a} \in \Lambda_n(\calx/G)} \frac{1}{\vert G_{sa_1} \vert} \right) \\
& = & \sum_{n \geq 0} (-1)^n \cdot \left( \sum_{\overline{a} \in \Lambda_n(\calx)/G} \frac{1}{\vert G_{a} \vert} \right) \\ & = & \sum_{n \geq 0} (-1)^n \cdot \left( \sum_{\overline{a} \in \Lambda_n(\calx)/G} \frac{\vert \text{orbit}(a) \vert}{\vert G \vert} \right)  \\ & = & \frac{1}{\vert G \vert} \sum_{n \geq 0} (-1)^n \cdot \left( \sum_{\overline{a} \in \Lambda_n(\calx)/G} \vert \text{orbit}(a) \vert \right) \\ & = & \frac{1}{\vert G \vert} \sum_{n \geq 0} (-1)^n \vert \Lambda_n(\calx) \vert
\\ & = & \frac{1}{\vert G \vert} \sum_{n \geq 0} (-1)^n c_n(\calx)
\\ & = & \frac{\chi^{(2)}(\calx)}{\vert G \vert}.
\end{eqnarray*}
\end{proof}

\begin{example} \label{exa:necessary_conditions_for_developability:single_arrow}
By the classical theorem of Bass--Serre, any injective group
homomorphism
\begin{equation} \label{equ:one_arrow_complex_of groups}
G_0 \to G_1
\end{equation}
is a developable complex of groups.  The $L^2$-Euler characteristic
of the homotopy colimit of \eqref{equ:one_arrow_complex_of groups}
is $1/\vert G_1 \vert$ by
Example~\ref{exa:hocolim_formula_for_I_with_terminal_object}.
Theorem~\ref{thm:Euler_characteristic_of_hocolim_of_quotient_complex}
then says we must have
\[
\frac{\vert G \vert}{\vert G_1 \vert} = \chi^{(2)}(\calx)
=\chi(B\calx;\IC)
\]
if \eqref{equ:one_arrow_complex_of groups} is to be developable from
a scwol $\calx$ by an action of $G$ in the sense of
Definition~\ref{def:group_action_on_a_scwol}. Thus
\eqref{equ:one_arrow_complex_of groups} is not developable from any
scwol $\calx$ whose geometric realization has Euler characteristic
0, such as $\{j \rightrightarrows k\}$. Nor can
\eqref{equ:one_arrow_complex_of groups} be developed from any scwol
$\calx$ with $\chi(B\calx;\IC)$ negative. The integer $\vert G \vert
$ must also be divisible by $\vert G_1 \vert$, since
$\chi(B\calx;\IC)$ is always an integer.  Moreover, the Euler
characteristic of $\calx$ must be less than or equal to $\vert G
\vert$. This trivial example illustrates how one can find necessary
conditions on $\calx$ and $G$ if a given complex of groups is to be
developable from $\calx$ and $G$.
\end{example}

\begin{example}[Euler characteristics of transport groupoid in finite
case] \label{exa:transport_groupoid_finite_case} Let $X$ be a finite
set and $G$ a finite group acting on $X$. Let $R$ be a ring such that the
orders of subgroups of $G$ are invertible in $R$.
Considering $X$ as a
scwol, we clearly have an action in the sense of
Definition~\ref{def:group_action_on_a_scwol}. The associated complex
of groups $F:X/G \to \GROUPS$ assigns to $\text{orbit}(\sigma)$ the
stabilizer $G_\sigma$. The homotopy colimit $\hocolim_{X/G} F$ is
equivalent to the transport groupoid $\calg^G(X)$ of
Example~\ref{exa:transport_groupoid}, so
\[
\chi\left(\calg^G(X);R\right)=\chi(\hocolim_{X/G} F;R)=\chi(X/G;R)=\vert
X/G \vert.
\]
For the $L^2$-Euler characteristic, on the other hand, we have
\[
\chi^{(2)}\left(\calg^G(X)\right)=\chi^{(2)}(\hocolim_{X/G}
F)=\frac{\chi^{(2)}(X)}{\vert G \vert}=\frac{\vert X\vert}{\vert G
\vert},
\]
a formula obtained by Baez--Dolan~\cite{Baez-Dolan(2001)}.
\end{example}

We also generalize the following formula of Haefliger for the Euler
characteristic of the homotopy colimit of a (not necessarily
developable) complex of groups.

\begin{theorem}[Corollary 3.5.3 of Haefliger~\cite{Haefliger(1992)}] \label{the:Haefligers_corollary}
Let $G(X)$ be a complex of groups over a finite ordered simplicial
cell complex $X$. Assume that each $G_\sigma$ is the fundamental
group of a finite aspherical cell complex. Then $BG(X)$ has the
homotopy type of a finite complex and its Euler-Poincar\'e
characteristic is given by\footnote{Haefliger's original formula
has, instead of the lower link $L^\sigma$, the upper link
$L_\sigma$, which is the full subcategory of the scwol $C(X)
\downarrow \sigma$ on all objects except $1_\sigma$. However, this
is merely a typo, for if we use the upper link $Lk_\sigma$ and
consider the example $C(X)=\{k \leftarrow j \to \ell\}$ with pseudo
functor $G(X)(\ell):=\IZ$ and $G(X)(j):=G(X)(k):=\{0\}$, then
$\chi(BG(X))=\chi(S^1)=0$ but
$\sum(1-\chi(Lk_\sigma))\chi(G_\sigma)=1$.}
\[
\chi(BG(X))=\sum_{\sigma \in \ob(C(X))}
(1-\chi(Lk^\sigma))\chi(G_\sigma).
\]
\end{theorem}

The terms in Haefliger's theorem have the following meanings. An
{\em ordered simplicial cell complex} $X$ is by definition the nerve
of a skeletal scwol, denoted $C(X)$. The notation $BG(X)$ signifies the
geometric realization of the nerve of the homotopy colimit of the pseudo functor
$G(X)\colon C(X) \to \GROUPS$. An {\it aspherical} cell complex is
one for which all homotopy groups beyond the fundamental group
vanish. The {\it lower link $Lk^\sigma$ of the object $\sigma$} is
the full subcategory of the scwol $\sigma \downarrow C(X)$ on all
objects except $1_\sigma$.

\begin{theorem}[Extension of Corollary 3.5.3 of Haefliger~\cite{Haefliger(1992)}] \label{the:extension_of_Haefligers_corollary}
Let $\cali$ be a finite skeletal scwol and $F\colon \cali \to
\GROUPS$ a complex of groups such that for each object $i$ of
$\cali$, the group $F(i)$ is of type (FF$_\IZ$). Then
\[
\chi(B \hocolim_\cali F )=\sum_{i \in \ob(\cali)}
(1-\chi(BLk^i))\chi(B F(i)),
\]
where $B$ indicates geometric realization composed with the nerve
functor.
\end{theorem}

\begin{proof}
All hypotheses of
Theorem~\ref{the:homotopy_colimit_formula}\ref{the:homotopy_colimit_formula:chi}
are satisfied. The skeletal scwol $\cali$ is directly finite by
Lemma~\ref{lem:scwol_directly_finite_EI} and admits a finite
$\cali$-$CW$-model for its $\cali$-classifying space by
Theorem~\ref{the:finite_models_for_finite_scwols}. Each group
$\calc(i)$ is automatically directly finite, and assumed to be of
type (FF$_\IZ$). The bar construction model $E^{\barcon}\cali$ in
Remark~\ref{rem:Ebarcalc} has an $n$-cell based at $i$ for each path
of $n$-many non-identity morphisms in $\cali$
\begin{equation*}
i \to i_1 \to i_2 \to \cdots \to i_n.
\end{equation*}
Each such path in $\cali$ corresponds uniquely to a path of
$(n-1)$-many non-identity morphisms in the scwol $Lk^{i}$ beginning
at the object $i \to i_1$. Thus
\begin{eqnarray*}
1-\chi(BLk^{i}) & = & 1- \sum_{m \geq 0} (-1)^m c_m(Lk^{i}) \\
& = & 1- \sum_{m \geq 0} (-1)^m \text{card}\{\text{$(m+1)$-paths in
$\cali$ beginning at $i$}\} \\
& = & 1- \sum_{n \geq 1} (-1)^{n-1} \text{card}\{\text{$n$-paths in
$\cali$ beginning at $i$}\} \\
& = & 1+ \sum_{n \geq 1} (-1)^{n} \text{card}\{\text{$n$-paths in
$\cali$ beginning at $i$}\} \\
& = & \sum_{n \geq 0} (-1)^{n} \text{card}\{\text{$n$-paths in
$\cali$ beginning at $i$}\}. \\
\end{eqnarray*}
Then by
Theorem~\ref{the:homotopy_colimit_formula}~\ref{the:homotopy_colimit_formula:directly_finite},
Theorem~\ref{the:homotopy_colimit_formula}~\ref{the:homotopy_colimit_formula:(FF)},
Theorem~\ref{the:coincidence}, and
Theorem~\ref{the:homotopy_colimit_formula}~\ref{the:homotopy_colimit_formula:chi},
we have
\begin{eqnarray*}
\chi(B \hocolim_\cali F ) & = & \chi(\hocolim_\cali F ) \\
& = & \sum_{n \geq 0} (-1)^n \cdot \sum_{\lambda \in \Lambda_n}
\chi(F(i_\lambda)) \\
& = & \sum_{i \in \ob(\cali)}\left( 1-\chi(BLk^{i}) \right) \cdot
\chi(F(i)) \\
& = & \sum_{i \in \ob(\cali)}\left( 1-\chi(BLk^{i}) \right) \cdot
\chi(BF(i)).
\end{eqnarray*}
\end{proof}

\begin{remark}
The assumptions in our Theorem~\ref{the:extension_of_Haefligers_corollary} on the groups $F(i)$ are related to
the assumptions in Theorem~\ref{the:Haefligers_corollary} on the groups $G_\sigma$ in that any finitely presentable group of type (FF$_\IZ$) admits a finite model for its classifying space.
\end{remark}


\typeout{-----------------------------Appendix-------------------------------------------}

\section{Appendix} \label{sec:appendix}

Let $G$ be a group acting on a scwol $\calx$ in the sense of
Definition~\ref{def:group_action_on_a_scwol}. In connection with
Theorem~\ref{the:reduction_to_skeletal_case}, we remark that it is
possible to choose a skeleton $\Gamma_0$ of $\calx$, a
$G$-equivariant functor $r\colon \calx \to \Gamma_0$, and a natural
isomorphism $\eta\colon ir \cong \id_\calx$ so that
\begin{itemize}
\item
the inclusion $i_0\colon \Gamma_0 \to \calx$ is $G$-equivariant,
\item
$ri_0=\id_{\Gamma_0}$, and
\item
for every object $x \in \ob(\calx)$ and each $g \in G$, we have
$\eta_{gx}=g\eta_x$.
\end{itemize}

To prove this, we first choose the object set of $\Gamma_0$ via an
equivariant section of the projection $\pi\colon \ob(\calx) \to
\iso(\calx)$, which assigns to each object of $\calx$ its
isomorphism class of objects. Let $\Theta$ denote the set of
$G$-orbits of $\iso(\calx)$. For each $G$-orbit $\theta \in \Theta$,
we use the axiom of choice to select an element $\overline{x}_\theta
\in \theta$. For each $\theta$, select then a $\pi$-preimage
$s(\overline{x}_\theta):=x_\theta$ of $\overline{x}_\theta$.  On the
orbit of each $\overline{x}_\theta$  we define the section $s$ by
$s(g\overline{x}_\theta):=gx_\theta$. This is well defined, for if
$g_1\overline{x}_\theta=g_2\overline{x}_\theta$, then $g_1x_\theta
\cong g_2x_\theta$, and $g_1x_\theta = g_2x_\theta$ by
Lemma~\ref{lem:consequences_of_group_action}~\ref{lem:consequences_of_group_action:(i)}.
Define the skeleton $\Gamma_0$ to be the full subcategory of $\calx$
on the objects in the image of the equivariant section $s\colon
\iso(\calx) \to \ob(\calx)$.

For each $\overline{x}_\theta$, and each $x \in
\overline{x}_\theta$, choose an isomorphism $\eta_x \colon x_\theta
\to x$. For $gx$, we define $\eta_{gx}$ as $g\eta_x$. Next, we
define a functor $r\colon \calx \to \Gamma_0$ on objects $x \in
\ob(\calx)$ by $r(x):=s\pi(x)$ and on morphisms $f\colon x \to y$ by
$r(f):= \eta_y \circ f \circ \eta_{x}^{-1}$. Then $\eta$ is clearly
a natural isomorphism, the inclusion $i_0\colon \Gamma_0 \to \calx$
is $G$-equivariant, and $ri_0=\id_{\Gamma_0}$.



\bibliographystyle{abbrv}

\end{document}